\UseRawInputEncoding
\documentclass[11pt]{amsart}
\usepackage{mathrsfs, amsmath, amssymb, array,graphicx, booktabs,diagbox,listings,color,textcomp,picinpar,wrapfig,cutwin,longtable,slashbox,bm}

\usepackage{booktabs}
\usepackage{multirow}
\usepackage{makecell}
\usepackage{hhline}
\usepackage[table]{xcolor}

\usepackage[page,title,titletoc,header]{appendix}
\usepackage[font=small,labelfont=bf]{caption}
\DeclareCaptionFont{tiny}{\tiny}
\DeclareCaptionFont{normalsize}{\normalsize}

\usepackage{geometry}
\definecolor{listinggray}{gray}{0.9}

\definecolor{lbcolor}{rgb}{0.9,0.9,0.9}
\definecolor{myblue1}{RGB}{202,206,250}
\definecolor{mygrey}{RGB}{191,191,191}
\definecolor{darkgray}{gray}{0.6}

\usepackage{floatrow}
\newfloatcommand{capbtabbox}{table}[][\FBwidth]
\newcolumntype{P}[1]{>{\centering\arraybackslash}p{#1}}

\usepackage{algorithm}
\usepackage[noend]{algpseudocode}

\def\baselinestretch{1.2}

\newtheorem{theorem}{Theorem}[section]
\newtheorem{lemma}[theorem]{Lemma}
\newtheorem{corollary}[theorem]{Corollary}
\newtheorem{proposition}[theorem]{Proposition}

\theoremstyle{definition}
\newtheorem{definition}[theorem]{Definition}

\theoremstyle{remark}
\newtheorem{remark}[theorem]{Remark}
\numberwithin{equation}{section}

\newcommand{\customizedfootnotetext}[2]{
	\let\oldthefootnote\thefootnote
	\renewcommand{\thefootnote}{#1}
	\footnotetext{#2}
	\let\thefootnote\oldthefootnote
}

\begin{document}

\title{Orientable hyperbolic 4-manifolds over the 120-cell}

\author{Jiming Ma}
\address{School of Mathematical Sciences \\Fudan University\\Shanghai 200433, China} 
\email{majiming@fudan.edu.cn}

\author{Fangting Zheng}
\address{Department of Mathematical Sciences, Xi'an Jiaotong-Liverpool University, Suzhou 215123, China  }
\email{Fangting.Zheng@xjtlu.edu.cn}



\date{Jun. 8, 2022}

\thanks{Jiming Ma was supported by  NSFC  12171092. Fangting Zheng was supported by NSFC 12101504.}

\begin{abstract}
Since there is no hyperbolic Dehn filling theorem for higher dimensions, it is challenging to construct explicit hyperbolic manifolds of small volume in dimension at least four. Here, we build up closed  hyperbolic 4-manifolds of volume $\frac{34\pi^2}{3}\cdot 16$ by using the small cover theory. In particular, we classify all of the orientable four-dimensional small covers over the right-angled 120-cell up to homeomorphism; these are all with even intersection forms.
\end{abstract}

\maketitle


\section{Introduction}

\subsection{\textbf{Hyperbolic manifolds}}
By the Thurston-Perelman's geometrization theorem, the majority of $3$-manifolds are known to be hyperbolic $3$-manifolds. Moreover,  according to the Jorgensen-Thurston theory, it is straightforward to construct closed hyperbolic $3$-manifolds $M(\epsilon)$ by generically Dehn filling cusped hyperbolic $3$-manifolds $M$ along the tori. However, the Jorgensen-Thurston theory does not hold in $\mathbb{H}^n$ when $n$ is greater than or equal to four. More precisely, when conducting a sufficiently large Dehn filling on cusp ends of a higher dimensional cusped hyperbolic manifold, we can obtain a manifold that admits a complete, finite volume Einstein metric; in this case, however, the fundamental group is relatively hyperbolic rather than hyperbolic in general \cite{Anderson:2006}. Hence, it is necessary to develop some methods to describe the higher dimensional situations. For example, Martelli and Riolo have studied the four-dimensional hyperbolic Dehn filling to another extent \cite{martelliRiolo:2017}. 

Moreover, besides the well-known four color theorem, it is helpful to understand the manifold theory through programming, for example see \cite{conder_dobcsanyi:2005,conder_maclachlan:2005,Everitt_Ratcliffe_Tschantz:2011,gabai_meryerhoff_milley:2009,gabai_meyerhoff_thurston:2003,jtcube,Linowitz_Voight:2015,long:2008,marshall_martin:2012,NV:15}. For the hyperbolic case, few broadly applicable techniques are currently available; a limited number of explicit examples have been reported in the literature. For example, Ratcliffe and Tschantz have determined the structure of the congruence two subgroups  $\Gamma_2^4$ of the group $\Gamma^4$ of integral, positive, Lorentzian $5\times 5$ matrices. They have classified all of the corresponding hyperbolic four-manifolds that have a minimum volume \cite{rt:2000}.  Later, Saratchandran used the Ratcliffe-Tschantz census to study the topology of cusped hyperbolic $4$-manifolds \cite{Sara1,Sara2, Sara3,Sara4}. 

Theoretically, by initially constructing a hyperbolic orbifold and then applying the Selberg lemma, a hyperbolic manifold can be constructed. However, this approach is often inefficient since there is no universal method to determine the covering degree. For instance, given a Fuchsian group $\Gamma$,  the minima of indices of torsion-free subgroups of $\Gamma$ is bounded above by the twice of the least common multiple (LCM) of orders of the finite subgroups of $\Gamma$ \cite{EEK:1982}. However, Jones and Reid have shown that the complex version can vary substantially. More precisely, the authors reported a sequence $\Gamma_k$ of co-compact Kleinian groups for which the ratio of the minimum index to the LCM could be arbitrarily large \cite{JR:1998}.

\subsection{\textbf{Small cover}}
Small covers have been studied by Davis and Januszkiewicz in \cite{dj:1991}, see also \cite{Vesnin:1987,Vesnin:1991,Izmestiev:01}.  
A small cover of an $n$-dimensional simple polytope $P$ is a manifold which arises as an orbifold cover of $P$ of minimal index $2^n$.  In other words, an $n$-dimensional  small covers admits  a locally standard $\mathbb{Z}_2^{n}$-action and the orbit space is an  $n$-dimensional simple polyhedron.  It turns out that any small cover of $P$ is determined by a map of the (abstract) right-angled Coxeter group $\Gamma$ associated to $P$ into the group $\mathbb{Z}_2^n$. Such map is called a \emph{characteristic function} or a \emph{coloring}. For an (abstract) right-angled Coxeter group $\Gamma$ associated to some three-dimensional simple polytope  $P$, $\Gamma$ is guaranteed to have a torsion-free subgroup of index eight according to the four color theorem. In other words, we can always construct a small cover over a three-dimensional simple polytope. In higher dimensions, however, this problem presents significant challenges and has driven the study of the \emph{Buchstaber invariant} \cite{BP}.

The algebraic and topological invariants of a small cover can be calculated by the combinatorics of the orbit polytope  and the coloring.  For example, the $mod~2$ Betti number $\beta_{i}^{(2)}$ of a small cover $M$ over the polytope $P$ agrees with $h_{i}$, where $h=(h_{0}, h_{1},\ldots, h_{n})$ is the $h$-vector of the polytope $P$ \cite{dj:1991}.  According to some recent results of the homology groups and the cohomology ring of small covers \cite{2licai:2013,DLiu:2018,mz1}, a small cover of dimension three or four has at most $2$-torsion in its homology group. Moreover, every oriented four-dimensional small cover has  zero signature and even intersection form.  Therefore the isomorphism class of the intersection form is determined by the second Betti numbers. See Section \ref{section:smallcover} and \ref{section:classify} for additional details.

Suppose that the simple polytope $P$ can be realized in hyperbolic $n$-space as a right-angled polytope. Then any small cover over $P$ is a hyperbolic
manifold.
For example, the first closed hyperbolic 3-manifold that appeared in the literature,
the L\"obell manifold, is a small cover \cite{Vesnin:1987},  see also \cite{Vesnin:17}.
Ratcliffe and Tschantz classified the congruence-2 hyperbolic 4-manifolds of minimal
volume. They provided a list of the resulting 1171 non-compact manifolds, of
which $23$ are orientable, together with their homology groups. There is a
straightforward generalization of small covers for the non-compact case. The problem studied by Ratcliffe and Tschantz in \cite{rt:2000} can be reformulated as the classification of the small covers of a non-compact right-angled hyperbolic 4-polytope discovered by Potyagailo and Vinberg \cite{pv:2005}. 
The analogous problem in dimension $5$ has then been solved similarly by Ratcliffe
and Tschantz in \cite{rt:2004}.
In this paper, we use the small cover theory to construct a family of closed hyperbolic 4-manifolds, which are obtained  after 16-degree coverings.

In comparison to the results of Ratcliffe and Tschantz, we consider the  minimal index torsion free subgroups of the Coxeter group associated to the right-angled 120-cell. The Coxeter  group used here has many more generators and relations than the group $\Gamma^4$ of integral, positive, Lorentzian $5\times 5$ matrices; this presents additional challenges in the identification of minimal index torsion free subgroups. For other hyperbolic 4-manifolds of small volume. See \cite{davis:1985,rt:2001,conder_maclachlan:2005,long:2008}.

Garrison and Scott have considered the simplest $\mathbb{Z}_2^4$-coloring case over the right-angled 120-cell \cite{Scott:02}, where the image of the characteristic functions is $\{e_1,~e_2,~e_3,~e_4,~e_1+e_2+e_3+e_4\}$. They have found that there is only one, up to homeomorphism, such non-orientable manifold. The manifold and the corresponding characteristic function are denoted as \emph{Garrison-Scott manifold} and \emph{Garrison-Scott coloring}, respectively. We obtain some other hyperbolic manifolds from the right-angled 120-cell and prove the following:

 \begin{theorem}\label{theorem:orientable}
There are exactly 56 orientable small covers over the right-angled hyperbolic 120-cell up to homeomorphism.
\end{theorem}
\begin{theorem}\label{theorem:nonorientable eight-coloring}
There are exactly two 8-coloring  non-orientable  small covers over the right-angled hyperbolic 120-cell up to homeomorphism.
\end{theorem}

In addition, one 9-coloring manifold is successfully constructed. See Section \ref{section:9col_cr} and \ref{section:9-coloring} for details. Therefore, we obtain a total of 59 new closed hyperbolic $4$-manifolds with volume $\frac{544\pi^2}{3}$ through two algorithms. See Table \ref{table:ori}, \ref{table:nonori} and \ref{table:9col} for details. The \emph{color-recursion algorithm} is much more universal, but is too inefficient to solve our problems. The \emph{block-pasting algorithm} is far more effective as it takes advantage of the symmetry information of the 120-cell; however, it still exhibits an unaffordable running time when considering $k$-colorings, where $k$ is greater than or equal to nine. We have identified all of the orientable small covers over the 120-cell, which are realized by at most eight colors. See Section 6 for more details. 

\subsection{Intersection form}
The Davis manifold \cite{davis:1985} is the first example of a hyperbolic $4$-manifold with  a calculated even type intersection form  $\oplus_{36}\begin{pmatrix}
 0 & 1\\
 1&0\\
 \end{pmatrix}$ \cite{rt:2001}. Quite recently, Martelli-Riolo-Slavich exhibited the first concrete closed oriented hyperbolic four-manifold with odd intersection form \cite{MRS:2019}. Here, we construct $60$ more examples with even intersection forms.
\begin{theorem}  \label{theorem:manifold}
	For the 60 manifolds $M_1,\cdots,M_{56},M_{57},M_{58},M_{59},M_{60}$ claimed in Tables \ref{table:ori}--\ref{table:9col}, we have their homologies and intersection forms as below:
{\renewcommand\baselinestretch{1.33}\selectfont	
	\begin{table}[H]
		\small
		\begin{tabular}{|c|c|c|c|c|}
			\hline
			$H_i(M;R)$& $M=M_j, 1\le j\le 56$ & $M=M_j, 57\le j\le 60$ & $M=\widetilde{M_j}, j=59$& $M=\widetilde{M_j}, j=57,58,60$\\
			\cline{2-5}
			&\multicolumn{2}{c}{$R=\mathbb{Z}$}&\multicolumn{2}{|c|}{$R=\mathbb{Q}$}\\
			\hline	
			$i=0$ & $\mathbb{Z}$ & $\mathbb{Z}$& $\mathbb{Q}$ & $\mathbb{Q}$\\
			\hline	
			$i=1$ & $\mathbb{Z}_2^{116}$& $\mathbb{Z}_2^{116}$ &  $\mathbb{Q}^{115}$ &$\mathbb{Q}^{51}$\\
			\hline
			$i=2$ & $\mathbb{Z}^{134}\oplus\mathbb{Z}_2^{116}$& $\mathbb{Z}^c\oplus\mathbb{Z}_2^{(250-c)}$&   $\mathbb{Q}^{500}$ &$\mathbb{Q}^{372}$  \\
			\hline
			$i=3$ & trivial   & $\mathbb{Z}^{(c-135)}\oplus\mathbb{Z}_2$&   \ $\mathbb{Q}^{115}$ &$\mathbb{Q}^{51}$  \\
			\hline
			$i=4$ & $\mathbb{Z}$ & $\mathbb{Z}_2$ &$\mathbb{Q}$  & $\mathbb{Q}$\\
			\hline
			\hline
			$Q(M)$ & {$\oplus_{67}\begin{pmatrix}
			0 & 1\\
			1&0\\
			\end{pmatrix}$ } &\diagbox[dir=SW]{\kern+4em}{\kern+4em}& {$\oplus_{250}\begin{pmatrix}
			0 & 1\\
			1&0\\
			\end{pmatrix}$ }& {$\oplus_{186}\begin{pmatrix}
			0 & 1\\
			1&0\\
			\end{pmatrix}$ }\\
			\hline
		\end{tabular}
		\caption{Homologies and intersection forms of the constructed manifolds}
		\label{table:topology}
	\end{table}
}
\noindent where, $c$ is the sum of the first reduced Betti numbers of the full-subcomplex $K_\omega$ for $\omega \in row~ \Lambda$. The $Q(M)$ means the intersection form of the orientable $4$-manifold $M$. Among these sixty small covers over the 120-cell $\mathcal{E}$, the manifolds $M_1,M_2,\cdots,M_{56}$ are orientable, whereas $M_{57},M_{58},M_{59},M_{60}$  are non-orientable. The manifold $\widetilde{M_j}$ is the orientable double cover of $M_j$, for $57\le j\le 60$.
\end{theorem}

\textbf{Outline of the paper}: In Section 2 and Section 3, we present some preliminaries on the algebraic theory  of small covers and the combinatorics of the 120-cell, respectively. In Section 4, we use the color-recursion algorithm, which is broadly applicable but of low efficiency, to address the  case where the total number of colors assigned to the facets of the right-angled 120-cell to construct small covers is $4,5,6~ \text{or}~ 7$. Moreover, one new small cover constructed by 9 colors is obtained by this approach. Based on the high symmetry of the 120-cell, we then introduce a substantially more efficient block-pasting algorithm in Section 5 and study the more complicated cases, where the total number of colors  assigned to the facets of the right-angled 120-cell to construct small covers is $8~ \text{or}~ 9$. We classify all of the small covers over the right-angled 120-cell that are obtained in previous sections by homeomorphism in Section 6, thereby proving Theorem  \ref{theorem:orientable} and Theorem \ref{theorem:nonorientable eight-coloring}. Finally, we consider the intersection forms of all of these manifolds and complete the proof of Theorem \ref{theorem:manifold} in Section 7.

\textbf{Acknowledgments}: The authors appreciate greatly the referees' valuable and constructive  comments on an earlier draft of this paper, which make the paper much more precise and readable. We are also grateful to Leonardo Ferrari who kindly pointed out some typos in figures of adjacent matrices.

\subsection{\textbf{Small covers}}\label{section:smallcover}
Given an $n$-dimensional simple polytope $P$, the manifold corresponding to a $2^{n}$-index torsion-free subgroup $\Gamma$ of the (abstract) right-angled Coxeter group associated to  $P$ is a \emph{small cover over $P$}  \cite{dj:1991}. There is another equivalent but more practical way to describe a small cover by using the language of \emph{coloring} \cite{dj:1991}: 

Let $\mathcal{F}(P)=\{F_1,F_2,\ldots,F_m\}$ be  the set of co-dimensional one faces of $P$.  Let us define the \emph{$\mathbb{Z}_2^r$-coloring characteristic function}, $n \leq r\leq m$, to be a function $$\lambda:\mathcal F(P)=\{F_1,F_2,\cdots,F_m\}\longrightarrow \mathbb{Z}_2^r,$$ that satisfies the \emph{non-singularity condition}. That is $\lambda(F_{i_1}),\lambda(F_{i_2}),\cdots,\lambda(F_{i_n})$ should generate a  subgroup  of $\mathbb{Z}_2^r$, which is isomorphic to $\mathbb{Z}_2^n$, when the facets $F_{i_1}, F_{i_2}, \cdots, F_{i_n}$ share a common vertex.  Note that the function $\lambda$ is not guaranteed to exist. The corresponding \emph{characteristic matrix} is obtained by placing the vector $\lambda(F_i)$ column by column. For the sake of brevity, a $\mathbb{Z}_2^r$-coloring characteristic function $\lambda$ is referred to as \emph{a coloring}.  Every vector in $\mathbb{Z}_2^r$ is called a \emph{color}. The vector $\lambda(F_i)\in \mathbb{Z}_2^r$ is called a \emph{color} of $F_i$. In particular, if $r=4$, there is a total of $16$ colors. Say a coloring is a \emph{$k$-coloring}  if exactly $k$ different colors are used for its construction. These $k$ colors form a \emph{coloring set} $\{c_1,c_2,\cdots,c_k\}$ of $\lambda$.
 
If the characteristic function $\lambda$ is defined successfully, we can recover a smooth $n$-manifold $M(P,\lambda):=P\times \mathbb{Z}_2^r/\sim $, which is named \emph{real toric manifold over $P$}, by the following equivalence relation:

\begin{center}
	$(x,g_1)\sim (y,g_2){\Longleftrightarrow }
	\left\{
	\begin{array}{rcl}x=y \hspace{0.45cm} {\rm and}  \hspace{0.5cm} g_1=g_2  \hspace{0.7cm} &\mbox~~~~~{\rm if}& x\in\ int P, \\
	x=y ~~ {\rm and} ~~~ g_1^{-1}g_2 \in G_f&\mbox ~~~~~{\rm if}& x \in \partial P,
	\end{array}\right.
	$
\end{center}
where, $f=F_{i_1}\cap\cdots\cap F_{i_{n-l}}$ is the unique face of co-dimensional $n-l$ that contains $x$ as an interior point, and $G_f$ is a subgroup generated by $\lambda(F_{i_1}), \lambda(F_{i_2}), \cdots, \lambda(F_{i_{n-l}})$. For brevity, we  only refer to the colorings when the polytope is given instead of talking about both colorings and manifolds. If $r=n$, then the corresponding manifold is called a \emph{small cover}.  

\begin{figure}[H]
	\scalebox{0.58}[0.58]{\includegraphics {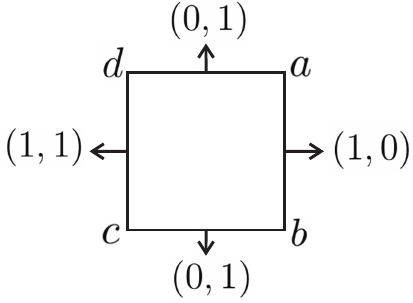}}
	\caption{A coloring over the square}\label{figure: klein_example}
\end{figure}

For example, define a $\mathbb{Z}_2^2$-coloring characteristic function $\lambda$ on the 2-dimensional cube $I^2$ as shown in Figure \ref{figure: klein_example}, where $(1,0)=e_1$ and $(0,1)=e_2$ are the standard basis of $\mathbb{Z}_{2}^{2}$. The small cover $M(I^2,\lambda)$ is the Klein bottle.
\noindent 

Since the dual of the boundary of a simple polytope $P$ is a simplicial complex $K$ \cite{BP}, the definition of real toric manifolds introduced above has a dual version. By substituting  the facet set $\mathcal{F}(P)$ with the vertex set $\mathcal{V}$ of the simplicial complex $K$, we can also define the \emph{characteristic function} $\lambda$ on $K$, namely $$\lambda:\mathcal V(K)=\{v_1,v_2,\cdots,v_m\}\longrightarrow \mathbb{Z}_2^r$$  The non-singularity condition changes as follows: if the convex hull $conv\{v_{i_1}, v_{i_2}, \cdots, v_{i_n}\}$ of $n$ vertices $v_{i_1},v_{i_2},\cdots,v_{i_n}$  is a facet of $K$, the images $\lambda(v_{i_1}),\lambda(v_{i_2}),\dots,\lambda(v_{i_n})$ shall generate a  subgroup isomorphic to $\mathbb{Z}_2^n$. 

To keep a concise notation, we encode every $\mathbb{Z}_2^*$-color by an integer through the binary-to-decimal conversion. For example, in the $\mathbb{Z}_2^3$-coloring case, we use the numbers $1,2,\cdots ,7$ to represent the following seven colors $(1,0,0)$, $(0,1,0)$,  $(1,1,0)$, $(0,0,1)$ $(1,0,1)$, $(0,1,1)$, and  $(1,1,1)$, respectively. Thus, a characteristic matrix can also be viewed as a \emph{characteristic vector}.  For example, the characteristic vector $C$ of characteristic matrix
\begin{center}
	$(\lambda(\mathcal{F}_1),\lambda(\mathcal{F}_2),\lambda(\mathcal{F}_3),\lambda(\mathcal{F}_4),\lambda(\mathcal{F}_5))=
	\begin{pmatrix}
	1 & 0 & 0 & 0 & 1\\
	0 & 1 & 0 & 0 & 1\\
	0 & 0 & 1 & 0 & 1\\
	0 & 0 & 0 & 1 & 1\\
	\end{pmatrix}$
\end{center}
is $(1, 2, 4, 8, 15)$. These integers are also called \emph{colors} if there is no ambiguity. In the previous example, we may say the color of facet $F_5$ is 15 or $(1,1,1,1)$. The characteristic function $\lambda$, the characteristic matrix $\Lambda$, and the characteristic vector $C$ can be inferred from each other easily. The characteristic vector $C$ is the most  concise representational form among them.

\subsection{$GL_n(\mathbb{Z}_2)$-equivalence and orientability}

\begin{figure}[H]
	\scalebox{0.45}[0.45]{\includegraphics {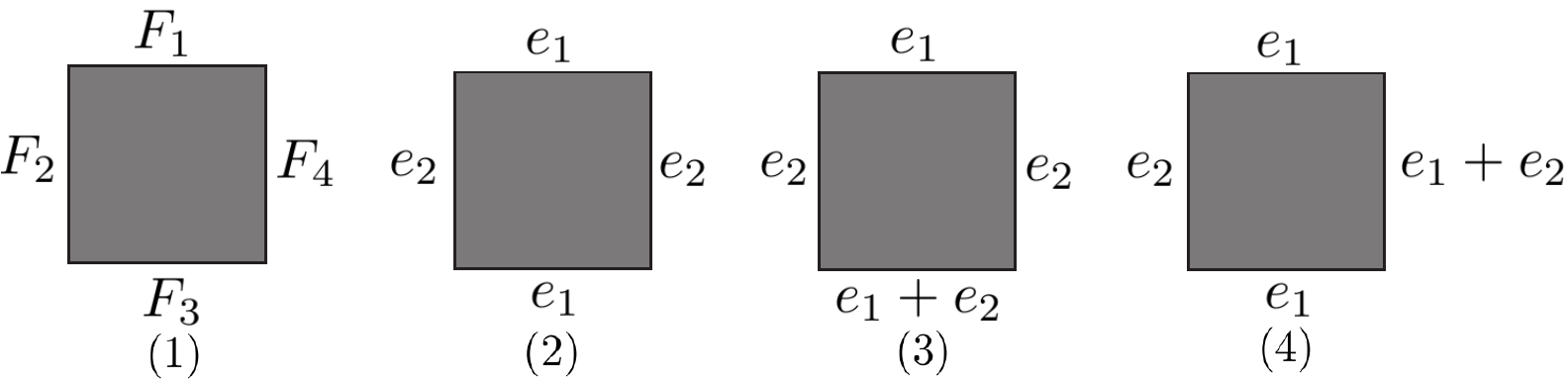}}
	\caption{Three $GL_2(\mathbb{Z}_2)$-equivalence representatives.}\label{figure:gl_cube}
\end{figure}

Let $P$ be an $n$-dimensional simple polytope with its first facets $F_1,F_2,\cdots,F_n$ sharing a vertex. Suppose $\lambda$ is a $\mathbb{Z}_2^r$-coloring characteristic function over $P$ that satisfies $\lambda(F_i)=e_i$, for $1\leq i\leq n$, where $\{e_i\}_{i=1}^n$ is the standard basis of $\mathbb{Z}_2^n$. Colors on other facets are linear combinations of the colors on the first $n$ facets. Another coloring $\widetilde{\lambda}$ is \emph{$GL_n(\mathbb{Z}_2)$-equivalent} to $\lambda$ if there exists some $g\in GL_n(\mathbb{Z}_2)$ such that $\widetilde{\lambda}=g\circ\lambda$. This relation can be easily checked to determine if it is an equivalence relation. Thus, all the characteristic functions over a given simple polytope can be classified under this equivalence relation. Then, $\lambda$ is referred to as the \emph{$GL_n(\mathbb{Z}_2)$-equivalence representative} in its equivalence class.  For example, assume $P$ is the 2-dimensional cube $I^2$ and its facets are ordered as shown in Figure \ref{figure:gl_cube} (1). There are a total of 18 small covers over $P$ as shown in Table \ref{table:square}, which can be grouped into three $GL_2(\mathbb{Z}_2)$-equivalence classes in total. The \emph{$GL_2(\mathbb{Z}_2)$-equivalence representatives} of these three classes are ($e_1$, $e_2$, $e_1$, $e_2$), ($e_1$, $e_2$, $e_1+e_2$, $e_2$) and ($e_1$, $e_2$, $e_1$, $e_1+e_2$), respectively, as shown in Figure \ref{figure:gl_cube}.
\begin{table}[H]
	\begin{tabular}{|c|c|c|}
		\hline
		($e_1$, $e_2$, $e_1$, $e_2$)&($e_1$, $e_2$, $e_1+e_2$, $e_2$)&($e_1$, $e_2$, $e_1$, $e_1+e_2$)\\
		
		($e_2$, $e_1$, $e_2$, $e_1$)&($e_2$, $e_1$, $e_1+e_2$, $e_1$)&($e_2$, $e_1$ $e_2$, $e_1+e_2$)\\
		
		($e_1+e_2$, $e_1$, $e_1+e_2$, $e_1$)&($e_1+e_2$, $e_1$, $e_2$, $e_1$)&($e_1+e_2$, $e_1$, $e_1+e_2$, $e_2$) \\

		($e_1$, $e_1+e_2$, $e_1$, $e_1+e_2$)&($e_1+e_2$, $e_2$, $e_1$, $e_2$)&($e_1$, $e_1+e_2$, $e_1$, $e_2$)\\

		($e_2$, $e_1+e_2$, $e_2$, $e_1+e_2$)&($e_2$, $e_1+e_2$, $e_1$, $e_1+e_2$)&($e_2$, $e_1+e_2$, $e_2$, $e_1$)\\

		($e_1+e_2$, $e_2$, $e_1+e_2$, $e_2$)&($e_1$, $e_1+e_2$, $e_2$, $e_1+e_2$)&($e_1+e_2$, $e_2$, $e_1+e_2$, $e_1$)\\
		\hline
	\end{tabular}
	\caption{ All the characteristic vectors of small covers over the $2$-cube $I^2$. }
	\label{table:square}
\end{table}

 By the definition of small cover, the $GL_n(\mathbb{Z}_2)$-equivalence characteristic functions always give homeomorphic $n$-dimensional small covers. Besides, small covers from different $GL_n(\mathbb{Z}_2)$-equivalence classes may also be  homeomorphic. See Section \ref{section:classify} for more details. In the following sections, we always first find the $GL_n(\mathbb{Z}_2)$-equivalence representatives, and then group them by homeomorphism. For a $GL_4(\mathbb{Z}_2)$-equivalence representative $k$-coloring with coloring set $\{1,2,4,8,c_5,\cdots,c_k\}$, the $(k-4)$ colors $c_5,\cdots,c_k$ are called
 \emph{added colors}. Especially, there is no added color in the $4$-coloring case. In the following, the coloring  with the coloring set $C=\{1, 2, 4, 8, c_5,\cdots, ck\}$ is called a $(c_5,\cdots,c_k)$-coloring. The case that considers all of the $(c_5,\cdots,c_k)$-coloring is called the $(c_5,\cdots,c_k)$-coloring case.

H. Nakayama and Y. Nishimura have discussed the orientability of a small cover in \cite{NakayamaN:05}. Their main theorem is presented below.

\begin{theorem} \label{theorem: NakayamaN} \emph{(H. Nakayama-Y. Nishimura \cite{NakayamaN:05})} For a simple $n$-dimensional polytope $P$, and for a basis $\{e_1,\cdots,e_n\}$ of $\mathbb{Z}_2^n$, a homomorphism $\epsilon:\mathbb{Z}_2^n\rightarrow\mathbb{Z}_2=\{0,1\}$ is defined by $\epsilon(e_i)=1$, for each $i=1,\cdots, n$. A small cover $M(L,\lambda)$ is orientable if and only if there exists a basis $\{e_1,...,e_n\}$ of $\mathbb{Z}_2^n$ such that the image of $\epsilon\lambda$ is $\{1\}$.
\end{theorem}

Therefore, a small cover is orientable if and only if there exists a basis of $\mathbb{Z}_2^n$ such that all the colors used, i.e., the columns of the characteristic matrix, have an odd number of $1$ entries. In the space of $\mathbb{Z}_2^4$, the vectors $(1,0,0,0)$, $(0,1,0,0)$, $(0,0,1.0)$, $(0,0,0,1)$, $(1,1,1,0)$, $(1,0,1,1)$, $(1,1,0,1)$, and $(0,1,1,1)$, which are the binary form of $1$, $2$, $4$, $8$, $7$, $11$, $13$, and $14$, are the only eight elements with entry sums of $1$ mod $2$. These eight vectors are called \emph{orientable colors}. The seven non-trivial colors left are $(1,1,0,0)$, $(1,0,1,0)$, $(0,1,1,0)$, $(1,0,0,1)$, $(0,1,0,1)$, $(0,0,1,1)$ and $(1,1,1,1)$, which are the binary forms of $3$, $5$, $6$, $9$, $10$, $12$, and $15$. An \emph{orientable basis} in $\mathbb{Z}_2^4$ is defined to be a basis in $\mathbb{Z}_2^4$ that consists of four linear-independent orientable colors. Especially, the standard basis in $\mathbb{Z}_2^4$, i.e., $(1,0,0,0)$, $(0,1,0,0)$, $(0,0,1,0)$ and $(0,0,0,1)$, is an orientable basis. Note that, for an orientable color, the number of entries with value $1$ is always odd. In other words, when changing from one orientable basis to another orientable one, we actually add or remove even many $1$s from the previous characteristic matrix to form the new one. Hence, parity of the number of $1$s in each column is maintained under different orientable basis.  Therefore, we have the following corollary.

\begin{corollary}\label{corollary: NakayamaN}
	Suppose $(1,2,4,8,a_1,\cdots,a_m)$ is a characteristic vector of some $GL_4(\mathbb{Z}_2)$-equivalence representative $\lambda$ over a simple polytope $P$ Then the corresponding small cover $M(P,\lambda)$ is non-orientable if and only if there is some $a_i \in \{3,5,6,9,10,12,15\}$. .  
\end{corollary}


\section{Combinatorics of the 120-cell}\label{section:120cell}
In this section, we provide some fundamental combinatorial background information about the 120-cell $\mathcal{E}$ with ordered facets. We refer to \cite{S.Schleimer H.Segerman:15} for a nice exposition about some well-known facts on the 120-cell. This part is important in constructing and classifying small covers over the right-angled 120-cell. Most of the content in this section is well-known for experts; we present it here for the convenience of the readers.
\subsection{\textbf{Dodecahedron}}

The regular dodecahedron, or pentagonal dodecahedron, is abbreviated as  \emph{dodecahedron} and denoted as $\mathcal{D}$ in the following discussion. It has 12 regular pentagonal faces, 20 vertices and 30 edges, see Figure \ref{figure:block1} in Section 5, and is represented by the Schl\"{a}fli symbol \{5, 3\}. 

The \emph{full symmetry group} $\mathbb{A}(P^n)$ of an $n$-dimensional polytope $P^n$ is the group of all isometries under which the polytope $P^n$ is invariant.
The group operation is the composition. It is the same as the combinatorial automorphism group of the $n$-dimensional polytope $P^n$ if $P^n$ is regular. The group $\mathbb{A}(P^n)$ includes a subgroup of orientation-preserving isometries $\mathbb{A}^+(P^n)$ and a subset of  orientation-reversing isometries $\mathbb{A}^-(P^n)$. For the symmetry group of the dodecahedron $\mathcal{D}$, there is a well-known lemma (e.g., \cite{S.Schleimer H.Segerman:15}):

\begin{lemma}
	$|\mathbb{A}(\mathcal{D})|=120$ and $|\mathbb{A}^+(\mathcal{D})|=60$.
\end{lemma}

The group of orientation-preserving symmetries of the dodecahedron $\mathcal{D}$ is known as the alternating group on five letters, which is denoted as $A_5$. The full symmetry group is isomorphic to $A_5\times \mathbb{Z}_2$, with the second factor being generated by the antipodal map.

\subsection{\textbf{$SO(3)$ and rigid motion}}

The special orthogonal group, denoted by $SO(n)$, contains all $n\times n$ orthogonal matrices  of determinant 1. This is also the group of rigid motions of $\mathbb{R}^n$ fixing the origin. Such a group is also called the rotation group. See \cite{S.Schleimer H.Segerman:15} for detailed discussions in studying the symmetries of $\mathcal{D}$ from the rigid motion's point of view.

\subsection{\textbf{Quaternions}}

Let $<1,i,j,k>$ be the usual basis for $\mathbb{R}^4$. The quaternion space  $\mathcal{H}=\mathbb{R}\oplus\mathbb{I}=\mathbb{R}\oplus i\mathbb{R}\oplus j\mathbb{R}\oplus k\mathbb{R}$, where $\mathbb{I}=i\mathbb{R}\oplus j\mathbb{R}\oplus k\mathbb{R}$ is the subspace of the purely imaginary quaternions. The sum of two elements of $\mathcal{H}$ is defined to be the sum of them as elements of $\mathbb{R}^4$. Similarly the product of an element of $\mathcal{H}$ and a real number is defined to be the same as the scalar multiplication in $\mathbb{R}^4$. The associative quaternion multiplication is defined by first establishing the products of basis elements as follows:
\begin{center}
$i^2=j^2=k^2=ijk=-1$.
\end{center}
\noindent Then, the distributive law is applied to define the general situation.

Let $q=a+bi+cj+dk$ be a quaternion. The conjugate of $q$ is the quaternion $q^{*}=a-bi-cj-dk$. The conjugate of a product of two quaternions is the product of the conjugates in reverse order. That is, if $p$ and $q$ are quaternions, then $(pq)^* = q^*p^*$. The square root of the product of a quaternion with its conjugate is called its \emph{norm} and is denoted by $\vert\vert q\vert\vert$, namely $\vert\vert q\vert\vert=\sqrt{qq^{*}}=\sqrt{q^{*}q}=\sqrt {a^2+b^2+c^2+d^2}$. The inverse of a quaternion, under quaternion multiplication, is $\frac{q^*}{||q||^2}$. Besides, using the norm we can define the 3-sphere $\mathbb{S}^3$ as the set of unit quaternions $\{q\in \mathcal{H}~\vert ~\vert\vert q\vert\vert=1\}$. The $3$-sphere $\mathbb{S}^3$ forms a group under quaternion multiplication. The point $1+0i+0j+0k$ serves as the identity; the associativity and inversion are inherited from $\mathcal{H}$. For two points $p$ and $q$ in $\mathcal{H}$, the norm of the difference between them is defined as $d(p, q)=\vert\vert p-q\vert\vert$. The metric on $\mathcal{H}$ naturally induces the round metric $d_S(p,q)=arccos(<p,q>)$ on the $3$-sphere. According to \cite{S.Schleimer H.Segerman:15}, the group structure and geometry of $\mathbb{S}^3$ interact with each other in the following way:

\begin{lemma}\label{lemma:120symmery group}
Both left and right actions of $\mathbb{S}^3$ on $\mathcal{H}$ are via orientation-preserving isometries in $\text{Isom}^+(\mathcal{H}^3)$. The same holds for the 3-sphere's action on itself.
\end{lemma}

\begin{proof}
For $p\in\mathbb{S}^3$ and $q,r\in\mathcal{H}$, we have $d_{\mathcal{H}}(pq,pr)=\vert pq-qr\vert=\vert p\vert \vert q-r\vert=\vert q-r\vert=d_{\mathcal{H}}(q,r)$. Similarly, $d_{\mathcal{H}}(qp,rp)=d_{\mathcal{H}}(q,r)$. Since $\mathbb{S}^3$ is connected and the identity $1$ acts trivially, the actions are orientation-preserving.
\end{proof}

Furthermore, for $q\in \mathbb{S}^3$, we define a twist action $\phi_q:\mathcal{H}\rightarrow\mathcal{H}$,
$p\mapsto \phi_q(p)=qpq^{-1}$, which is again via isometries. Note that $\phi_q$ preserves $\mathbb{I}\subset\mathcal{H}$ setwise. Thus we can define $\psi_q:=\phi_q|_\mathbb{I}:\mathbb{I}\rightarrow\mathbb{I}$. Recall that $SO(3)$ is the group of the three-by-three matrices with determinant one. Take $\langle i, j, k\rangle$ as a basis for $\mathbb{I}$, the group $SO(3)$ can be identify with $\text{Isom}_0^+(\mathbb{I})$, the group of orientation-preserving isometries of $\mathbb{I}$ fixing the origin. 
Hence the map $\psi_q$ is an element of $SO(3)$ and the reduced map $\psi:\mathbb{S}^3\rightarrow SO(3)$ is a double cover \cite[Lemma 4.8, Lemma 4.9]{S.Schleimer H.Segerman:15}. 
If $G\subset SO(3)$ is a group, then $G^*=\psi^{-1}(G)$ is called \emph{the binary group corresponding to G}. In particular, for the group of orientation-preserving isometries of the dodecahedra $\mathcal{D}$ $\mathbb{A}^+(\mathcal{D})\subset SO(3)$, we have the \emph{binary dodecahedral group} $\mathcal{D}^*\triangleq\psi^{-1}(\mathbb{A}^+(\mathcal{D}))\subset \mathbb{S}^3$, which has 120 elements in total. More precisely, $\mathcal{D}^*$, consist of 24 Hurwitz units $\{ \pm1,\pm i,\pm j,\pm k,\frac{1}{2} (\pm 1\pm i\pm j \pm k )\}$ and 96 quaternions  obtained from $\frac{1}{2}( 0\pm i\pm \phi^{-1}j \pm \phi k )$ by the even permutations of all four coordinates \cite{coxeter:1969}. 

 \subsection{\textbf{The 120-cell $\mathcal{E}$}} 
 We are going to formulate the construction of the 120-cell by means of group $\mathcal{D}^*\subset \mathbb{S}^3$.
 \begin{definition}
 	Suppose $V$ is a finite set of points in a metric space $X$, the \emph{Voronoi cell} about a point $q\in V$ is the set $\text{Vor}(q)=\{r\in X~|~\text{for~~ all}~~ p\in V,~~d_X(q,r)\leq d_X(r,p)\}$.
 \end{definition}
\noindent Two cells $\text{Vor}(p)$ and $\text{Vor}(q)$ are said to be adjacent if $\text{Vor}(p)\cap \text{Vor}(q)$ is an $(n-2)$-cell. Let $\mathcal{T}_{120}$ be the tiling of the $3$-sphere by the cell $\{\text{Vor}(q)~\vert~ q\in \mathcal{D}^*\}$ and define $\mathcal{C}=\text{Sym}(\mathcal{T}_{120})$. Here is a useful lemma:
\begin{lemma}\cite[Lemma 5.3,Lemma 5.4]{S.Schleimer H.Segerman:15}\label{120sym}
	\begin{enumerate}
		\item The left action of $\mathcal{D}^*$ on $\mathcal{T}_{120}$ is transitive on the $3$-dimensional cells.
		\item The twisted action of $\mathcal{D}^*$ fixes $\text{Vor}(1)$ setwise.
		\item Both above actions give homomorphsims of $\mathcal{D}^*$ to $C$. 
		\item Each cell $\text{Vor}(q)$ is a regular spherical dodecahedron.
	\end{enumerate}
\end{lemma}
 
 \begin{definition}
 	The Euclidean 120-cell $\mathcal{E}$ is the convex hull, taken in $\mathcal{H}$, of the vertices of $\mathcal{T}_{120}$.
 \end{definition}

 \begin{remark}
 	The convex hull in $\mathbb{R}^4$ of these 120 points is the regular 4-dimensional polytope 600-cell, which is dual to the 120-cell $\mathcal{E}$. If we forget the group structure of $\mathcal{D}^*$ and regard it as a set, it is the vertex set of the 600-cell.
 \end{remark}
 
 \begin{remark}
 	  The right-angled 120-cell $\mathcal{E}$ is a 4-dimensional compact hyperbolic polytope of dihedral angle $\frac{\pi}{2}$ and volume $\frac{34}{3}\pi^{2}$ \cite{Martelli:2018}. It has the same face lattice with the Euclidean 120-cell.
 \end{remark}
 
\begin{figure}[H]
 	\scalebox{0.22}[0.22]{\includegraphics {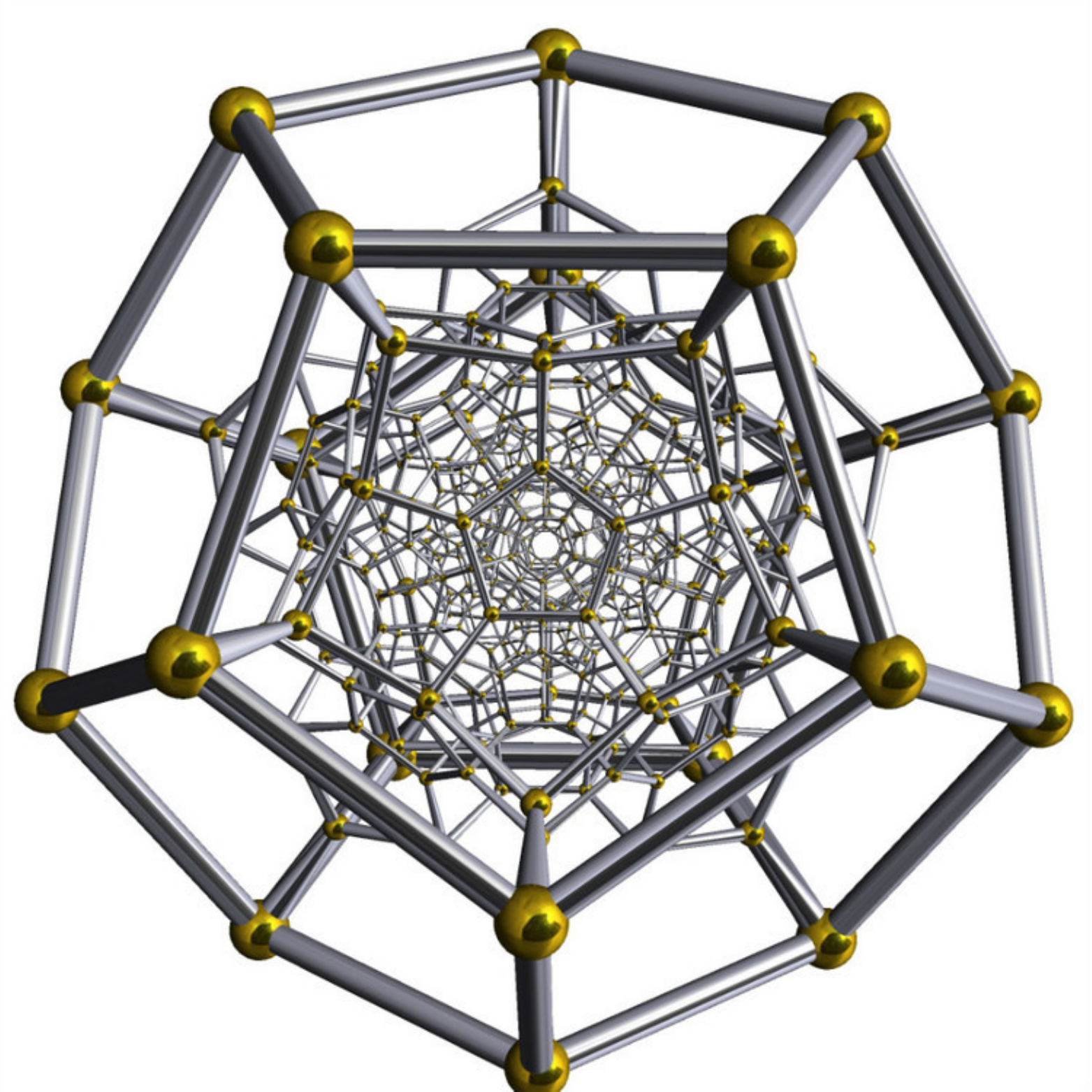}}
 	\vspace*{0.25cm}
 	\caption{Facet boundary of the 120-cell, from the web.}\label{figure:cell}
 \end{figure}

 The 120-cell has 120 dodecahedral facets, 720 pentagonal ridges, 1200 edges, and
 600 vertices. There are four facets meeting at each vertex. Figure \ref{figure:cell} shows a layer structure of the facets of the 120-cell. The first layer is the unbounded dodecahedral facet that contains $\infty$. The second layer consists of 12 dodecahedral facets. The adjacencies between the second and third layers can be seen in Figure \ref{figure:cell} with little effort. It may be easily checked that the 120-cell is a regular polytope, for example see \cite[Theorem 5.7]{S.Schleimer H.Segerman:15}. 

 An $n$-polytope $P^n$ is called \emph{simple} if every $r$-face is contained in exactly $n-r$ facets.
 A simple polytope is called a \emph{flagtope} if every collection of its pairwise intersecting facets has a nonempty intersection. If a simple polytope is a flagtope, all of the adjacency information is included in the adjacency matrix.  The 120-cell $\mathcal{E}$ can be verified to be a flagtope with little effort. Moreover, the Euclidean 120-cell defined above is inscribed in the 3-sphere $\mathbb{S}^3$ and we can define a projection $\Pi:\partial\mathcal{E}\rightarrow \mathbb{S}^3,x\mapsto \frac{x}{||x||}$. We fix the orders of the 120 facets of 120-cell and 120 elements of $\mathcal{D}^*$ by $F_1,F_2,\cdots F_{120}$ and $q_1,q_2,\cdots, q_{120}$ as shown in Tables \ref{table:coor1}--\ref{table:coor3}, where $\Pi(F_i)=\text{Vor}(q_i),q_i\in\mathcal{D}^*$. For brevity, we re-denote $\text{Vor}(q_i)$ as $\text{Vor}(i)$ and define a function $ord:\mathcal{D}^*\rightarrow \{1,2,\cdots,120\},~ q_i\mapsto ord(q_i)=i$ to obtain the label. The symbols $F_i$, $q_i$ and $ord$ will be used throughout the rest of the paper with this meaning, unless stated otherwise. Suppose $X(\mathcal{E})=(a_{ij})_{m\times m}$ is the \emph{adjacency matrix} of the 120-cell, then the entries can be defined as
	\begin{center}
		$a_{ij}=\left\{
		\begin{array} {rcl}1 &\mbox{
			if $Vor(i)$ and $Vor(j)$ are adjacent in $\mathcal{T}_{120}$} \\ 0 & \mbox{otherwise}
		\end{array}\right.$.
	\end{center}

We calculate distances from elements in $\mathcal{D}^*$, i.e., the vertices dual to the facets of the 120-cell, to the identity $1+0i+0j+0k$ under the round metric. Then, the facets of $\mathcal{E}$ can be divided into nine spherical layers $L_i$ according to the distances as shown in Table \ref{table:layers}.  A detailed geometric illustration is available in \cite[Section 6]{S.Schleimer H.Segerman:15}.

\begin{figure}[H]
	\scalebox{0.56}[0.56]{\includegraphics {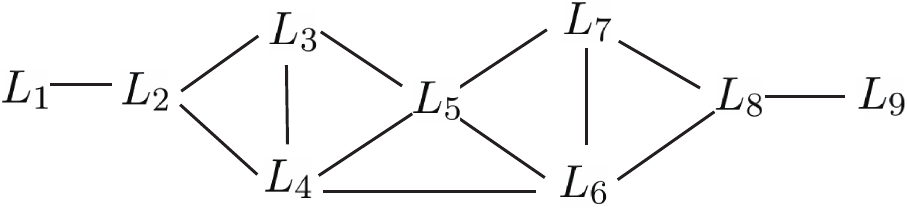}}
	\caption{Adjacency relations among the nine layers.}\label{figure:AD}
\end{figure}

\begin{table}[H]
	\caption{Nine-layer structure of the facets of the 120-cell $\mathcal{E}$.}
	\label{table:layers}
	\scalebox{0.98}{
		\linespread{1}\selectfont
	\begin{tabular}{|c|c|c|}
		\hline
		Layer&Distance&$\#$ Voronoi cells\\

		\hline	
		$L_1$ & 0 & 1 \\
		\hline
		$L_2$ & $\pi$/5& 12 \\
		\hline
		$L_3$ & $\pi$/3 & 20 \\
		\hline
		$L_4$ & $2\pi$/5 & 12 \\
		\hline
		$L_5$ &$\pi$/2 & 30 \\
		\hline
		$L_6$ & $3\pi$/5& 12 \\
		\hline
		$L_7$ & $2\pi$/3& 20 \\
		\hline
		$L_8$ & $4\pi$/5 & 12 \\
		\hline
		$L_9$ & $\pi$& 1 \\	
		\hline
		
	\end{tabular}
}
\end{table}
Note that the distances from the center of the facet $F_1$ to the centers of those facets adjacent to it, namely the 12 facets in the second layer, are always $\frac{\pi}{5}$, under the round metric on $\mathbb{S}^3$. Due to the symmetry, the distance between the two centers of every adjacent facets should always be $\frac{\pi}{5}$. Therefore, we can probe the adjacency matrix $X(\mathcal{E})$ by examining the distance under the round metric $d_S(q_i,q_j)$ for all $p,~ q\in \mathcal{D}^*$.  
For a better demonstration of the adjacency matrix $X(\mathcal{E})$, we first depict the adjacency relations among different spherical layers $L_i$ ($1\le i\le 9$). We say two layers $L_{i}$ and $L_j$ are \emph{adjacent} if there are two facets $F_{i_p}\in L_i$ and $F_{j_q}\in L_j$ such that $F_{i_p}$ and $F_{i_q}$ are adjacent. Adjacency relations are shown in Figure \ref{figure:AD}, where two layers are connected by edge if they are adjacent.     

All of the adjacency matrices  $X_{L_i,L_j}$, abbreviated as $X_{ij}$ of the adjacent layers, for $1\le i< j\le 9$, are shown in Figures \ref{figure:aj1} -- \ref{figure:aj8} in Appendix.  For clarity, we represent the matrix by a colored grid; the uncolored entries are all of value 0. The diagonal items are marked specially with dark shading and other positions with values of 1 are filled with light shading. The figures can be obtained from Table \ref{table:coor1}-\ref{table:coor3} in Appendix. The whole symmetric adjacency matrix $X(\mathcal{E})$ for the 120-cell $\mathcal{E}$ is
\begin{center}
	$X_{\mathcal{E}}=\begin{pmatrix}
	X_{11} & X_{12} & 0 & 0 & 0 & 0 & 0 & 0 & 0 \\
	X_{21}& X_{22} & X_{23} & X_{24} & 0 & 0 & 0 & 0 & 0\\
     0& X_{32} & X_{33} & X_{34} &  X_{35} & 0 & 0 & 0 & 0\\
	0& X_{42} & X_{43} & X_{44} &  X_{45} & X_{46} & 0 & 0 & 0\\
	0&0 & X_{53} & X_{54} &  X_{55} & X_{56} &X_{57} & 0 & 0\\
	0&0 &0 & X_{64} &  X_{65} & X_{66} &X_{67} & X_{68}  & 0\\
	0&0 &0 & 0 &  X_{75} & X_{76} &X_{77} & X_{78}  & 0\\
	0&0 &0 & 0 &  0 & X_{86} &X_{87} & X_{88}  & X_{89}\\
	0&0 &0 & 0 &  0 & 0 &0 & X_{98}  & X_{99}\\
	\end{pmatrix}$,
\end{center}
where $X_{ji}=X_{ij}^T$.

It can be read from the matrix $X_{12}$ that the first four facets of the 120-cell $\mathcal{E}$ are adjacent to each other. Since the 120-cell is a flagtope, these four facets would share a common vertex. We may always color the facet $F_i$ by $e_i$, $i=1,~2,~3,$ and $4$, where $\{e_i\}_{i=1}^4$ is the standard basis set of $\mathbb{Z}_2^4$.

\subsection{\textbf{Symmetry group of 120-cell $\mathcal{E}$}}\label{section:symmetry of 120}The symmetry group, i.e., the combinatorial automorphism group, of the 120-cell $\mathcal{E}$, denoted as $\mathbb{A}(\mathcal{E})$, is well known as the Coxeter group $H_4$. However, that is not enough for our purpose. We want to describe the group $\mathbb{A}(E)$ as the subgroup of the permutation group $S_{120}$, writing down precisely each element in $\mathbb{A}(E)$  as a permutation of 120 facets. Each permutation tells how the facets of the 120-cell are permuted under certain symmetry. The approach we adopt to find all of the permutations is similar to the method used in the proof of \cite[Theorem 5.7]{S.Schleimer H.Segerman:15}. The permutation form of $\mathbb{A}(E)$ is of great importance for classifying small covers over the right-angled 120-cell by homeomorphism in Section \ref{section:classify}. 

Let $\mathbb{X}^n\in \{\mathbb{S}^n,\mathbb{E}^n,\mathbb{H}^n\}$ be one of the three standard $n$-dimensional geometric spaces of constant curvature $1$, $0$ or $-1$ and $P^n$ be an $n$-dimensional polytope in the $n$-space $\mathbb{X}^n$. A collection of faces $Q_0\subset Q_1\subset \cdots\subset Q_{n-1}\subset Q_n=P^n$ is called a \emph{flag of the polytope} $P^n$ if $Q_i$ has dimension $i$. Suppose that $\mathcal{Q}= \{Q_i\}$ is a flag in $P^n$. The corresponding \emph{flag polytope}
$P^n_\mathcal{Q}$ is defined to be the convex hull of the centres of the $Q_i$.  
Similarly, we can also define the \emph{flag of a tessellation $\mathcal{T}$} of $n$-space $\mathbb{X}^n$. A collection of cells  $Q_0\subset Q_1\subset \cdots\subset Q_{n-1}\subset Q_n$ is called a \emph{flag of the tessellation $\mathcal{T}$} if $Q_i$ is an $i$-cell in the tiling $\mathcal{T}$. Suppose that $\mathcal{Q}= \{Q_i\}$ is a flag in $\mathcal{T}$. The corresponding \emph{flag polytope}
$\mathcal{T}_\mathcal{Q}$ is the convex hull in $\mathbb{X}^n$ of the centres of the $Q_i$.

The orientation preserving
symmetry group of the 120-cell is the quotient of the product of two copies
of the binary icosahedral group by the group generated by $(-1,-1)$. More precisely, 
the flags of the 120-cell $\mathcal{E}$ in $\mathcal{H}$ are 4-dimensional simplices with one vertex at the origin, which are in one-to-one correspondence with the 14,400 spherical flag tetrahedra of $\mathcal{T}_{120}$. Next, we fix a right-handed spherical flag tetrahedron $T$ of $\text{Vor}(1)$. By Lemma \ref{120sym}, we may use the left action of $\mathcal{D}^*$ to transfer
every tetrahedron $T^{'}\in\mathcal{T}_{120}$ into $\text{Vor}(1)$. Then, we use the twisted action of $\mathcal{D}^*$ to send $T^{'}$ to $T$. We now want to find the distinct 14,400 elements. 

First, we can obtain a symmetry element of $\mathbb{A}^+(\mathcal{E})$ by multiplying one element of the set  $\mathcal{D}^*$ from left and one element of the set  $\mathcal{D}^*$ from the right (the two factors may be the same) to all of the centers of the Voronoi cells in $\mathcal{T}_{120}$, i.e., all of the elements in $\mathcal{D}^*$. Namely, for $p,q\in\mathcal{D}^*$, define $\psi_{p,q}:\mathcal{D}^*\rightarrow \mathcal{D}^*, q_i\mapsto \psi_{p,q}(q_i)=pq_ip$. For every pair $(p,q)\in \mathcal{D}^*\times \mathcal{D}^*$, we have a permutation $[ord(\psi_{p,q}(q_1)), ord(\psi_{p,q}(q_2)),\cdots,ord(\psi_{p,q}(q_{120}))]$ of an orientation-preserving symmetry element in $\mathbb{A}^+(\mathcal{E})$. Note that the symmetry obtained by left-acting $a_1+b_1 i +c_1 j+d_1k$ and right-acting $-a_2+b_2 i +c_2 j+d_2k$ is the same as the symmetry obtained by left-acting $-a_1-b_1 i -c_1 j-d_1k$ and right-acting $a_2-b_2 i -c_2 j-d_2k$. Then, among all $120\times120=14400$ results, there are only 7200 left after removing duplicates. Secondly, the conjugation map $a+bi+cj+dk\rightarrow a-bi-cj-dk$ is a product of three reflections, so is orientation reversing in $\mathcal{H}$. It preserves $\mathbb{S}^3$ and is again orientation reversing there. Since $\mathcal{D}^*$ is a group of quaternions, it is closed under conjugation. Since the tiling $\mathcal{T}_{120}$ is metrically defined in terms of $\mathcal{D}^*$, it is also invariant under conjugation. Hence, if we do the quaternion conjugation to each of the above 7200 results, we can obtain the corresponding orientation-reversing ones. So, there are also 7200 orientation-reversing symmetries of the 120-cell $\mathcal{E}$. Then, we have searched 14400 different elements in the symmetry group $\mathbb{A}(\mathcal{E})$. From the well-known result that $\vert\mathbb{A}(\mathcal{E})\vert=14400$, it can be claimed that all of the elements in $\mathbb{A}(\mathcal{E})$ have been depicted clearly by showing how its facets are permuted.

\section{Constructing small covers with the color-recursion algorithm}

\subsection{Color-recursion algorithm}
The \textbf{color-recursion algorithm } is broadly applicable. Here, we use it on the 120-cell $\mathcal{E}$ to generate small covers.  The 120-cell is dual to the 600-cell. By the dual non-singularity condition introduced in Section 2.1, if the convex hull $conv\{v_{i_1},v_{i_2},v_{i_3},v_{i_4}\}$ of four vertices $v_{i_1},v_{i_2},v_{i_3},v_{i_4}$ of the 600-cell is a tetrahedral facet of the 600-cell, the images $\lambda(v_{i_1}),\lambda(v_{i_2}),\lambda(v_{i_3}),\lambda(v_{i_4})$ shall generate a subgroup isomorphic to $\mathbb{Z}_2^4$. The color of the vertex is the color on the corresponding dual  facet. We now want to exclude the cases that fail to meet the non-singularity condition by systematically ``checking simplexes". For an $n$-dimensional simple polytope, every facet meets exactly $n-1$ other facets at a vertex.  Thus, we need to conduct the $i$-simplex criterion checks, where $1\leq i\leq n-1$. For the 4-dimensional polytope 120-cell $\mathcal{E}$, we only need to check the cases of 1-simplex, 2-simplex, and 3-simplex.

Given a coloring set $C=\{1,2,4,8,c_1,\cdots,c_k\}$ and suppose the first $(i-1)$ facets are colored, we are going  to color the facet $F_i$. Now, the initial candidate set of admissible colors for $F_i$, denoted by $d^0(F_i)$, is $C$. For a color $c\in C$, we undertake the following three tests to conclude whether we can color the facet $F_i$ by the color $c$ or not

\textbf{(a). 1-simplex checking}

\begin{figure}[H]
	\scalebox{0.4}[0.4]{\includegraphics {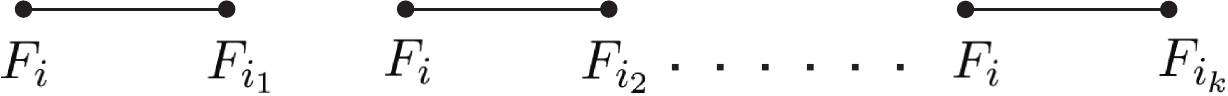}}
	\caption{The 1-simplex check.} \label{figure:1simplex}
\end{figure}

By the adjacency matrix $X_{\mathcal{E}}$, we can obtain the set $\mathcal{A}^1(F_i)$ of the facets that are adjacent to $F_i$ and have been colored. If $\mathcal{A}^1(F_i)=\emptyset$, all of the $2^n-1$ colors are qualified. If $\mathcal{A}^1(F_i)=\{F_{i_1},F_{i_2},...,F_{i_k}\}\ne \emptyset$,  there are $k$ 1-simplexes in the 600-cell as shown in Figure \ref{figure:1simplex}. We still use the notation $F_i$ rather than introducing new notation $v_i$ of the dual vertex because they are one-to-one dual to each other. This convention applies throughout the section. The non-singularity condition requires that these 1-simplexes do not have the same colors on their two vertices.
Let $c(F_{i_j})\in\mathbb{Z}^n_2$ be the color of colored facet $F_{i_j}\in\mathcal{A}^1(F_i)$. Now the set of admissible colors for $F_i$ is $$d^1(F_i)=d(F_i)-\{c(F_{i_1}),c(F_{i_2}),\cdots,c(F_{i_k})\}.$$
If $c\in d^1(F_i)$, we move on to the next checking. Otherwise, we check another candidate of color in the coloring set $C$. 

\textbf{(b). 2-simplex checking}

\begin{figure}[H]
	\scalebox{0.38}[0.38]{\includegraphics {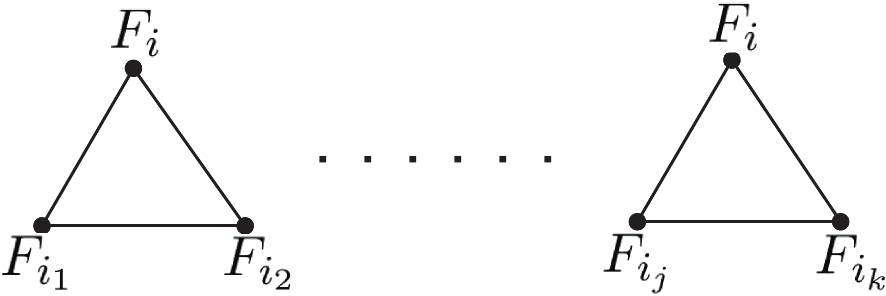}}
	\caption{The 2-simplex check.}\label{figure:2simplex}
\end{figure}

Next, we find from  $\mathcal{A}^1(F_i)$ the set $\mathcal{A}^2(F_i)$ of pairs of adjacent facets in $\mathcal{A}^1(F_i)$. In turn, each pair together with the facet $F_i$ provides us with three vertices in the 600-cell. Then, the complete graph on these vertices give us the 1-skeleton of a 2-simplex as shown in Figure \ref{figure:2simplex}. \noindent The non-singular condition requires that no color of the vertices for every 2-simplex can be represented as the $\mathbb{Z}_2$-coefficient linear combination of colors of the other two. If $\mathcal{A}^2(F_i)=\emptyset$, all of the candidates are qualified. Otherwise, assume $\mathcal{A}^2(F_i)=\{(F_{i_1},F_{i_2}),\cdots,(F_{i_j},F_{i_k})\}$. Now, the set of admissible colors for $F_i$ is  $$d^2(F_i)=d^1(F_i)-Span\{c(F_{i_1}),c(F_{i_2})\}-\cdots-Span\{c(F_{i_j}),c(F_{i_k})\}.$$
If $c\in d^2(F_i)$, we move on to the next checking. Otherwise, we check another candidate of color in the coloring set $C$. 

\textbf{(c). 3-simplex checking}

\begin{figure}[H]
	\scalebox{0.38}[0.38]{\includegraphics {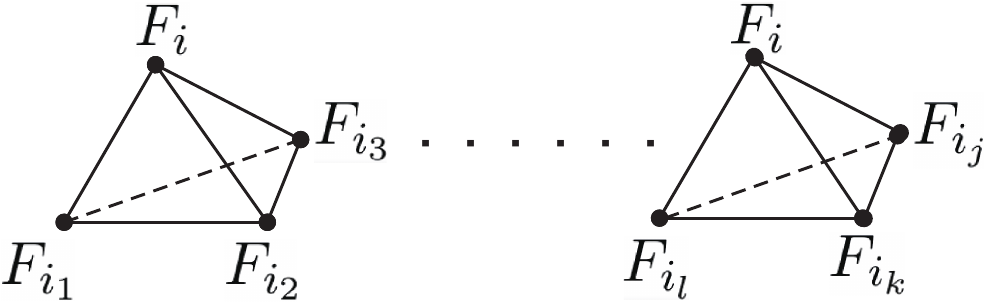}}
	\caption{The 3-simplex check.}\label{figure:3simplex}
\end{figure}

Finally, we find from  $\mathcal{A}^2(F_i)$ the set $\mathcal{A}^3(F_i)$ of facet triples such that every two of the three facets are adjacent. In turn, each triple together with the facet $F_i$ provides us with four vertices. Then, the complete graph on these vertices gives us the 1-skeleton of a tetrahedron as shown in Figure \ref{figure:3simplex}.

If $\mathcal{A}^3(F_i)=\emptyset$, all of the colors are qualified. Otherwise, assume $\mathcal{A}^3(F_i)=\{(F_{i_1},F_{i_2},F_{i_3}),\cdots,$ $(F_{i_j},F_{i_k},F_{i_l})\}$, then
$$d^3(F_i)=d^1(F_i)-Span\{c(F_{i_1}),c(F_{i_2}),c(F_{i_3})\}-\cdots-Span\{c(F_{i_j}),c(F_{i_k}),c(F_{i_l})\}.$$
If $c\in d^3(F_i)$, we can color $F_i$ by the color $c$. Otherwise, we turn to check another candidate of color in the coloring set $C$. 

We keep doing the above procedure until we successfully color $F_{120}$ by admissible colors. If so, the resulting characteristic vectors are printed to an assigned file. 
We call it the ``color-recursion" algorithm because we use the recursion method to generate the coloring. It automatically turns to color the $F_{i+1}$-th facet when the facet $F_i$  is colored; or it traces back to take an alternative color among the coloring set. In particular, if all of the candidate colors in the coloring set encounter $d^1(F_i)=\emptyset$, $d^2(F_i)=\emptyset$ or $d^3(F_i)=\emptyset$,  we may go back to change the color of the $(i-1)$-th facet and do a new round of checking. If there are still no other admissible choices for $F_{i-1}$, we move one more facet backwards, namely to the $(i-2)$-th facet, to check. Since the first four mutually adjacent facets have been colored by the standard basis in $\mathbb{Z}_2^4$, we start from facet $F_5$. Therefore, tracing-back ends when we fail in supplying a suitable color for the facet $F_5$. The pseudocode for the color-recursion algorithm is given in Algorithm 1.

Moreover, by means of the stabilizer of a vertex,  we can actually select out the representatives of the added colors for every $k$-coloring case to substantially reduce the computational load.
\begin{lemma}\label{lemma:S4}
	The stabilizer of every vertex in the symmetry group of the regular polytope 120-cell $\mathcal{E}$ is $S_4$, the symmetry group over four letters.
\end{lemma}
For example, suppose the first four facets that  intersect at a common vertex $v$ are colored by $e_1,e_2,e_3,$ and $e_4$. Then one more color $e_1+e_3$, the combination of the first and the third colors, is added to the coloring set. Namely, the coloring set is $\{1,2,4,8,3\}$. If we act upon the first four facets by a symmetry in  $S_4$ to generate, for example,  a permutation as follows
\begin{equation}\label{matrix:2}
\begin{pmatrix}
\begin{smallmatrix}
e_1&\	e_2&\	e_3&\	e_4\\
e_2&\	e_1&\	e_3&\	e_4\\
\end{smallmatrix}
\end{pmatrix},
\end{equation}
then the added color, which is the combination of the first and third colors, changes from $e_1+e_3$ to $e_2+e_3$. Namely, the coloring set becomes $\{1,2,4,8,5\}$. That is, the manifolds that can be constructed by $\{1,2,4,8,3\}$ and $\{1,2,4,8,5\}$ are actually the same up to homeomorphism. 
By the same reason, the coloring sets $\{1,2,4,8,3\}$, $\{1,2,4,8,5\}$, $\{1,2,4,8,6\}$, $\{1,2,4,8,9\}$,

		\begin{algorithm}[H]
			\footnotesize
			\label{algorithm1}
			\caption{Pseudocode for the color-recursion algorithm}
			\begin{algorithmic}
				
				\Function{Check}{$i, s, vector$} \Comment{{\color{darkgray}to check whether the  $i$-th facet can be colored by  $c$ or not. }}
				\State $A\gets$ the set of labels of facets that adjacent to the $i$-th facet 
				\ForAll {$j$ in $A$}
				\If {$j>i$}
				\State \textbf{break}
				\EndIf
				\If {the $j$-th entry of $vector$ is $s$ } \Comment{{\color{darkgray}1-simplex checking}}
				\State\Return \textbf{FALSE}
				\EndIf
				\ForAll {$p$ in $\{x~\vert~x>j,~x\in A\}$}
				\If {$p>i$}
				\State \textbf{break}
				\EndIf
				\If {the $j$-th and $p$-th facets are adjacent \& the three colors, namely the color $s$ as well as the colors of the $j$-th and the $p$-th facets, failed to satisfy the 2-simplex criterion}
				\State 	
				\Comment{{\color{darkgray}2-simplex checking}}
				\State\Return \textbf{FALSE}
				\EndIf
				\ForAll {$q$ in $\{x~\vert~x>p,~x\in A\}$}
				\If {$q>i$}
				\State \textbf{break}
				\EndIf
				\If {the $j$-th , $p$-th and $q$-th facet are mutually adjacent \& the four colors,  namely the color $s$ as well as the colors of the $j$-th, $p$-th and $q$-th facets, failed to satisfy the 3-simplex criterion}
				\State 	\Comment{{\color{darkgray}3-simplex checking}}
				\State\Return \textbf{FALSE}
				\EndIf
				\EndFor
				\EndFor
				\EndFor
				\State\Return \textbf{TRUE}
				\EndFunction
				
				\vspace{0.5cm}
				\Function{ColorRecursion}{$i$,$vector$}
				\If {$i=120$}
				\State \textbf{print} vector 
				
				\Comment{{\color{darkgray} When the facet $F_{120}$ is colored successfully, the characteristic vector of the $GL_4(\mathbb{Z}_4)$-equivalence representative would be printed to an assigned txt file. }}
				\EndIf
				\ForAll {$s$ in the given coloring set $C$}
				\State the $i$-th item of vector is valued by $s$
				\If {\Call{Check}{$i$,$s$,vector}}
				\State {\Call{ColorRecursion}{$i+1$,vector}}
				\EndIf
				\EndFor
				\EndFunction
				
				\vspace{0.5cm}
				\Require
				\Statex $X(\mathcal{E})_{120\times 120}$: Adjacency matrix of the right-angled 120-cell $\mathcal{E}$
				\Statex $v_{1\times 120}$:$(1,2,4,8,0,0,\cdots,0,0)$\footnotemark[1] \Comment{{\color{darkgray} All but the first four entries of $v$ are with value $0$}}
				\Statex coloring set $C$: $\{1,2,4,8,c_1,\cdots,c_k\}$
				\vspace{0.5cm}
				
				\vspace{0.5cm}
				\Ensure
				\Statex	\textbf{T1}=current system time
				\Statex \Call{ColorRecursion}{5,$v$}\footnotemark[2]
				
				\Comment{{\color{darkgray}All the $GL_4(\mathbb{Z}_2)$-equivalence  representatives with coloring set $C$ will be printed}}
				\Statex\textbf{T2}=current system time
				\Statex\textbf{print} 'running time:',\textbf{T2}-\textbf{T1}		
			\end{algorithmic}
		\end{algorithm}
		
		\customizedfootnotetext{1,2}{input $v_{1\times 120}=(0,0,\cdots,0)$ with all the entries valued 0 and run $\textsc{ColorRrecursion}(1,v)$. We obtain all of the characteristic vectors. }

\noindent $\{1,2,4,8,10\}$ and $(1,2,4,8,12\}$ are in the same family.  Therefore, through Lemma \ref{lemma:S4}, we can classify the $k$-colorings over the 120-cell $\mathcal{E}$ by considering the stabilizer group action. The representatives of the added colors are summarized in Table \ref{table:rep} in Appendix.  This  means we have the following corollary, which is a standard claim for the 5-coloring case and can be adopted by any $k$-coloring case based on Table \ref{table:rep}.
\begin{corollary}\label{coro:9}
	For every small cover $M(\mathcal{E},\lambda)$ with $\vert Im(\lambda)\vert=5$, there exists another characteristic function $\lambda^*$ with $Im(\lambda^*)\in\{\{1,2,4,8,3\},\{1,2,4,8,7\},\{1,2,4,8,15\}\}$ such that $M(\mathcal{E},\lambda^*)\cong M(\mathcal{E},\lambda)$.
\end{corollary}

Garrison-Scott has considered the case of (1,2,4,8,15)-coloring in \cite{Scott:02}. The added color $15$ is actually one of the three representatives when $k=5$ in Table \ref{table:rep}. Up to $GL_4(\mathbb{Z}_2)$-equivalence, there are only ten characteristic functions in this $(1,2,4,8,15)$-coloring case. Moreover, they are mutually homeomorphic. We discuss the details in Section \ref{section:classify}.

\subsection{\{4,~5,~6,~7,~8\}-coloring small covers over the 120-cell $\mathcal{E}$}
Applying the color-recursion algorithm to $4$-coloring, $5$-coloring, $6$-coloring, and $7$-coloring situations, abbreviated as $\{4,5,6,7\}$-coloring cases, we have the following Proposition \ref{prop:4567}:

\begin{proposition}\label{prop:4567}
	There is only one 5-coloring small cover, i.e., the Garrison-Scott manifold, up to homeomorphism, and no $\{4,6,7\}$-coloring small cover over the right-angled 120-cell $\mathcal{E}$.
\end{proposition}

	\begin{table}[H]
		\caption{Code-recursion algorithm computation of 8-coloring small covers over the 120-cell.}
		\label{table:compare1}
		\scalebox{0.8}{
			\linespread{1}\selectfont
			\begin{tabular}{|p{1em}|c|c|c|c|p{1em}|c|p{1em}|c|c|c|p{1em}|c|p{1em}|}  
				\hline
				& \multicolumn{4}{c|}{\multirowcell{3}{Added colors}} & & \multirow{3}{5em}{\centering Cardinality of the class} & & \multicolumn{3}{c|}{\multirowcell{2}{Running time of \\ Color-Recursion Algorithm}}  & & \multirow{3}{8em}{\centering Number of resulting vectors} & \\
				& \multicolumn{4}{c|}{}                            & &                                          & & \multicolumn{3}{c|}{}                             & &                                          & \\ \cline{9-11}
				& \multicolumn{1}{p{2em}}{} & \multicolumn{1}{p{2em}}{} & \multicolumn{1}{p{2em}}{} & \multicolumn{1}{p{2em}|}{} & & \multicolumn{1}{p{8em}|}{} & & \multicolumn{1}{p{4em}|}{\centering Seconds} & \multicolumn{1}{p{4em}|}{\centering Hours} & Days & & & \\
				\hhline{|~|-|-|-|-|~|-|~|-|-|-|~|-|~|}
				& 3 & 5  & 6  & 7  & & 4  & & 0.4      & 　    & 　  & & 0  & \\ \hhline{|~|-|-|-|-|~|-|~|-|~|~|~|-|~|}
				& 3 & 5  & 6  & 9  & & 12 & & 0.7      &       &     & & 0  & \\ \hhline{|~|-|-|-|-|~|-|~|-|~|~|~|-|~|}
				& 3 & 5  & 6  & 11 & & 12 & & 1.4      &       &     & & 0  & \\ \hhline{|~|-|-|-|-|~|-|~|-|-|-|~|-|~|}
				& 3 & 5  & 6  & 15 & & 4  & & \cellcolor{gray!50}145175.0 & \cellcolor{gray!50}40.3  & \cellcolor{gray!50}1.7 & & 0  & \\ \hhline{|~|-|-|-|-|~|-|~|-|-|-|~|-|~|}
				& 3 & 5  & 7  & 9  & & 12 & & 0.7      & 　    & 　  & & 0  & \\ \hhline{|~|-|-|-|-|~|-|~|-|~|~|~|-|~|}
				& 3 & 5  & 7  & 10 & & 24 & & 1.4      &       &     & & 0  & \\ \hhline{|~|-|-|-|-|~|-|~|-|~|~|~|-|~|}
				& 3 & 5  & 7  & 11 & & 24 & & 1.5      &       &     & & 0  & \\ \hhline{|~|-|-|-|-|~|-|~|-|-|~|~|-|~|}
				& 3 & 5  & 7  & 14 & & 12 & & \cellcolor{gray!10}72839.0  & \cellcolor{gray!10}20.2  &     & & 0  & \\ \hhline{|~|-|-|-|-|~|-|~|-|-|~|~|-|~|}
				& 3 & 5  & 7  & 15 & & 12 & & \cellcolor{gray!10}70508.0  & \cellcolor{gray!10}19.6  &     & & 0  & \\ \hhline{|~|-|-|-|-|~|-|~|-|-|-|~|-|~|}
				& 3 & 5  & 9  & 14 & & 4  & & \cellcolor{gray!50}446368.0 & \cellcolor{gray!50}124.0 & \cellcolor{gray!50}5.2 & & 8  & \\ \hhline{|~|-|-|-|-|~|-|~|-|-|-|~|-|~|}
				& 3 & 5  & 9  & 15 & & 4  & & \cellcolor{gray!50}754788.6 & \cellcolor{gray!50}209.7 & \cellcolor{gray!50}8.7 & & 22 & \\ \hhline{|~|-|-|-|-|~|-|~|-|-|-|~|-|~|}
				& 3 & 5  & 10 & 12 & & 3  & & 1.6      & 　    & 　  & & 0  & \\ \hhline{|~|-|-|-|-|~|-|~|-|~|~|~|-|~|}
				& 3 & 5  & 10 & 13 & & 24 & & 1.5      & 　    & 　  & & 0  & \\ \hhline{|~|-|-|-|-|~|-|~|-|-|-|~|-|~|}
				& 3 & 5  & 10 & 15 & & 12 & & \cellcolor{gray!50}117424.0 & \cellcolor{gray!50}32.6  & \cellcolor{gray!50}1.3 & & 0  & \\ \hhline{|~|-|-|-|-|~|-|~|-|-|-|~|-|~|}
				& 3 & 5  & 11 & 13 & & 12 & & 1.5      & 　    & 　  & & 0  & \\ \hhline{|~|-|-|-|-|~|-|~|-|-|~|~|-|~|}
				& 3 & 5  & 11 & 14 & & 24 & & \cellcolor{gray!10}38227.0  & \cellcolor{gray!10}10.6  &     & & 0  & \\ \hhline{|~|-|-|-|-|~|-|~|-|-|~|~|-|~|}
				& 3 & 5  & 11 & 15 & & 24 & & \cellcolor{gray!10}50797.0  & \cellcolor{gray!10}14.1  &     & & 0  & \\ \hhline{|~|-|-|-|-|~|-|~|-|-|-|~|-|~|}
				& 3 & 5  & 14 & 15 & & 12 & & \cellcolor{yellow!40}*        & \cellcolor{yellow!40}*     & \cellcolor{yellow!40}*   & & \cellcolor{yellow!40}*  & \\ \hhline{|~|-|-|-|-|~|-|~|-|-|-|~|-|~|}
				& 3 & 7  & 11 & 12 & & 6  & & \cellcolor{gray!10}1903.0   & \cellcolor{gray!10}0.5   & 　  & & 0  & \\ \hhline{|~|-|-|-|-|~|-|~|-|-|~|~|-|~|}
				& 3 & 7  & 11 & 13 & & 12 & & \cellcolor{gray!10}34232.0  & \cellcolor{gray!10}9.5   &     & & 0  & \\ \hhline{|~|-|-|-|-|~|-|~|-|-|~|~|-|~|}
				& 3 & 7  & 11 & 15 & & 6  & & \cellcolor{gray!10}9643.0   & \cellcolor{gray!10}2.7   &     & & 0  & \\ \hhline{|~|-|-|-|-|~|-|~|-|-|~|~|-|~|}
				& 3 & 7  & 12 & 13 & & 12 & & 3.1      & 　    &     & & 0  & \\ \hhline{|~|-|-|-|-|~|-|~|-|-|~|~|-|~|}
				& 3 & 7  & 12 & 15 & & 12 & & \cellcolor{gray!10}10023.0  & \cellcolor{gray!10}2.8   &     & & 0  & \\ \hhline{|~|-|-|-|-|~|-|~|-|-|-|~|-|~|}
				& 3 & 7  & 13 & 14 & & 12 & & \cellcolor{gray!50}136406.0 & \cellcolor{gray!50}37.9  & \cellcolor{gray!50}1.5 & & 0  & \\ \hhline{|~|-|-|-|-|~|-|~|-|-|-|~|-|~|}
				& 3 & 7  & 13 & 15 & & 24 & & \cellcolor{gray!10}36482.0  & \cellcolor{gray!10}10.1  & 　  & & 0  & \\ \hhline{|~|-|-|-|-|~|-|~|-|-|-|~|-|~|}
				& 3 & 13 & 14 & 15 & & 6  & & \cellcolor{gray!50}193959.0 & \cellcolor{gray!50}53.9  & \cellcolor{gray!50}2.2 & & 0  & \\ \hhline{|~|-|-|-|-|~|-|~|-|-|-|~|-|~|}
				& 7 & 11 & 13 & 14 & & 1  & & \cellcolor{yellow!40}*        & \cellcolor{yellow!40}*     & \cellcolor{yellow!40}*   & & \cellcolor{yellow!40}*  & \\ \hhline{|~|-|-|-|-|~|-|~|-|-|-|~|-|~|}
				& 7 & 11 & 13 & 15 & & 4  & & \cellcolor{yellow!40}*        & \cellcolor{yellow!40}*     & \cellcolor{yellow!40}*   & & \cellcolor{yellow!40}*  & \\
				\hline
			\end{tabular}
		}
	\end{table}
The coloring cases are computed in this work. The machine is equipped with Windows 7 Ultimate. Its processor is an AMD A8-7100 Redeon R5, with eight computing cores $4C+4G$ (4 CPUs) of a 1.8 GHz clockspeed, and the RAM is 16384 MB. Python is our programming language. It costs 11.069 seconds to recalculate the (15)-coloring case, namely the Garrison-Scott manifold, by this algorithm. We list this characteristic vector in Table \ref{table:5col}. It took nearly 11 hours to accomplish the whole 7-coloring case. However,  when we move on to address the 8-coloring case, the low efficiency and huge time consumption makes it unaffordable. See Table \ref{table:compare1} for the detailed data.

In the cases of $(3,5,9,14)$-coloring and $(3,5,9,15)$-coloring, some new manifolds are constructed. The * here represents the case that can not be calculated within an acceptable amount of time. It is estimated from experimental testing that at least half a year is needed to complete the $(7,11,13,14)$-coloring case. By Corollary \ref{corollary: NakayamaN}  and Proposition \ref{prop:4567}, all of the orientable $\mathbb{Z}_2^4$-colorings belong to the $(7,11,13,14)$-coloring case. Thus, we cannot effectively obtain all the orientable small covers over the right-angled 120-cell through this algorithm.
\subsection{A 9-coloring small cover over the 120-cell $\mathcal{E}$}\label{section:9col_cr}
We also try to construct new small covers over the right-angled 120-cell through another approach. Take the unique 5-coloring characteristic vector and color it after freeing the fixed colors from the back to the front. We begin by fixing the first 119 colors of the Garrison-Scott characteristic vector and using the color-recursion algorithm to run over all of the admissible colors on the 120-th facet. Later, we fix the colors of the first 118 facets and color the last two, etc. We conduct this procedure by computer and it is not until we move ahead to free the colors of the facets behind the first $33$ facets that we obtain a 9-coloring, which  gives a new small cover. The first 33 facets cover exactly the first three layers of the boundary of the facet-ordered 120-cell, as shown in Table \ref{table:layers}.

We stop here because these results required approximately 13 days and need to consume much more time to process further. Moreover, the 9-coloring characteristic function we obtain here is in the same homeomorphism class with the 9-coloring results obtained so far by the new block-pasting algorithm that is proposed in Section \ref{section:bp}. 

In total, we obtain three more new manifolds with the color-recursion algorithm up to homeomorphism. We explain the homeomorphism classification in Section \ref{section:classify}.

\section{Constructing small covers with the proposed block-pasting algorithm}\label{section:bp}
Now we change our strategy and propose a new algorithm. We firstly truncate the set of facets by blocks and produce the truncated characteristic vectors concurrently. Then, we match them up and finally obtain the characteristic vectors.
\subsection{Optimized algorithm method}
\begin{figure}[H]
	\scalebox{0.25}[0.25]{\includegraphics {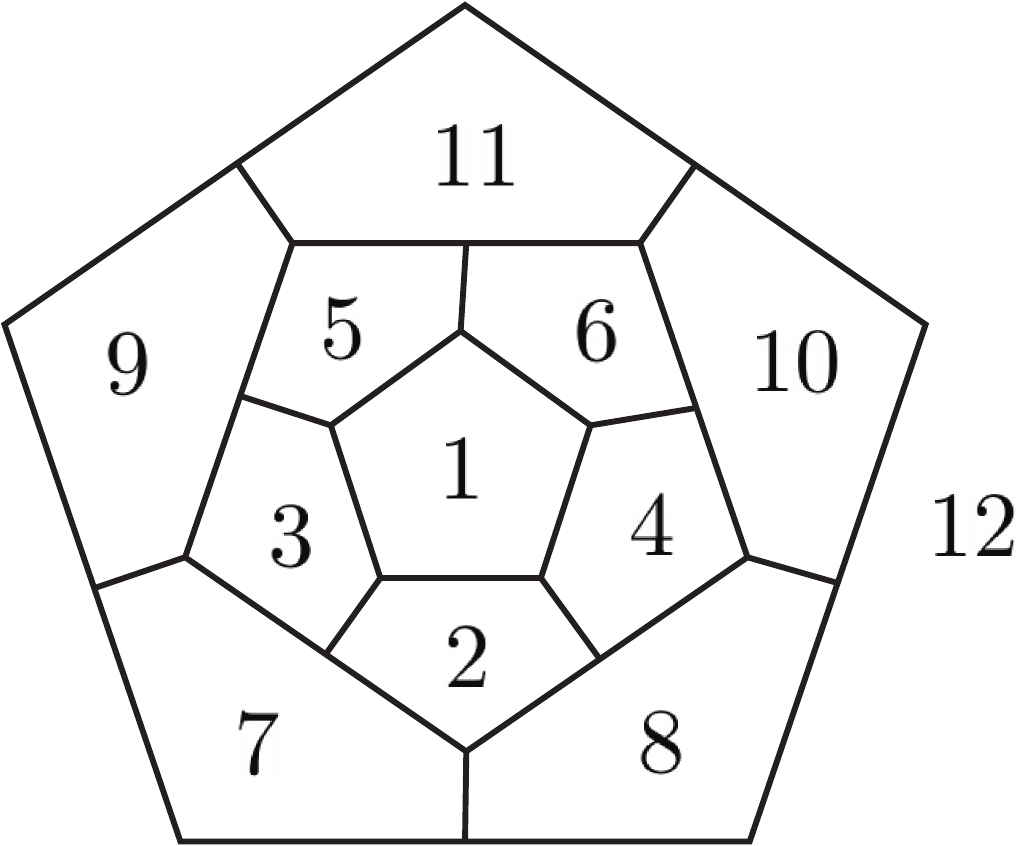}}
	\caption{Labeled facets of a dodecahedron (the numbers are the subscripts of the corresponding facets).}\label{figure:block1}
\end{figure}

\begin{figure}[H]
	\scalebox{0.25}[0.25]{\includegraphics {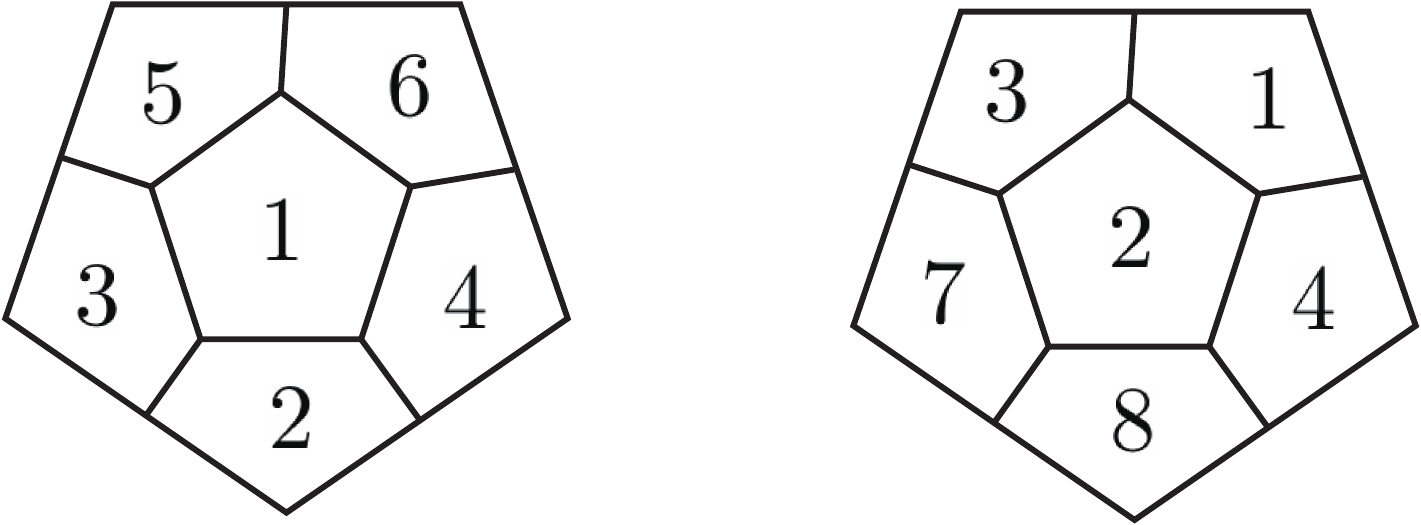}}
	\caption{Block of two boundary pentagons $F_1$ and $F_2$ (the numbers are the subscripts of the corresponding facets). }\label{figure:block2}
\end{figure}

The new algorithm named \textbf{block-pasting} is designed to address the computational efficiency issue  by making use of the high symmetry of the 120-cell.  For each  facet $F_i$ of an $n$-dimensional polytope $P^n$, we define the \textbf{i-th ~~block}, denoted as $b_i$, to be the ordered set of facets consisting of $F_i$ and all of the facets adjacent to it with increasing subscripts. We may also use $b_i$ to denote the ordered set of subscripts when $b_i$ is referred to as a set of integers. For example, for a facet-labeled dodecahedron $P$ in Figure \ref{figure:block1},

\noindent the block  $b_1=\{F_1,F_2,F_3,F_4,F_5,F_6\}$ (or $\{1,2,3,4,5,6\}$) and the block $b_2=\{F_1,F_2,F_3,F_4,F_7,F_8\}$ (or $\{1,2,3,4,5,6\}$) are as shown in Figure \ref{figure:block2}.

\noindent Recall that every dodecahedral facet of the 120-cell is adjacent to exactly 12 other dodecahedral facets. So, all of the 120 blocks of the 120-cell $\mathcal{E}$ contain 13 dodecahedral facets. 

Each symmetry of the 120-cell clearly induces a permutation of the facets. In Section \ref{section:symmetry of 120}, we describe the symmetry group $\mathbb{A}(\mathcal{E})$ of the 120-cell as a subgroup  of the permutation group $S_{120}$. For every  $1\le i\le 120$, there is a set made of $120$ induced permutations, any of which sends the facet $F_1$ to the facet $F_i$. We arbitrarily pick one from each of these 120 sets and  form a representative set $\mathcal{J}$ with 120 elements. Now, the whole block $b_1$ may be transferred to another block $b_i$ by some unique element in the set $\mathcal{J}$, for $1\le i\le 120$. The set of the selected 120 such permutations can be regarded as a $120\times 120$ matrix $M_{\mathcal{J}}$, where every row represents a permutation.

We are now only concerned with the results of these permutations  on the first block $b_1$. Hence, we only need the first $120\times 13$ submatrix $M^{13}_{\mathcal{J}}$ of the matrix $M_{\mathcal{J}}$, where every row represents a result of one permutation on $b_1$. Parts of the matrix $M^{13}_{\mathcal{J}}$ are shown in Table \ref{table:per}.

\begin{table}[H]
		\footnotesize
	\begin{tabular}{c|P{0.7cm}P{0.7cm}P{0.7cm}P{0.7cm}P{0.7cm}P{0.7cm}P{0.7cm}P{0.7cm}P{0.7cm}P{0.7cm}P{0.7cm}P{0.7cm}P{0.7cm}}
		\hline
		 & \textbf{1} & \textbf{2} & \textbf{3}& \textbf{4} & \textbf{5} & \textbf{6}& \textbf{7} & \textbf{8}&\textbf{9} & \textbf{10} & \textbf{11} & \textbf{12} & \textbf{13}
\\
\hline
\textbf{1}& 1 & 2 & 3 & 4 & 5 & 6 & 7 & 8 & 9 & 10 & 11 & 12 & 13
\\
\hline
\textbf{2} & 2 & 34 & 28 & 26 & 10 & 1 & 4 & 18 & 12 & 14 & 3 & 30 & 8
\\
\hline
\textbf{3} & 3 & 26 & 36 & 16 & 1 & 13 & 11 & 28 & 8 & 4 & 32 & 2 & 20
\\
\hline
\textbf{4} & 4 & 14 & 26 & 38 & 7 & 11 & 23 & 2 & 1 & 22 & 16 & 10 & 3
\\
\hline
\textbf{5} & 5 & 12 & 1 & 10 & 35 & 27 & 15 & 9 & 29 & 31 & 7 & 19 & 6
\\
\hline
\textbf{6}& 6 & 1 & 11 & 7 & 29 & 37 & 27 & 13 & 21 & 5 & 17 & 9 & 33
\\
\hline
\textbf{7} & 7 & 10 & 4 & 22 & 27 & 17 & 39 & 1 & 6 & 15 & 23 & 5 & 11
\\
\hline
\textbf{8} & 8 & 28 & 20 & 3 & 12 & 9 & 1 & 40 & 24 & 2 & 13 & 18 & 25
\\
\hline
\textbf{9} & 9 & 8 & 13 & 1 & 19 & 29 & 5 & 25 & 41 & 12 & 6 & 24 & 21
\\
\hline
\textbf{10} & 10 & 30 & 2 & 14 & 15 & 7 & 22 & 12 & 5 & 42 & 4 & 31 & 1
\\
\hline
\textbf{11} & 11 & 4 & 16 & 23 & 6 & 33 & 17 & 3 & 13 & 7 & 43 & 1 & 32
\\
\hline
\textbf{12} & 12 & 18 & 8 & 2 & 31 & 5 & 10 & 24 & 19 & 30 & 1 & 44 & 9
\\
\hline
\textbf{13} & 13 & 3 & 32 & 11 & 9 & 21 & 6 & 20 & 25 & 1 & 33 & 8 & 45
\\
\hline
$\vdots$ & $\vdots$ & $\vdots$ & $\vdots$ & $\vdots$ & $\vdots$ & $\vdots$ & $\vdots$ & $\vdots$ & $\vdots$ & $\vdots$ & $\vdots$ & $\vdots$ & $\vdots$
\\
 \hline        
\textbf{120} & 120  & 108 & 109 & 110 & 111 &114  & 115 & 112  & 113 & 118  & 119 & 116 & 117
\\
\hline
\end{tabular}
\caption{ The matrix $M^{13}_{\mathcal{J}}$. }\label{table:per}
\end{table}

We use $X_{\mathcal{E}(b_i)}$ to represent the $13\times 13$ submatrix of $X_{\mathcal{E}}$ obtained by selecting all of the $k$-th rows and the $k$-th columns, where $k \in b_i$. Given a coloring set $C=\{1,2,4,8,c_1,\cdots,c_k\}$. We use the color-recursion algorithm to construct the \emph{truncated characteristic vectors} of length $13$. We also need to ensure the non-singularity condition by checking simplexes as mentioned in Section \ref{section:9col_cr}. The only difference is that the set $\mathcal{A}^1(F_i)$ is now read from the truncated adjacency matrix $X_{\mathcal{E}(b_i)}$, rather than the whole adjacency matrix $X_{\mathcal{E}}$. The value of the $q$-th entry of a resulting characteristic vector is the color of the $q$-th facet in  $b_i$. 

For example, we obtain $$b_3=\{F_1,F_2,F_3,F_4,F_8,F_{11},F_{13},F_{16},F_{20},F_{26},F_{28},F_{32},F_{36}\}$$ from the adjacency matrix $X_{\mathcal{E}}$ by probing those facets that are adjacent to or equal to facet $F_3$. Then $X_{\mathcal{E}(b_3)}=(a_{ij})_{13\times 13},~\text{where}~ a_{ij}=a_{t_i,t_j}\in X_{\mathcal{E}},~ t_i~\text{and}~t_j \text{~are~ the}~ i\text{-th}~\text{and}~j\text{-th}$ elements in the ordered set $\{1,2,3,4,8,11,13,16,20,26,28,32,36\}.$ Now, we input the coloring set  $C=\{1,2,4,8,7,11,13,14\}$, the vector $v_{1\times 13}=(0,0,\cdots,0)$ and the truncated adjacency matrix $X_{\mathcal{E}(b_i)}$ into the color-recursion algorithm and run the function $\textsc{ColorRecursion}(1,v)$. Suppose $$v=(2,7,1,15,15,4,8,8,11,2,11,7,4)$$ is one of the resulting truncated coloring vector. It implies that the facets $$F_1,F_2,F_3,F_4,F_8,F_{11},F_{13},F_{16},F_{20},F_{26},F_{28},F_{32},F_{36}$$ can be colored by $2,7,1,15,15,4,8,8,11,2,11,7,4 $, respectively.

The set of all truncated coloring vectors with coloring set $\{1,2,4,8,c_5,\cdots,c_k\}$ is called the \emph{unfix set of $b_i$ defined on $(c_5,c_2,\cdots,c_k)$}, where $c_5,\cdots,c_k$ are the added colors. We may not emphasize the added colors if there is no ambiguity. Next, we denote the unfix set of $b_i$ by $B_i$. In particular, when the first four entries are fixed to be $e_1$, $e_2$, $e_3$, and $e_4$, we obtain  the \emph{fixing set of $b_i$}; it is denoted by $\overline{B}_i$. By the color-recursion algorithm, we can calculate the sets $B_1$ and $\overline{B}_1$ according to the adjacency information given by $X_{\mathcal{E}(b_1)}$.

The set $B_i$ (or $\overline{B}_i$) can be regarded as an $l\times 13$ matrix, where $l$ is the cardinality of $B_i$ (or $\overline{B}_i$). In the following, we do not distinguish these two viewpoints and may refer to $B_i$ (or $\overline{B}_i$) as either a set or a matrix. From Lemma \ref{lemma:S4}, every vertex of the 120-cell $\mathcal{E}$ enjoys the same stabilizer. Thus, in a certain $k$-coloring case, the unfix sets of all the other blocks $B_2,B_3,\cdots,B_{120}$ can be obtained from $B_1$ by permuting the columns of $B_1$. The permutation rule with respect to certain $B_i$, $2\le i\le120$, can be inferred from the adjacency matrix $X_{\mathcal{E}}$ and the set $\mathcal{J}$ formed by 120 induced permutations as defined before. 

 For example, the block $b_1$ and $b_2$ are $$(1,2,3,4,5,6,7,8,9,10,11,12) ~~\rm{and}~~ (1,2,3,4,8,10,12,14,18,26,28,30,34),$$ respectively. We select the permutation, from the set $\mathcal{J}$, that transfers $$(1,2,3,4,5,6,7,8,9,10,11,12,13)~~\rm{to}~~ (2,34,26,14,12,1,10,28,8,30,4,18,3).$$ By matching the positions, the unfix set of $b_2$ can be obtained by rearranging the column's order of $B_1$ as illustrated in Table \ref{table:match}. The rearranged order now is $(6,1,13,11,9,7,5,4,12,3,8,10,2)$. This means the first column of matrix $B_2$ corresponds to the $6$-th column of matrix $B_1$, and the second column of $B_2$ corresponds to the first column of matrix $B_1$, etc. Similarly, we can obtain all the orders of rearrangements  for $B_i$, $2\le i\le 120$.

\begin{table}[H]
	\footnotesize
	
	\begin{tabular}{P{0.7cm}P{0.7cm}P{0.7cm}P{0.7cm}P{0.7cm}P{0.7cm}P{0.7cm}P{0.7cm}P{0.7cm}P{0.7cm}P{0.7cm}P{0.7cm}P{0.7cm}}
		\hline
		\textbf{1} & \textbf{2} & \textbf{3}& \textbf{4}& \textbf{8} & \textbf{10} & \textbf{12}& \textbf{14} & \textbf{18}&\textbf{26} & \textbf{28} & \textbf{30}  & \textbf{34}\\
		\hline
		1&2&3&4&5&6&7&8&9&10&11&12&13\\
		\hline
		2&34&26&14&12&1&10&28&8&30&4&18&3\\
		\hline
	\end{tabular}
	\caption{ Finding the order for rearranging. }
	\label{table:match}
\end{table}

So far, we have prepared 120 sets $\overline{B}_1$, $B_2$, ..., $B_{120}$. Now, every $l\times 13$ matrix is extended to an $l\times 120$ matrix by simply putting each $s$-th column to the position of the $b_i(s)$-th column, where $b_i(s)$ means the $s$-th element in the set $b_i$, and filling in the value zero in the other positions.
We continue to use the same notation $B_i$ (or $\overline{B}_i$) for the extended matrix.

Now, we are going to paste the set $\overline{B}_1$ and $B_2$ based on some assigned columns of them. More precisely, a row from $\overline{B}_1$ is matched up with a row of $B_2$ where every two entries specified by the same index $i$, where $i\in b_1\cap b_2$, are with the same values. The index set $b_1\cap b_2$ is called a \emph{linking key} for the pasting. The resulting new row is actually the sum of these two rows. The set of all these new rows is denoted by $\overline{B}_1\cup^*B_2$.

\begin{figure}[H]
	\scalebox{0.25}[0.25]{\includegraphics {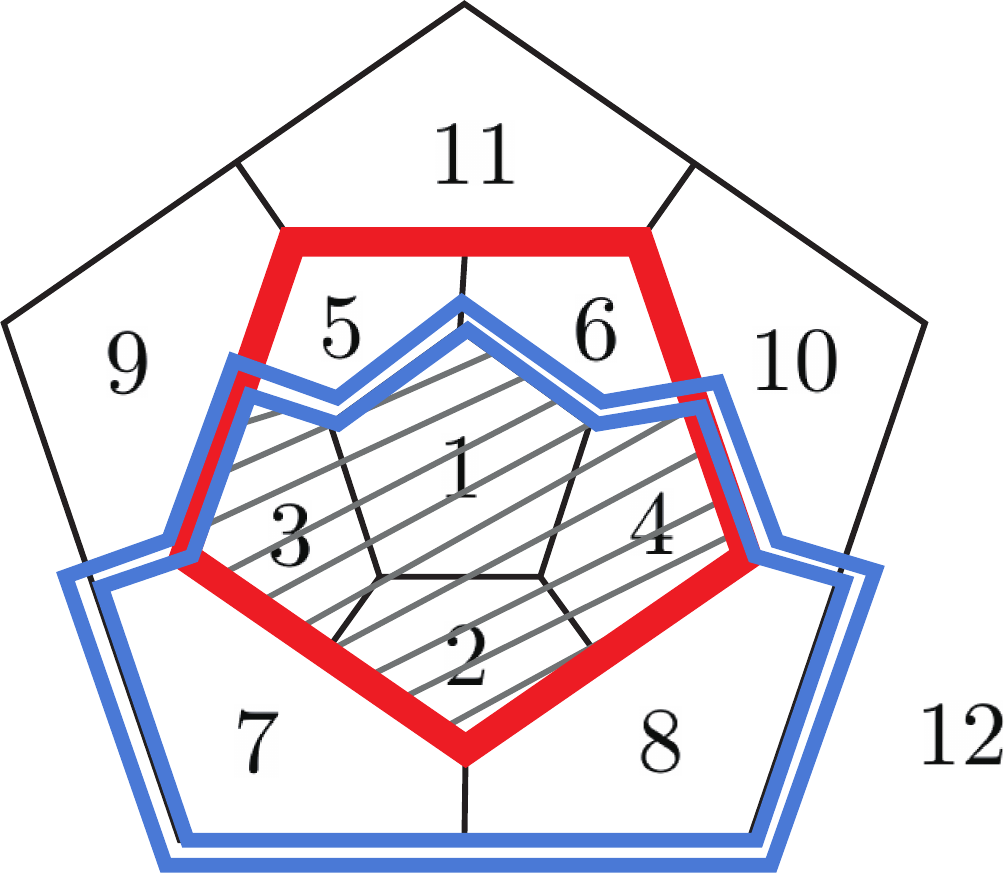}}
	\caption{Pasting two colored blocks of boundary pentagons (the numbers are the subscripts of the  corresponding facets).}\label{figure:paste}
\end{figure}

We use the dodecahedron $\mathcal{D}$ to explain this process. Figure \ref{figure:paste} is a projective illustration of the dodecahedron. The area circled by the bold line is $b_1$, the block of $F_1$; the area circled by the double line is $b_2$, the  block of $F_2$. We suppose, for example, that

$\overline{B}_1=\{x_1,x_2\}=\{(1,2,4,4,2,6,0,0,0,0,0,0), (1,2,4,5,2,6,0,0,0,0,0,0)\}$,

$B_2=\{y_1,y_2,y_3\}=\{(1,2,4,4,0,0,1,7,0,0,0,0),(1,2,4,4,0,0,6,5,0,0,0,0)$,

\hspace{3.5cm}$(1,2,3,4,0,0,1,7,0,0,0,0)\} .$

In this example, $x_1$ has the same color with $y_1$ and $y_2$ on the first, second, third and fourth facets. Thus, $y_1$ and $y_2$ can paste to $x_1$, forming the coloring vectors $(1,2,4,4,2,6,1,7,0,0,0,0)$ and $(1,2,4,4,2,6,6,5,0,0,0,0)$, respectively. In contrast, $x_2$ cannot be pasted to any element of $B_2$ as there are no vectors of the color $5$ on the 4-th position. Therefore, $$\overline{B}_1\cup^* B_2=\{(1,2,4,4,2,6,1,7,0,0,0,0),(1,2,4,4,2,6,6,5,0,0,0,0)\}.$$

We then move on to paste the sets $\overline{B}_1 \cup^*B_2$ and $B_3$. We follow the same procedure except that the index set that specifies the matching line, namely the linking key, is $b_1\cup b_2\cap b_3$. We conduct this procedure until we finish pasting the final set $B_{120}$ and obtain the resulting set  $S=\overline{B}_1\cup^*B_2\cdots \cup^*B_{120}$. The set of linking keys of this procedure is $$\{b_1\cap b_2, b_1\cup b_2\cap b_3,\cdots,b_1\cup b_2\cup\cdots\cup b_{i-1}\cap b_i,\cdots ,b_1\cup b_2\cup \cdots \cup b_{119}\cap b_{120}\}.$$
 
Note that the first set is a fixing set, which is different from other 119 sets. This means after pasting the final facet $F_{120}$, we actually obtain results of $GL_4(\mathbb{Z}_2)$-equivalence representatives. The pseudocode for the block-pasting algorithm is given in Algorithm 2. 

\begin{algorithm}
	\small
	\label{algorithm3}
	\caption{Pseudocode for the block-pasting algorithm}
	\begin{algorithmic}
		\Require
		\Statex $\mathcal{L}$: Set of Linking keys
		\Statex $T$: Selected 120 permutations that permute $F_1$ to the original position of $F_i$.
		\Statex $B_1$: Unfix set of $b_1$
		\Statex $\overline{B_1}$: Fixing set of $b_1$
		
		\Ensure All of the characteristic vectors that satisfy nonsingularity condition will be saved in an assigned txt file.
		
		\Function{BlockPaste}{$\mathcal{R},\mathcal{O},T,B_1,\overline{B_1}$}
		\State $B_2=t_2(B_1)$ \Comment{{\color{darkgray}obtain $B_2$ from $B_1$ through the second permutation $t_2$ in $T$}}
		\ForAll {$i$ in row($\overline{B_1}$)}
		\State $S\gets$ $i*_{l_1}B_2$ 
		
		\Comment{{\color{darkgray} here * represents merging by the first linking key $l_1:=b_1\cap b_2\in \mathcal{L}$}}
		
		\ForAll {$j$ in $[3,4]$}
		\State generate $B_j$ from $B_1$ through the $j$-th permutation in $T$
-		\State $S\gets S*_{j-1} B_j$
		\EndFor
		
		\If {$S=\emptyset$}
		\State \textbf{break}
		\EndIf
		
		\If {$\vert S\vert>120000$} \Comment{{\color{darkgray}120000 is set according to the internal storage of our machine}}
		\State {piece up the set $S$ into $t$ many subset such that the cardinality of every subset $S_k$, ($k=1,2,...,t$), is less than $120000$}
		
		\ForAll {subset $S_k$}
		\ForAll {$j$ in $[5,6,\cdots,120]$}
		\State generate $B_j$ and $S_k\gets S_k*_{l_{j-1}}B_j$
		
		\If {$S_k=\emptyset$}
		\State \textbf{break}
		\EndIf
		\EndFor		
		\EndFor
		
		\State \textbf{$S\gets\bigcup S_k$\quad and \quad save $S$ into an assigned ``i{\_}result.txt"}
		\EndIf
		
		\If {$\vert S\vert\le120000$}
		\ForAll {$j$ in $[5,6,\cdots,120]$}
		\State generate $B_j$ and $S\gets S*_{l_{j-1}}B_j$
		
		\If {$S=\emptyset$}
		\State \textbf{break}
		\EndIf
		\EndFor
		\State \textbf{save $S$ into an assigned ``i{\_}result.txt"}
		\EndIf

		\EndFor
		\State \textbf{Return null}
		\EndFunction

	\end{algorithmic}
\end{algorithm}

Compared to the first algorithm, the second one is far more efficient. The color-recursion algorithm is using the method of ``series circuit", where the resulting characteristic function is produced one by one. Whereas, the block-pasting algorithm adopts the idea of ``parallel circuit", where different parts of a resulting characteristic vector are generated simultaneously and then  pasted together. However, it demands huge internal storage because the step of gluing may cause a surge in data. By experimental testing, at least 2000 GB of internal storage is needed to fulfill this algorithm in the  $8$-coloring case. Hence, we start the block matching from one single element in $\overline{B_1}$ instead of dealing with the whole set at a time to reduce the storage requirements. Moreover, we control the amount of truncated coloring vectors that are being pasted simultaneously with the unfix set $B_5$ to make the computation more practical. More precisely, we set up an upper threshold of 120,000. This quantity is set in accordance with the device storage. This means the result set with the number of elements larger than 120,000 is automatically pieced into $p$ subsets with cardinalities all less than 120,000. The number $p$ is the smallest satisfied integer that is calculated by the computer automatically. 

\subsection{Verification and results}
We have five evidences to validate our computing results:
\begin{enumerate}
	\item Apply the color-recursion algorithm to compute the number of  $GL_n(\mathbb{Z}_2)$-equivalence classes of $\mathbb{Z}_2^n$-colorings of the $n$-cube $\emph{I}^n$ for $n=2,3,4,5,6$ as shown in Table \ref{table:choi}. The results coincide with the work of Suyong Choi \cite{suyong:08}.
	\begin{table}[H]
		\small
		\begin{tabular}{c|c|c|c|c}
			\hline
			$\emph{I}^2$ & $\emph{I}^3$ & $\emph{I}^4$ & $\emph{I}^5$& $\emph{I}^6$ \\		
			\hline
			3&25& 543&29,281&3,781,503\\
			\hline
		\end{tabular}
		\caption{The number of $GL_n(\mathbb{Z}_2)$-equivalence classes of  $\mathbb{Z}_2^n$-colorings of the $n$-cube $\emph{I}^n$}
		\label{table:choi}
	\end{table}
	\item Apply the color-recursion algorithm to compute the number of  $GL_3(\mathbb{Z}_2)$-equivalence classes of $\mathbb{Z}_2^3$-colorings  of the 3-dimensional L\"obell polytope $\mathcal{L}_{n}$ for $n=5,6,7$ as shown in Table \ref{table:fu}. The results coincide with the work of Xin Fu \cite{Xin:2017}.
	\begin{table}[H]
		\small
		\begin{tabular}{c|c|c}
			\hline
			$\mathcal{L}_5$ & $\mathcal{L}_6$ & $\mathcal{L}_7$  \\		
			\hline
			2,165&18,073& 97,454\\
			\hline
		\end{tabular}
		\caption{The number of $GL_3(\mathbb{Z}_2)$-equivalence classes of  $\mathbb{Z}_2^3$-colorings of L\"obell polytope $\mathcal{L}_{n}$}
		\label{table:fu}
	\end{table}
	\item Apply the color-recursion algorithm to compute the number of  $GL_4(\mathbb{Z}_2)$-equivalence classes of $\mathbb{Z}_2^4$-colorings  of the neighborly $4$-polytope with eight facets as shown in Table \ref{table:bm}. The results coincide with the results done by Barali\'c--Milenkovi\'c \cite{Baralic:2017}
		\begin{table}[H]
			\small
			\begin{tabular}{c|c}
				\hline
				$P_0^4(8)$ & $P_1^4(8)$  \\		
				\hline
				7&3\\
				\hline
			\end{tabular}
			\caption{The number of $GL_4(\mathbb{Z}_2)$-equivalence classes of $\mathbb{Z}_2^4$-colorings of neighborly $4$-polytope with eight facets}
			\label{table:bm}
		\end{table}

\item We cross compare the results of the color-recursion and the block-pasting algorithms for the 120-cell. Both of them produce 10 characteristic vectors in the setting of 5-coloring and no result in the $\{4,6,7\}$-coloring cases. Then, by comparing Table \ref{table:compare1} and Table \ref{table:compare2} we can see these two algorithms generate the same resulting characteristic vectors in the cases of $(3,5,9,14)$-coloring and $(3,5,9,15)$-coloring.

\item Although we haven't ruled out all of the characteristic vector of 9-coloring setting yet, the two algorithms both calculate a common 9-coloring small cover as reported in Table \ref{table:9col}.
\end{enumerate}

Through the block-pasting algorithm, we finally solve the three remaining cases as shown in Figure \ref{table:compare2}. There are $480\times 128=61,440$ orientable $GL_4(\mathbb{Z}_2)$-equivalence small covers over the right-angled 120-cell with respect to the $(7,11,13,14)$-coloring case. The multiplicative equation means that there are 480 truncated coloring vectors in $B_1$ defined on $(7,11,13,14)$ that can be finally pasted up; each of them equally has 128 resulting vectors.  Moreover, we can now easily compare the efficiency and productivity of the color-recursion and the proposed block-pasting algorithms through Table \ref{table:compare1} and Table \ref{table:compare2}. The latter can  be up to 60 times faster in some cases. The $(3,5,9,14)$-coloring case costs two hours with the block-pasting algorithm; it takes five days with the color-recursion algorithm. The performance improvement is mainly because after calculating $B_1$, the time taken to generating $B_2$, $B_3$,..., $B_{120}$ through permutations in the block-pasting algorithm is negligible.

	\begin{table}[H]
		\caption{Block-pasting algorithm computation of 8-coloring small covers over the 120-cell.}
		\label{table:compare2}
		\scalebox{0.7}{
			\linespread{1}\selectfont
			\begin{tabular}{|p{1em}|c|c|c|c|p{1em}|c|p{1em}|c|p{1em}|c|p{1em}|c|c|c|p{1em}|c|p{1em}|}  
				\hline
				& \multicolumn{4}{c|}{\multirowcell{3}{Added colors}}                                                            & & \multirow{3}{5em}{\centering Cardinality of the class} & & \multirowcell{3}{$\vert B_1 \vert$} & & \multirowcell{3}{$\vert \overline{B_1} \vert$} & & \multicolumn{3}{c|}{\multirowcell{2}{Running time of \\ Color-Recursion Algorithm}}  & & \multirow{3}{8em}{\centering Number of resulting vectors} & \\
				& \multicolumn{4}{c|}{}                                                                                          & &                                                        & &                                     & &                                                & & \multicolumn{3}{c|}{}                                                                & &                                                           & \\ \cline{13-15}
				& \multicolumn{1}{p{2em}}{} & \multicolumn{1}{p{2em}}{} & \multicolumn{1}{p{2em}}{} & \multicolumn{1}{p{2em}|}{} & & \multicolumn{1}{p{8em}|}{} & & \multicolumn{1}{p{4em}|}{} & & \multicolumn{1}{p{4em}|}{} & & \multicolumn{1}{p{4em}|}{\centering Seconds} & \multicolumn{1}{p{4em}|}{\centering Hours} & Days & &                                                           & \\
				\hhline{|~|-|-|-|-|~|-|~|-|~|-|~|-|-|-|~|-|~|}
				& 3 & 5 & 6 & 7 & & 4 & & 0 & & & & 　 & 　 & 　 & & 　 & \\ \hhline{|~|-|-|-|-|~|-|~|-|~|~|~|~|~|~|~|~|~|}
				& 3 & 5 & 6 & 9 & & 12 & & 0 & & & & & & & & & \\ \hhline{|~|-|-|-|-|~|-|~|-|~|~|~|~|~|~|~|~|~|}
				& 3 & 5 & 6 & 11 & & 12 & & 0 & & & & & & & & & \\ \hhline{|~|-|-|-|-|~|-|~|-|~|-|~|~|~|~|~|~|~|}
				& 3 & 5 & 6 & 15 & & 4 & & 160 & & 819600 & & & & & & & \\ \hhline{|~|-|-|-|-|~|-|~|-|~|-|~|~|~|~|~|~|~|}
				& 3 & 5 & 7 & 9 & & 12 & & 0 & & 　 & & & & & & & \\ \hhline{|~|-|-|-|-|~|-|~|-|~|~|~|~|~|~|~|~|~|}
				& 3 & 5 & 7 & 10 & & 24 & & 0 & & & & & & & & & \\ \hhline{|~|-|-|-|-|~|-|~|-|~|~|~|~|~|~|~|~|~|}
				& 3 & 5 & 7 & 11 & & 24 & & 0 & & & & & & & & & \\ \hhline{|~|-|-|-|-|~|-|~|-|~|-|~|~|~|~|~|~|~|}
				& 3 & 5 & 7 & 14 & & 12 & & 160 & & 819600 & & & & & & & \\ \hhline{|~|-|-|-|-|~|-|~|-|~|-|~|~|~|~|~|~|~|}
				& 3 & 5 & 7 & 15 & & 12 & & 160 & & 93600 & & & & & & & \\  \hhline{|~|-|-|-|-|~|-|~|-|~|-|~|-|-|-|~|-|~|}
				& 3 & 5 & 9 & 14 & & 4 & & 640 & & 698280 & & \cellcolor{gray!50}7534.0 & \cellcolor{gray!50}2.1 & & & \cellcolor{gray!50}8 & \\ \hhline{|~|-|-|-|-|~|-|~|-|~|-|~|-|-|-|~|-|~|}
				& 3 & 5 & 9 & 15 & & 4 & & 640 & & 698280 & & \cellcolor{gray!50}42828.7 & \cellcolor{gray!50}11.9 & & & \cellcolor{gray!50}22 & \\ \hhline{|~|-|-|-|-|~|-|~|-|~|-|~|-|-|-|~|-|~|}
				& 3 & 5 & 10 & 12 & & 3 & & 0 & & 　 & & 　 & 　 & & & 　 & \\ \hhline{|~|-|-|-|-|~|-|~|-|~|~|~|~|~|~|~|~|~|}
				& 3 & 5 & 10 & 13 & & 24 & & 0 & & & & & & & & & \\ \hhline{|~|-|-|-|-|~|-|~|-|~|-|~|~|~|~|~|~|~|}
				& 3 & 5 & 10 & 15 & & 12 & & 160 & & 202560 & & & & & & & \\ \hhline{|~|-|-|-|-|~|-|~|-|~|-|~|~|~|~|~|~|~|}
				& 3 & 5 & 11 & 13 & & 12 & & 0 & & 　 & & & & & & & \\ \hhline{|~|-|-|-|-|~|-|~|-|~|-|~|~|~|~|~|~|~|}
				& 3 & 5 & 11 & 14 & & 24 & & 160 & & 202560 & & & & & & & \\ \hhline{|~|-|-|-|-|~|-|~|-|~|-|~|~|~|~|~|~|~|}
				& 3 & 5 & 11 & 15 & & 24 & & 160 & & 202560 & & & & & & & \\ \hhline{|~|-|-|-|-|~|-|~|-|~|-|~|-|-|-|~|-|~|}
				& 3 & 5 & 14 & 15 & & 12 & & 1280 & & 698280 & & \cellcolor{yellow!40}29771.0 & \cellcolor{yellow!40}8.3 & & & \cellcolor{yellow!40}30 & \\ \hhline{|~|-|-|-|-|~|-|~|-|~|-|~|-|-|-|~|-|~|}
				& 3 & 7 & 11 & 12 & & 6 & & 112 & & 93600 & & 　 & 　 & & & 　 & \\ \hhline{|~|-|-|-|-|~|-|~|-|~|-|~|~|~|~|~|~|~|}
				& 3 & 7 & 11 & 13 & & 12 & & 112 & & 698280 & & & & & & & \\ \hhline{|~|-|-|-|-|~|-|~|-|~|-|~|~|~|~|~|~|~|}
				& 3 & 7 & 11 & 15 & & 6 & & 236 & & 93600 & & & & & & & \\ \hhline{|~|-|-|-|-|~|-|~|-|~|-|~|~|~|~|~|~|~|}
				& 3 & 7 & 12 & 13 & & 12 & & 0 & & 　 & & & & & & & \\ \hhline{|~|-|-|-|-|~|-|~|-|~|-|~|~|~|~|~|~|~|}
				& 3 & 7 & 12 & 15 & & 12 & & 256 & & 93600 & & & & & & & \\ \hhline{|~|-|-|-|-|~|-|~|-|~|-|~|~|~|~|~|~|~|}
				& 3 & 7 & 13 & 14 & & 12 & & 256 & & 698280 & & & & & & & \\ \hhline{|~|-|-|-|-|~|-|~|-|~|-|~|~|~|~|~|~|~|}
				& 3 & 7 & 13 & 15 & & 24 & & 256 & & 202560 & & & & & & & \\ \hhline{|~|-|-|-|-|~|-|~|-|~|-|~|~|~|~|~|~|~|}
				& 3 & 13 & 14 & 15 & & 6 & & 320 & & 819600 & & & & & & & \\ \hhline{|~|-|-|-|-|~|-|~|-|~|-|~|-|-|-|~|-|~|}
				& 7 & 11 & 13 & 14 & & 1 & & 2165 & & 2909760 & & \cellcolor{yellow!40}691201.2 & \cellcolor{yellow!40}192.0 & 8\cellcolor{yellow!40}.0 & & \cellcolor{yellow!40}480*128=61,440 & \\ \hhline{|~|-|-|-|-|~|-|~|-|~|-|~|-|-|-|~|-|~|}
				& 7 & 11 & 13 & 15 & & 4 & & 2165 & & 698280 & & \cellcolor{yellow!40}129606.0 & \cellcolor{yellow!40}36.0 & \cellcolor{yellow!40}1.5 & & \cellcolor{yellow!40}30 & \\
				\hline
			\end{tabular}
		}
	\end{table}

\section{Small covers over the 120-cell up to homeomorphism}\label{section:classify}
\subsection{DJ-equivalence rigidity}
Two $G$-manifolds $M_1$ and $M_2$ are  \emph{DJ-equivalent} if there is a homeomorphism  $f: M_1 \rightarrow M_2$ and an automorphism $\theta$ of $G$, such that $f(gx)=\theta (g)f(x)$,  for all $g\in G$ and
$x \in M_1$.

For an $n$-dimensional simple polytope $P$, let $\mathbb{A}=\mathbb{A}(P)$ denote the symmetry group of $P$ and $G=GL_n(\mathbb{Z}_{2})$ denote the general linear group of degree $n$ over $\mathbb{Z}_{2}$. Two small covers $M(P,\lambda)$ and $M(P,\mu)$ are \emph{DJ-equivalent} if and only if $\lambda=g \circ \mu \circ a_*$, where $g \in GL_n(\mathbb{Z}_{2})$ and $a_*$ is the permutation of facets induced by some some $a \in \mathbb{A}$ \cite{Scott:02}.

For example, denote the characteristic functions in (2), (3) and (4) of Figure \ref{figure:gl_cube} by $\lambda_1$, $\lambda_2$, and $\lambda_3$, respectively. The manifolds $M(P,\lambda_2)$ and $M(P,\lambda_3)$ can be proved to be DJ-equivalent for $\lambda_3=g\circ\lambda_2\circ a_*$, where $g$ is a linear transformation sending the basis $\{e_1,e_2\}$ to $\{e_2,e_1+e_2\}$ and $a_*$ is the permutation of facets induced by a  $\pi/2$ planar rotation $a\in\mathbb{A}(P)$. However $M(P,\lambda_1)$ is different from $M(P,\lambda_2)$ and $M(P,\lambda_3)$ up to DJ-equivalence.

The DJ-equivalence may coincide with the homeomorphism. By the Mostow Rigidity Theorem and the geometric feature of the hyperbolic right-angled 120-cell $\mathcal{E}$, Garrison and Scott have proven the following DJ-equivalence rigidity theorem.

\begin{theorem}\label{theorem: Scott} \emph{(Garrison-Scott \cite{Scott:02})}
	Two small covers over the right-angled 120-cell are homeomorphic if and only if they are DJ-equivalent.
\end{theorem}

\begin{lemma} \label{groupby}
	Two hyperbolic small covers $M(\mathcal{E},\lambda_i)$ and $M(\mathcal{E},\lambda_j)$ over the right-angled 120-cell are homeomorphic if and only if $(\lambda_i\cdot\mathbb{A}(\mathcal{E}))\bigcap (GL_4(\mathbb{Z}_2)\cdot\lambda_j)\ne\emptyset$. Here $\lambda_i\cdot\mathbb{A}(\mathcal{E})$ and $ GL_4(\mathbb{Z}_2)\cdot\lambda_j$ are the orbit spaces of $\lambda_i$ and $\lambda_j$ under the group actions of  $\mathbb{A}(\mathcal{E})$ and $GL_4(\mathbb{Z}_2)$, respectively.
	
\end{lemma}
Although the proof of Lemma \ref{groupby} is quite trivial, this lemma is of great importance in designing a practical grouping algorithm.
The problem of the classification up to homeomorphism has been transferred into detecting the intersections of some orbit spaces of  different characteristic functions.
We aim at finding the non-empty intersection of two small sets by Lemma \ref{groupby} instead of checking the belonging of a coloring to a big set by Theorem \ref{theorem: Scott} in order to reduce the amounts of both storage and computation. Namely we use $(\lambda_i\cdot\mathbb{A}(\mathcal{E}))\bigcap (GL_4(\mathbb{Z}_2)\cdot\lambda_j)\ne\emptyset$ rather than $\lambda_1\in \mathbb{A}(\mathcal{E})\cdot\lambda_2\cdot GL_4(\mathbb{Z}_2)$ as the basic rule in designing the algorithm. For an arbitrary characteristic function $\lambda_i$, the cardinalities of $\mathbb{A}(\mathcal{E})\cdot\lambda_i$ and $GL_4(\mathbb{Z}_2)\cdot\lambda_j$ are 14,400 and 20,160 respectively, which are much smaller than that of the large-scale multiplicative case with  290,340,000 items before duplication. This reduced  level of computational effort is much more feasible and efficient in actual experiments.

\subsection{DJ-equivalence classes of small covers over the $n$-cubes $I^n$}We firstly use the grouping algorithm to calculate the number  of DJ-equivalence classes over the $n$-cube as an experimental test before addressing the 120-cell case as shown in Table \ref{table:cube}.

By the color-recursion algorithm, we can calculate that the number of $GL_n(\mathbb{Z}_2)$-equivalence classes of small covers over the $n$-cubes $I^n$ are 25, 543, and 29281 where, $n$=3, 4, and 5, respectively. These results coincide with the number of acyclic digraphs with $n$ labeled nodes, denoted by $R_n$, which are counted by Robinson and Stanley in \cite{Robinson:70,Stabley:73}. Based on the work of Mckay et al. \cite{mckay:04}, the number of $GL_n(\mathbb{Z}_2)$-equivalence classes over $I^n$ equals $R_n$; see Section 2.1 of \cite{suyong:08} for more detail. 
The diffeomorphism classes of small covers over an $n$-cube is given by Choi, Masuda, and Oum \cite{choi_masuda:2017}. There is no direct previous result about the number DJ-equivalence classes. However, It is at least reasonable that our result on the number of the DJ-equivalence classes is bounded by the cardinalities of the results presented in \cite{mckay:04,choi_masuda:2017}.

\begin{table}[H]
	\small
	\begin{tabular}{c|P{1.5cm}|P{1.5cm}|P{1.5cm}|P{1.5cm}}
		\hline
		& $I^2$ & $I^3$ & $I^4$ & ...   \\		
		\hline
		$\# GL_n(\mathbb{Z}_2)$-equivalence classes &3& 25& 543&  ...\\
		\hline
		$\#$ DJ-equivalence classes &2& 5& 19&  ...\\
		\hline
		$\#$ diffeomorphism classes &2& 4& 12&  ...\\
		\hline
	\end{tabular}
	\caption{Numbers of classes of $GL_n(\mathbb{Z}_2)$-equivalence,  DJ-equivalence and diffeomorphism  classes of small cover over the $n$-cube $I^n$. }
	\label{table:cube}
\end{table}

\subsection{Homeomorphism classes of 8-coloring small covers over the 120-cell $\mathcal{E}$}

By Lemma \ref{groupby}, we can calculate the homeomorphism classes of the $GL_4(\mathbb{Z}_2)$-equivalence representatives that are obtained by the proposed  block-pasting algorithm. We prove here the following two theorems; the characteristic vectors are reported afterwards.

\noindent \textbf{Theorem \ref{theorem:orientable}.} \emph{There are exactly 56 orientable small covers over the right-angled hyperbolic 120-cell up to homeomorphism.}

We apply the Grouping algorithm to classify the 61,440 $GL_4(\mathbb{Z}_2)$-equivalence classes of the $(7,11,13,14)$-coloring by DJ-equivalence and obtain 56 final results. They are all orientable according to Corollary \ref{corollary: NakayamaN}. The characteristic vectors are all of size of $1\times 120$ and the 2-digits denoting forms of the colors are replaced by letters successively, namely 10, 11, 12, 13, 14, and 15 are denoted by $a,~~b,~~c,~~d,~~e,$ and $f$, respectively. For example, the first result in Table \ref{table:ori} actually means the characteristic vector below:

\begin{center}
	{\footnotesize
		\noindent$(1,2,4,8,2,8,4,11,14,14,7,7,13,13,1,14,14,13,4,7,2,7,2,8,4,11,13,1,11,4,8,2,11,7,7,13,4,1,11,2$,
		
		\noindent$13,11,1,1,1,13,2,8,14,8,7,2,8,14,4,8,1,13,8,13,11,4,14,4,1,14,8,14,2,8,14,11,14,11,7,2,11,1,1,11$,
		
		\noindent$4,13,7,7,7,7,13,4,1,2,11,8,8,7,11,2,14,4,1,13,7,11,13,13,4,14,2,14,2,4,4,8,13,2,14,11,1,1,8,7)$。}
\end{center}

\noindent All of the 56 characteristic vectors are reported in Table \ref{table:ori}. The readers can obtain the closed orientable hyperbolic $4$-manifolds from the characteristic vectors  and the combinatorics of the 120-cell given in Tables \ref{table:coor1}--\ref{table:coor3}.

{\centering
{\footnotesize	
	\renewcommand\baselinestretch{1.5}\selectfont
	\begin{longtable}{cccc}
		
		\captionsetup{font=normalsize}
        \caption{All of the orientable small covers over the right-angled hyperbolic 120-cell $\mathcal{E}$ up to homeomorphism.}\label{table:ori}\\
		\hline 
		&1248284bee77dd1eed4727284bd1b482b77d41b2 &&	1248284bee77dd1eed4727284bd1b482b77d41b2\\
		
		1&db111d28e8728e481d8db4e41e8e28ebeb72b11b &2&	db111d2e8e72e8481d8db48418ee2e8beb72b11b\\
		
		&4d7777d412b887b2e41d7bdd4e2e2448d2eb1187&&	4d7777d412b887b2e41d7bdd4e2e2448d2eb1187\\

		\hline  
		
		&1248284bee77dd1eed4727284bb1d482b77d41d2 &&	1248284bee77dd1eed4727284bb1d482b77d41d2\\
		
		3&bb111b28e8728e481d8dd4e41e8e28ebeb72b11d &4&	bb111b2e8e72e8481d8dd48418ee2e8beb72b11d\\
		
		&4b7777d412b887b2e41b7ddd4e2e2448d2eb1187&&	4b7777d412b887b2e41b7ddd4e2e2448d2eb1187\\
		
		\hline 

		&1248284bee77dd1eed4747282bd1b482b77d21b4 &&	1248284bee77dd1eed4747282bd1b482b77d21b4\\
		
		5&db111d28e8748e481d8db4e21e8e28ebeb74b11b &6&	db111d2e8e74e8481d8db48218ee2e8beb74b11b\\
		
		&2d7777d212b887b4e41d7bdd4e2e2448d2eb1187&&	2d7777d212b887b4e41d7bdd4e2e2448d2eb1187\\

		\hline 
		
		&1248284bee77dd1eed4747282bd1b482b77d21b4 &&	1248284bee77dd1eed4747282bd1b482b77d21b4\\
		
		7&db111d28e8748e481d8db4e21e8e28ebeb74b11b &8&	db111d2e8e74e8481d8db48218ee2e8beb74b11b\\
		
		&2d7777d212b887b4e41d7bdd4e2e2448d2eb1187&&	2d7777d212b887b4e41d7bdd4e2e2448d2eb1187\\
		
		\hline 
		
		&1248284bee77dd1eed4e17287bd8b482be7d71b1 &&	1248284bee77dd1eed4e17287bb8d482be7d71d1\\
		
		9&db118d24247142488d8db427124e242bebe1b18b &10&	bb118b22427124488d8dd447142e224bebe1b18d\\
		
		&7d7e77d782b887b1e41debdd4e2e2448d2eb1187&&	7b7e77d782b887b1e41beddd4e2e2448d2eb1187\\
		
		\hline 
		
		&1248284bee77dd1eed4e17287bb8d482be7d71d1 &&	1248284bee77dd1eed4e77281bd8b482be7d11b7\\
		
		11&bb118b24247142488d8dd427124e242bebe1b18d &12&	db118d22427724488d8db441142e224bebe7b18b\\
		
		&7b7e77d782b887b1e41beddd4e2e2448d2eb1187&&	1d7e77d182b887b7e41debdd4e2e2448d2eb1187\\
		\hline 
		
		&1248284bee77dd1eed4e77281bd8b482be7d11b7 &&	1248284bee77dd1eed4e77281bb8d482be7d11d7\\
		
		13&db118d24247742488d8db421124e242bebe7b18b &14&	bb118b22427724488d8dd441142e224bebe7b18d\\
		
		&1d7e77d182b887b7e41debdd4e2e2448d2eb1187&&	1b7e77d182b887b7e41beddd4e2e2448d2eb1187\\
		
		\hline 
		
		&1248284bee77dd1eed4127284bd7b482b17d41b2 &&	1248284bee77dd1eed4127284bb7d482b17d41d2\\
		
		15&db117d2e8e72e8487d8db48418ee2e8beb12b17b &16&	bb117b28e8728e487d8dd4e41e8e28ebeb12b17d\\
		
		&4d7177d472b887b2e41d1bdd4e2e2448d2eb1187&&	4b7177d472b887b2e41b1ddd4e2e2448d2eb1187\\
		
		\hline 
		
		&1248284bee77dd1eed4127284bb7d482b17d41d2 &&	1248284bee77dd1eed4147282bd7b482b17d21b4\\
		
		17&bb117b2e8e72e8487d8dd48418ee2e8beb12b17d &18&	db117d28e8748e487d8db4e21e8e28ebeb14b17b\\
		
		&4b7177d472b887b2e41b1ddd4e2e2448d2eb1187&&	2d7177d272b887b4e41d1bdd4e2e2448d2eb1187\\
		
		\hline 
		&1248284bee77dd1eed4147282bd7b482b17d21b4&
	 &	1248284bee77dd1eed4147282bb7d482b17d21d4\\
		
		19&db117d2e8e74e8487d8db48218ee2e8beb14b17b &20&	bb117b28e8748e487d8dd4e21e8e28ebeb14b17d\\
		
		&2d7177d272b887b4e41d1bdd4e2e2448d2eb1187&&	2b7177d272b887b4e41b1ddd4e2e2448d2eb1187\\
		
		\hline 
		&1248284bee77dd1eed4147282bb7d482b17d21d4&
	&	1248284bee77dd1eed4817287bdeb482b87d71b1\\
		
		21&bb117b2e8e74e8487d8dd48218ee2e8beb14b17d &22&	db11ed2242712448ed8db447142e224beb81b1eb\\
		
		&2b7177d272b887b4e41b1ddd4e2e2448d2eb1187&&	7d7877d7e2b887b1e41d8bdd4e2e2448d2eb1187\\
		
		\hline 
		
		&1248284bee77dd1eed4817287bbed482b87d71d1 &&	1248284bee77dd1eed4817287bbed482b87d71d1\\
		
		23&bb11eb2242712448ed8dd447142e224beb81b1ed &24&	bb11eb2424714248ed8dd427124e242beb81b1ed\\
		
		&7b7877d7e2b887b1e41b8ddd4e2e2448d2eb1187&&	7b7877d7e2b887b1e41b8ddd4e2e2448d2eb1187\\
		
		\hline 
		
		&1248284bee77dd1eed4877281bdeb482b87d11b7 &&	1248284bee77dd1eed4877281bbed482b87d11d7\\
		
		25&db11ed2242772448ed8db441142e224beb87b1eb &26&	bb11eb2242772448ed8dd441142e224beb87b1ed\\
		
		&1d7877d1e2b887b7e41d8bdd4e2e2448d2eb1187&&	1b7877d1e2b887b7e41b8ddd4e2e2448d2eb1187\\
		
		\hline 
		
		&1248284bee77dd7eed4721284bd1b482b71d47b2&&	1248284bee77dd7eed4721284bd1b482b71d47b2\\
		
		27&db111d28e8128e481d8db4e47e8e28ebeb72b71b &28&	db111d2e8e12e8481d8db48478ee2e8beb72b71b\\
		
		&4d1777d412b881b2e47d7bdd4e2e2448d2eb1187&&	4d1777d412b881b2e47d7bdd4e2e2448d2eb1187\\
		
		\hline 
		
		&1248284bee77dd7eed4721284bb1d482b71d47d2&&	1248284bee77dd7eed4e11287bd8b482be1d77b1\\
		
		29&bb111b2e8e12e8481d8dd48478ee2e8beb72b71d &30&	db118d22421124488d8db447742e224bebe1b78b\\
		
		&4b1777d412b881b2e47b7ddd4e2e2448d2eb1187&&	7d1e77d782b881b1e47debdd4e2e2448d2eb1187\\

		\hline 
		
		&1248284bee77dd7eed4e11287bd8b482be1d77b1&&	1248284bee77dd7eed4e11287bb8d482be1d77d1\\
		
		31&db118d24241142488d8db427724e242bebe1b78b &32&	bb118b22421124488d8dd447742e224bebe1b78d\\
		
		&7d1e77d782b881b1e47debdd4e2e2448d2eb1187&&	7b1e77d782b881b1e47beddd4e2e2448d2eb1187\\
		
		\hline 
		
		&1248284bee77dd7eed4e71281bb8d482be1d17d7&&	1248284bee77dd7eed4121284bd7b482b11d47b2\\
		
		33&bb118b22421724488d8dd441742e224bebe7b78d &34&	db117d2e8e12e8487d8db48478ee2e8beb12b77b\\
		
		&1b1e77d182b881b7e47beddd4e2e2448d2eb1187&&	4d1177d472b881b2e47d1bdd4e2e2448d2eb1187\\
		
		\hline 
		
		&1248284bee77dd7eed4121284bb7d482b11d47d2&&	1248284bee77dd7eed4141282bd7b482b11d27b4\\
		
		35&bb117b2e8e12e8487d8dd48478ee2e8beb12b77d &36&	db117d28e8148e487d8db4e27e8e28ebeb14b77b\\
		
		&4b1177d472b881b2e47b1ddd4e2e2448d2eb1187&&	2d1177d272b881b4e47d1bdd4e2e2448d2eb1187\\
		
		\hline 
		
		&1248284bee77dd7eed4141282bd7b482b11d27b4&&	1248284bee77dd7eed4141282bb7d482b11d27d4\\
		
		37&db117d2e8e14e8487d8db48278ee2e8beb14b77b &38&	bb117b28e8148e487d8dd4e27e8e28ebeb14b77d\\
		
    	&2d1177d272b881b4e47d1bdd4e2e2448d2eb1187&&	2b1177d272b881b4e47b1ddd4e2e2448d2eb1187\\
    	
		\hline 
		
		&1248284bee77dd7eed4141282bb7d482b11d27d4&&	1248284bee77dd7eed4811287bbed482b81d77d1\\
		
		39&bb117b2e8e14e8487d8dd48278ee2e8beb14b77d &40&	bb11eb2424114248ed8dd427724e242beb81b7ed\\
		
		&2b1177d272b881b4e47b1ddd4e2e2448d2eb1187&&	7b1877d7e2b881b1e47b8ddd4e2e2448d2eb1187\\
		
		\hline 
		
		&1248284bee77dd7eed4871281bdeb482b81d17b7&&	1248284bee77dd7eed4871281bdeb482b81d17b7\\
		
		41&db11ed2242172448ed8db441742e224beb87b7eb &42&	db11ed2424174248ed8db421724e242beb87b7eb\\
		
		&1d1877d1e2b881b7e47d8bdd4e2e2448d2eb1187&&	1d1877d1e2b881b7e47d8bdd4e2e2448d2eb1187\\
		
		\hline 
		
		&1248284bee77dd7eed4871281bbed482b81d17d7&&	1248284bee77db1eed4727284dd1b482b77b41b2\\
		
		43&bb11eb2242172448ed8dd441742e224beb87b7ed &44&	dd111d2e8e72e8481b8db48418ee2e8bed72d11b\\
		
		&1b1877d1e2b881b7e47b8ddd4e2e2448d2eb1187&&	4d7777b412b887d2e41d7bbd4e2e2448d2eb1187\\
		
		\hline 
		
		&1248284bee77db1eed4e17287dd8b482be7b71b1&&	1248284bee77db1eed4e17287dd8b482be7b71b1\\
		
		45&dd118d22427124488b8db447142e224bede1d18b &46&	dd118d24247142488b8db427124e242bede1d18b\\
		
		&7d7e77b782b887d1e41debbd4e2e2448d2eb1187&&	7d7e77b782b887d1e41debbd4e2e2448d2eb1187\\
		
		\hline 
		
		&1248284bee77db1eed4e17287db8d482be7b71d1&&	1248284bee77db1eed4127284dd7b482b17b41b2\\
		
		47&bd118b22427124488b8dd447142e224bede1d18d &48&	dd117d2e8e72e8487b8db48418ee2e8bed12d17b\\
		
		&7b7e77b782b887d1e41bedbd4e2e2448d2eb1187&&	4d7177b472b887d2e41d1bbd4e2e2448d2eb1187\\
		
		\hline 
		
		&1248284bee77db1eed4817287dbed482b87b71d1&&	1248284bee77db1eed4877281dbed482b87b11d7\\
		
		49&bd11eb2242712448eb8dd447142e224bed81d1ed &50&	bd11eb2424774248eb8dd421124e242bed87d1ed\\
		
		&7b7877b7e2b887d1e41b8dbd4e2e2448d2eb1187&&	1b7877b1e2b887d7e41b8dbd4e2e2448d2eb1187\\
		
		\hline 
		
		&1248284bee77db7eed4721284db1d482b71b47d2&&	1248284bee77db7eed4741282db1d482b71b27d4\\
		
		51&bd111b2e8e12e8481b8dd48478ee2e8bed72d71d &52&	bd111b28e8148e481b8dd4e27e8e28ebed74d71d\\
		
		&4b1777b412b881d2e47b7dbd4e2e2448d2eb1187&&	2b1777b212b881d4e47b7dbd4e2e2448d2eb1187\\
		
		\hline 
		
		&1248284bee77db7eed4e11287dd8b482be1b77b1&&	1248284bee77db7eed4121284db7d482b11b47d2\\
		
		53&dd118d24241142488b8db427724e242bede1d78b &54&	bd117b2e8e12e8487b8dd48478ee2e8bed12d77d\\
		
		&7d1e77b782b881d1e47debbd4e2e2448d2eb1187&&	4b1177b472b881d2e47b1dbd4e2e2448d2eb1187\\
		
		\hline 
		
		&124828eed4bb7bde4872b1248d771e8d21b8d7bd&&	12482d48be77e1bebe4d2d2d471b74828dd14b72\\
		
		55&27111d2e4e21e448d17e784b84e2ee44db7721b8 &56&	1711b128e8d28e4db18e74e4be8b28e8e7d27bb7\\
		
		&71dbbbd18d2e4712e8bbdd787242e887442d11eb&&	41dd7714b28d8d72b4b1d71e4e2b244de2e81187\\	
		\hline	
	\end{longtable}}
}

Then, we apply the Grouping algorithm to classify the other 80 $GL_4(\mathbb{Z}_2)$-equivalence classes of the small covers over the right-angled 120-cell. They are all non-orientable by Theorem \ref{theorem: NakayamaN}.

\noindent \textbf{Theorem \ref{theorem:nonorientable eight-coloring}.} \emph{There are exactly two  8-coloring  non-orientable  small covers over the right-angled hyperbolic 120-cell up to homeomorphism.}

The resulting characteristic vectors are reported in Table \ref{table:nonori}, and we adopt the same denotation manner.

{\centering
		{\footnotesize	
			\renewcommand\baselinestretch{1.5}\selectfont
			\begin{longtable}{clcc}
				
				\captionsetup{font=normalsize}
				\caption{8-coloring non-orientable small covers over the 120-cell $\mathcal{E}$.}\label{table:nonori}\\
				\hline 
		
		&1248285e4ee9351324e1e21181e911481e9e5e82&&	124827b8fff4217d2f8f1d11bf11b811b4f7fbfd\\
		
		\ 57 \ &384245e395511e9e33195e31559e9e191335e9e3~\hspace*{0.3cm} &\ 58 \ &	d47278fd8417f1fb2d418fd7148bf1f14d2f7f4d\\
		
		&28e4248e1e42315811219e1184184224e5839ee1&&	dfb72741f8f2d71b11db11fb1187422f8b824ff1\\
		
				\hline

	\end{longtable}}
}

We also list the unique, up to homeomorphism, $5$-coloring characteristic vector in Table \ref{table:5col}; it was found by Garrison-Scott \cite{Scott:02}.

\begin{table}[H]
	\centering
	\footnotesize
	\renewcommand\baselinestretch{1.5}\selectfont
	\begin{tabular}{cc}
		\hline
		&12482848fff421822f8f12114f11481144f8f4f2\\
	59 \ 	&248288f28418f1f422418f281484f1f1422f8f42\\
		&2f482841f8f22814112411f41188422f84824ff1\\
		\hline
	\end{tabular}
	\captionsetup{font=normalsize}
	\caption{The unique 5-coloring small cover over the 120-cell $\mathcal{E}$.}\label{table:5col}
\end{table}

\subsection{9-coloring small covers over 120-cell $\mathcal{E}$}\label{section:9-coloring}

When considering the 9-coloring case, the computation far exceeds the capabilities of today's laptops. We summarize the details in Table \ref{table:9colblock}.

	\begin{table}[H]
		\caption{Block-pasting algorithm computation of 9-coloring small covers over the 120-cell.}\label{table:9colblock}
		\scalebox{0.7}{
			\linespread{1}\selectfont
			\begin{tabular}{|p{1em}|c|c|c|c|c|p{1em}|c|p{1em}|c|p{1em}|c|p{1em}|}  
				\hline
				& 
				\multicolumn{5}{c|}{\multirowcell{3}{Added colors}}                                          & & 
				\multirowcell{3}{$\vert B_1 \vert$} 
				& & 
				\multirowcell{3}{$\vert \overline{B_1} \vert$} 
				& & 
				\multirow{3}{8em}{\centering Number of resulting vectors} & \\
				& \multicolumn{5}{c|}{}                                                                                          & &                                                        & &                                     & &                                               & \\ 
				
				& \multicolumn{1}{p{2em}}{} & \multicolumn{1}{p{2em}}{} & \multicolumn{1}{p{2em}}{} & \multicolumn{1}{p{2em}}{} & \multicolumn{1}{p{2em}|}{} & & \multicolumn{1}{p{4em}|}{} & &\multicolumn{1}{p{4em}|}{} & &\multicolumn{1}{p{4em}|}{} & \\
				\hhline{|~|-|-|-|-|-|~|-|~|-|~|-|~|}
				
				& \cellcolor{gray!50}3 & \cellcolor{gray!50}5 & \cellcolor{gray!50}6 & \cellcolor{gray!50}7 & \cellcolor{gray!50}9 & & \cellcolor{gray!50}0 & &\multirow{2}{8em}{} & & \multirow{2}{8em}{} & \\ 
				\hhline{|~|-|-|-|-|-|~|-|~|~|~|~|~|}
				& \cellcolor{gray!50}3 & \cellcolor{gray!50}5 & \cellcolor{gray!50}6 & \cellcolor{gray!50}7 &\cellcolor{gray!50}11 & & \cellcolor{gray!50}0 & & & & & \\ \hhline{|~|-|-|-|-|-|~|-|~|-|~|-|~|}
				& 3 & 5 & 6 & 7 &15 & & 160 & &2709360 & &0 & \\
				\hhline{|~|-|-|-|-|-|~|-|~|-|~|-|~|}
				& \cellcolor{gray!50}3 & \cellcolor{gray!50}5 & \cellcolor{gray!50}6 & \cellcolor{gray!50}9 & \cellcolor{gray!50}10 & & \cellcolor{gray!50}0 & &\multirow{2}{8em}{} & & \multirow{2}{8em}{} & \\ 
				\hhline{|~|-|-|-|-|-|~|-|~|~|~|~|~|}
				& \cellcolor{gray!50}3 & \cellcolor{gray!50}5 & \cellcolor{gray!50}6 &\cellcolor{gray!50} 9 &\cellcolor{gray!50}11 & & \cellcolor{gray!50}0 & & & & & \\ \hhline{|~|-|-|-|-|-|~|-|~|-|~|-|~|}
				& 3 & 5 & 6 & 9 &14 & & 640 & &3635040 & &0 & \\
				\hhline{|~|-|-|-|-|-|~|-|~|-|~|-|~|}
				& 3 & 5 & 6 & 9 &15 & & 640 & &3635040 & &0 & \\
				\hhline{|~|-|-|-|-|-|~|-|~|-|~|-|~|}
				& 3 & 5 & 6 & 11 &13 & & 224 & &910080 & &0 & \\
				\hhline{|~|-|-|-|-|-|~|-|~|-|~|-|~|}
				& 3 & 5 & 6 & 11 &15 & & 628 & &3635040 & &0 & \\
				\hhline{|~|-|-|-|-|-|~|-|~|-|~|-|~|}
				& \cellcolor{gray!50}3 & \cellcolor{gray!50}5 & \cellcolor{gray!50}7 & \cellcolor{gray!50}9 & \cellcolor{gray!50}11 & & \cellcolor{gray!50}0 & & & & & \\
				\hhline{|~|-|-|-|-|-|~|-|~|-|~|-|~|}
				& 3 & 5 & 7 & 9 & 14 & & 640 & &3635040 & &0 & \\
				\hhline{|~|-|-|-|-|-|~|-|~|-|~|-|~|}
				& 3 & 5 & 7 & 9 & 15 & & 640 & &3635040 & &0 & \\
				\hhline{|~|-|-|-|-|-|~|-|~|-|~|-|~|}
				& \cellcolor{gray!50}3 & \cellcolor{gray!50}5 & \cellcolor{gray!50}7 & \cellcolor{gray!50}10 & \cellcolor{gray!50}11 & & \cellcolor{gray!50}0 & & & & & \\
				\hhline{|~|-|-|-|-|-|~|-|~|-|~|-|~|}
				
				& 3 & 5 & 7 & 10 & 12 & & 224 & &3635040 & &0 & \\
				\hhline{|~|-|-|-|-|-|~|-|~|-|~|-|~|}
				& 3 & 5 & 7 & 10 & 13 & & 224 & &910080  & & 0& \\
				\hhline{|~|-|-|-|-|-|~|-|~|-|~|-|~|}
				& 3 & 5 & 7 & 10 & 14 & & 628 & &3635040 & & 0& \\
				\hhline{|~|-|-|-|-|-|~|-|~|-|~|-|~|}
				& 3 & 5 & 7 & 10 & 15 & & 628 & &910080  & &0 & \\
				\hhline{|~|-|-|-|-|-|~|-|~|-|~|-|~|}
				& 3 & 5 & 7 & 11 & 13 & & 224 & &3635040 & & 0& \\
				\hhline{|~|-|-|-|-|-|~|-|~|-|~|-|~|}
				& 3 & 5 & 7 & 11 & 14 & & 628 & &3635040 & & 0& \\
				\hhline{|~|-|-|-|-|-|~|-|~|-|~|-|~|}
				& 3 & 5 & 7 & 11 & 15 & & 628 & &910080  & & 0& \\
				\hhline{|~|-|-|-|-|-|~|-|~|-|~|-|~|}
				& 3 & 5 & 7 & 14 & 15 & & 1280 & &910080  & &0 & \\
				\hhline{|~|-|-|-|-|-|~|-|~|-|~|-|~|}
				& \cellcolor{yellow!40}3 & \cellcolor{yellow!40}5 & \cellcolor{yellow!40}9 & \cellcolor{yellow!40}14 & \cellcolor{yellow!40}15 & & \cellcolor{yellow!40}5120 & &\cellcolor{yellow!40}11451840  & & \cellcolor{yellow!40}?& \\
				\hhline{|~|-|-|-|-|-|~|-|~|-|~|-|~|}
				
				& 3 & 5 & 10 & 12 & 15 & & 826 & &1764720  & &0 & \\
				\hhline{|~|-|-|-|-|-|~|-|~|-|~|-|~|}
				& 3 & 5 & 10 & 13 & 14 & & 826 & &1764720   & & 0& \\
				\hhline{|~|-|-|-|-|-|~|-|~|-|~|-|~|}
				& 3 & 5 & 10 & 13 & 15 & & 826 & &3635040 & & 0& \\
				\hhline{|~|-|-|-|-|-|~|-|~|-|~|-|~|}
				& 3 & 5 & 11 & 13 & 14 & & 826 & &3635040 & &0 & \\
				\hhline{|~|-|-|-|-|-|~|-|~|-|~|-|~|}
				& 3 & 5 & 11 & 13 & 15 & & 826 & &1764720 & & 0& \\
				\hhline{|~|-|-|-|-|-|~|-|~|-|~|-|~|}
				& 3 & 5 & 11 & 14 & 15 & & 1280 & &3635040  & &0 & \\
				\hhline{|~|-|-|-|-|-|~|-|~|-|~|-|~|}
				& 3 & 7 & 11 & 12 & 13 & & 560 & &3635040  & &0 & \\
				\hhline{|~|-|-|-|-|-|~|-|~|-|~|-|~|}
				& 3 & 7 & 11 & 12 & 15 & & 4490 & &2709360  & &0 & \\
				\hhline{|~|-|-|-|-|-|~|-|~|-|~|-|~|}
				
				& \cellcolor{yellow!40}3 & \cellcolor{yellow!40}7 & \cellcolor{yellow!40}11 & \cellcolor{yellow!40}13 & \cellcolor{yellow!40}14 & & \cellcolor{yellow!40}4490 & &\cellcolor{yellow!40}11451840  & &\cellcolor{yellow!40}? & \\
				\hhline{|~|-|-|-|-|-|~|-|~|-|~|-|~|}
				
				& 3 & 7 & 11 & 13 & 15 & & 4490 & &2709360  & & 0& \\
				\hhline{|~|-|-|-|-|-|~|-|~|-|~|-|~|}
				& 3 & 7 & 12 & 13 & 15 & & 680 & &3635040   & &0 & \\
				\hhline{|~|-|-|-|-|-|~|-|~|-|~|-|~|}
				& 3 & 7 & 13 & 14 & 15 & & 1256 & &910080  & & 0& \\
				\hhline{|~|-|-|-|-|-|~|-|~|-|~|-|~|}
				
				& \cellcolor{yellow!40}7 & \cellcolor{yellow!40}11 & \cellcolor{yellow!40}13 & \cellcolor{yellow!40}14 &\cellcolor{yellow!40} 15 & & \cellcolor{yellow!40}7800 & &\cellcolor{yellow!40}11451840   & &\cellcolor{yellow!40}? & \\			
				\hline
			\end{tabular}
		}
	\end{table}
	
So far, we can only report that, the 9-coloring characteristic functions exist within the following three cases of added colors:
$$(3,~5,~9,~14,~15), (3,7,11,13,14),  (7,11,13,14,15).$$
\noindent The three 9-coloring results we have found so far are all homeomorphic to the coloring in Table \ref{table:9col}, which is exactly what we obtain in Section \ref{section:9col_cr} by the color-recursion algorithm. It has the same colors with the 5-coloring resulting vectors in Table \ref{table:5col} in the first $33$ facets.

\begin{table}[H]
	\centering
	\footnotesize
	\renewcommand\baselinestretch{1.5}\selectfont
	\begin{tabular}{cc}
		\hline
		&12482848fff421822f8f12114f11481144f8f4f2\\
		60 \	&248289e39519e1e533519e391595e1e1533f8f42\\
		&2f482841f8f22814112411f41188422f84824ff1\\
		
		\hline
	\end{tabular}
	\captionsetup{font=normalsize}
	\caption{A known 9-coloring small cover over the right-angled 120-cell $\mathcal{E}$.}\label{table:9col}
\end{table}

Therefore, we find a total of 59 new characteristic functions of small covers over the 120-cell. The notation $M(\mathcal{E},\lambda_i)$, abbreviated as $M_i$, for $1\le i \le 60$, is used to denote the corresponding 60 known small covers over the 120-cell $\mathcal{E}$. This notation appears in Theorem \ref{theorem:manifold}.

\section{Intersection forms of closed hyperbolic $4$-manifolds}
 Davis and Januszkiewicz formulated how to obtain the $\mathbb{Z}_2$-coefficient cohomology groups of a small cover by the orbit polytope and its characteristic function \cite[Theorem 4.14]{dj:1991}. In 2013, Li Cai gave  a method to calculate the $\mathbb{Z}$-coefficient cohomology groups of the real moment-angle manifold over a simplicial complex $K$ \cite{2licai:2013}. Based on the work of Cai and Suciu-Trevisan's result on rational homology groups of real toric manifolds  \cite {Suciu:13, SuciuT:12}, Choi and Park gave a formula of the cohomology groups of real toric manifolds \cite{ChoiP:2013}, which can also be viewed as a combinatorial version of the Hochster Theorem \cite{Hochster:1977}.

Suppose that  $P$ is a simple polytope and $K$ is the simplicial complex dual to the boundary of $P$. Recall that we can also define the \emph{characteristic function} $\lambda$ on $K$ by substituting  the facet set $\mathcal{F}(P)$ with the vertex set $\mathcal{V}$ of the simplicial complex $K$. See Section 2.1 for more details. We denote the linear space $\mathbb{Z}_2^{\vert\mathcal{V}\vert}$ by $\mathbb{Z}_2^{\mathcal{V}}$. In addition, we can identify $\mathbb{Z}_2^{\mathcal{V}}$ with the power set $2^{\mathcal{V}}$ in the canonical way, where $\emptyset$ corresponds to the identity element and multiplication to the symmetric difference. Namely, we have a map $\varphi:\mathbb{Z}_2^{\mathcal{V}}\longrightarrow 2^{\mathcal{V}}$. 
Accordingly, every full-subcomplex $K_\omega$ of $K$, $\omega \subseteq {\mathcal{V}}$, is identified with an element of $2^{\mathcal{V}}$, i.e. $K_{\omega}$ is the full sub-complex of $K=(\partial L)^*$ by restricting to $\omega \subseteq \mathcal{V}$. The notations $\varphi$, $K$ and $K_{\omega}$ will be used throughout this section with this meaning.

Let $\lambda$ be a $\mathbb{Z}_2^r$-coloring characteristic function. Denote by $row\Lambda$ the \emph{row space} of the characteristic matrix $\Lambda$. The following Choi-Park theorem shows that the cohomology group of a real toric manifold $M(P,\lambda)$ is the direct sum of the cohomology groups of some full sub-complexes of the dual polytope $K=(\partial L)^*$. The sub-complexes are determined by the characteristic function.

\begin{theorem}  \label{theorem: ChoiP} \emph{(Choi-Park \cite{ChoiP:2013})} Assume $G$ is the coefficient ring $\mathbb{Q}$ or $\mathbb{Z}_q$ for a positive odd integer $q$. There is an additive isomorphism 
	$$H^p(M(P,\lambda);G)\cong\underset{\varphi^{-1}(\omega) \in row \Lambda}\oplus\widetilde{H}^{p-1}(K_\omega;G).$$
\end{theorem}

0.
\begin{corollary} \label{corollary: ChoiP}
	For a simple polytope $P$,  $$\beta^1(M(P,\lambda);\mathbb{Q})=\underset{\varphi^{-1}(\omega) \in row \Lambda}\sum{\widetilde\beta^{0}}(K_\omega;\mathbb{Q}), $$
	where $K_{\omega}$ is the full sub-complex of $K=(\partial P)^*$ by restricting to $\omega \subseteq \mathcal{V}$.
\end{corollary}

By means of Corollary \ref{corollary: ChoiP}, we may calculate the first Betti number of a real toric manifold using the combinatorial information of the colored polytope and the row space of its characteristic matrix.

Besides, Cai and Choi discussed the cohomology groups of real toric manifolds in \cite{caichoi:2017}. Let $K=(\partial P)^*$ and denote by $b,c$ and $d$ the sum of the zeroth, first and second reduced Betti numbers of $K_{\omega}$ for $\omega\in\text{row}\Gamma$. They showed the following table in their paper. 

\begin{table}[H]
	\small
	\begin{tabular}{|c|P{3cm}|P{4cm}|}
		\hline
		$H^i(M(P,\lambda);\mathbb{Z})$& Orientable& Non-orientable\\
		\hline	
	    $i=0$&$\mathbb{Z}$ &$\mathbb{Z}$\\
	    $i=1$&$\mathbb{Z}^b$ &$\mathbb{Z}^b$\\
	    $i=2$&$\mathbb{Z}^c\oplus\mathbb{Z}_2^{m-4-b}$ &$\mathbb{Z}^c\oplus\mathbb{Z}_2^{m-4-b}$\\
	    $i=3$&$\mathbb{Z}^b\oplus\mathbb{Z}_2^{m-4-b}$ &$\mathbb{Z}^d\oplus\mathbb{Z}_2^{m-5-d}$\\
	    $i=4$&$\mathbb{Z}$ &$\mathbb{Z}_2$\\
		\hline
	\end{tabular}
	\caption{Cohomologies of 4-dimensional small cover}
	\label{table:cohomoligies_4dsc}
\end{table}

From Table \ref{table:cohomoligies_4dsc}, Theorem 3.1 in \cite{dj:1991} and the fact that the Euler characteristic is independent of the choice of the coefficient ring, one sees that the integral cohomology group $H^*(M(P,\lambda);\mathbb{Z})$ is completely determined by the reduced cohomology group of $K_{\omega}$ for $\omega\in row \Lambda$ and the $h$-vector of $K$. As corollaries, small covers satisfy some topological restrictions. For instance, in dimension $n=4$, the only (if any) torsion in homology is 2-torsion.

Over the right-angled 120-cell, there is a canonical way to construct a natural  $\mathbb{Z}_2^5$-extension $\delta$ from a non-orientable $\mathbb{Z}_2^4$-coloring $\lambda$ such that the real toric manifold $M(\mathcal{E},\delta)$ is the orientable double cover of the non-orientable small cover $M(\mathcal{E},\lambda)$. See \cite{mz} for more details. Thus, we have four natural $\mathbb{Z}_2^5$-extensions $\delta_{57},~\delta_{58},~\delta_{59}$, and $\delta_{60}$ of the non-orientable $\mathbb{Z}_2^4$-colorings $\lambda_{57},~\lambda_{58},~\lambda_{59}$, and $\lambda_{60}$ respectively. We use $M(\mathcal{E},\delta_i)$, abbreviated as $\widetilde{M_i}$, to denote the corresponding real toric manifold, where $57\le i\le 60$.  These notations appear in Theorem \ref{theorem:manifold}.

\noindent \textbf{Proof of Theorem \ref{theorem:manifold} (homology group part):} By Theorem \ref{theorem: ChoiP}, the first Betti numbers of the known 60 small covers $M(\mathcal{E}, \lambda_i)$ can be computed by the row spaces of their characteristic matrices. By calculation, each subcomplex $K_{\omega_i}$ of every small cover, $\varphi^{-1}(\omega_i)$ is a row in $row\Lambda$ and $1\leq i\leq 15$ , has at most one connected component. Thus, the first Betti numbers are all zero.

Next, we use Theorem \ref{theorem: ChoiP} to calculate the first Betti numbers of the four orientable real toric manifolds $\widetilde{M_{57}}$, $\widetilde{M_{58}}$, $\widetilde{M_{59}}$ and $\widetilde{M_{60}}$. The results are summarized in Table \ref{table:betti}.

 \begin{table}[H]
 	\centering
 	\tiny
 	\renewcommand\baselinestretch{1.5}\selectfont
 	 \begin{tabular}{|c|P{0.1cm}|P{0.1cm}|P{0.1cm}|P{0.1cm}|P{0.1cm}|P{0.1cm}|P{0.1cm}|P{0.1cm}|P{0.1cm}|P{0.1cm}|P{0.1cm}|P{0.1cm}|P{0.1cm}|P{0.1cm}|P{0.1cm}|P{0.1cm}|P{0.1cm}|P{0.1cm}|P{0.1cm}|P{0.1cm}|P{0.1cm}|P{0.1cm}|P{0.1cm}|P{0.1cm}|P{0.1cm}|P{0.1cm}|P{0.1cm}|P{0.1cm}|P{0.1cm}|P{0.1cm}|P{0.1cm}|P{0.3cm}|}
 		\hline
 		\diagbox{$i$}{$\widetilde\beta^0(j)$}{$j$}  &1&2&3&4&5&6&7&8&9&10&11&12&13&14&15&16&17&18&19&20&21&22&23&24&25&26&27&28&29&30&31&$\beta^1$\\
 		
    	\hline
    	57&0&0&0&0&7&0&0&0&23&0&0&0&0&0&0&0&0&0&0&0&0&0&7&7&7&0&0&0&0&0&0&51\\
    	
    	\hline
    	58&0&0&0&0&23&0&0&0&0&0&0&7&0&7&7&0&0&0&0&0&0&0&0&0&0&0&0&0&0&7&0&51\\
    	
    	\hline
    	59&0&0&0&0&23&0&0&0&23&0&0&23&0&23&23&0&0&0&0&0&0&0&0&0&0&0&0&0&0&0&0&115\\
    	
 		\hline
 		60&0&0&0&0&7&0&0&0&23&0&0&7&0&7&7&0&0&0&0&0&0&0&0&0&0&0&0&0&0&0&0&51\\
 		\hline
 	\end{tabular}
 	\caption{$\beta^1$ of $M(\mathcal{E},\delta_i)$ }
 	\label{table:betti}
 \end{table}
\noindent Here, $\widetilde\beta^0(j)$ equals the 0-th Betti number of the $j$-th full sub-complex $K_{\omega_j}$, $\beta^1$ is the first Betti number of the  manifold $M(\mathcal{E},\delta_i)$, and $i$ represents the serial numbers of the characteristic vectors as claimed in Tables \ref{table:ori}--\ref{table:9col}. For example, considering $K_{\omega_5}$ of $\lambda_{60}$, the first full sub-complex with a non-trivial reduced zeroth homology group, the corresponding row vector is:

\begin{center}
	{\footnotesize
		\noindent$(
		0,0,0,0,0,0,0,0,1,1,1,0,0,0,0,0,0,1,0,1,0,0,0,0,0,1,0,0,0,0,0,0,0,0,1,0,1,0,1,0,0,0,0,0,0,1,0,1,1,1,0,$
		
		\noindent$
		1,0,0,0,1,1,1,1,0,1,0,1,1,0,1,1,1,0,0,0,0,1,1,1,1,0,1,0,0,0,1,0,0,0,0,0,0,1,0,1,0,0,0,0,0,0,0,0,0,0,0,$
	
	\noindent$1,0,0,0,0,0,0,0,0,1,0,0,0,0,0,1,1,0).$}
\end{center}

\noindent The orders of facets being selected by the non-trivial entries are:
\begin{center}
	{\footnotesize
		\noindent$(9,10,11,18,20,26,35,37,39,46,48,49,50,52,56,57,58,59,61,63,64,66,67,68,73,74,75,76,78,82,89,91,103,112,118,119).$}
\end{center}

According to the adjacency matrix $M_{\mathcal{E}}$, there are a total of eight connected components. This means the reduced zeroth Betti number of $K_{\omega_5}$ is seven as shown in Table \ref{table:betti}. The eight groups of facets with respect to  the eight connected components are shown in Table \ref{table:betti_facet}.

\begin{table}[H]
	\small
	\begin{tabular}{|c|c|c|c|P{3cm}|P{4cm}|c|c|c|}
		\hline
		&1st&2nd&3rd&4th&5th&6th&7th&8th\\
		\hline	
	   labels of&&&&18, 20, 26, 48, 49,&35, 37, 39, 46, 50, &&&\\
      included &9&10&11&57, 59, 61, 64, 66, &52, 56, 58, 63, 67 &112&118&119\\
     facets&&&&68, 74, 76, 78, 82 &73, 75, 89, 91, 103 &&&\\
     \hline
\end{tabular}
\caption{Labels of facets with respect to the eight connected components of $K_{\omega_5}$}
\label{table:betti_facet}
\end{table}

With Table \ref{table:cohomoligies_4dsc}, Theorem 3.1 of \cite{dj:1991} and some simple calculations, we complete the homology part of the proof of Theorem \ref{theorem:manifold}.   \qed

\vspace{0.5cm}

Suppose M is a compact, oriented, topological $4$-manifold. We have the concept of intersection form $Q(M)$ as follows \cite{kirby_calculus}.

\begin{definition}
	The symmetric bilinear form $Q(M): H^2(M, \partial M; \mathbb{Z}) \times H^2(M, \partial M; \mathbb{Z}) \rightarrow \mathbb{Z}$
	defined by $Q(M)(a, b)=\langle a\cup b, [M]\rangle = a \cdot b\in \mathbb{Z}$ is called the intersection form
of $M$, where $[M]\in H^4(M, \partial M;\mathbb{Z})$ is the fundamental class. 
\end{definition}
Note that  $Q(M)(a, b)=0$ if $a$ or $b$ is a torsion element. Hence $Q(M)$
descends to a pairing on homology mod torsion. By choosing a basis of
the free part of $H^2(M,\partial M; \mathbb{Z})$, we can represent $Q(M)$ by a matrix. 

By the Hirzebruch
 signature formula, every closed oriented hyperbolic four-manifold has zero signature. Therefore, the intersection form of a closed oriented hyperbolic $4$-manifold $M$ is either
 \begin{center}
 	$\oplus_{k}\begin{pmatrix}
 	0 & 1\\
 	1&0\\
 	\end{pmatrix}$  or $\oplus_{k}(1)\oplus_{k}(-1)$.
 \end{center}
The Davis manifold $M^*$ constructed by M. Davis in \cite{davis:1985} is a hyperbolic 4-manifold obtained by the $\frac{2\pi}{5}$-angled 120-cell. From the characteristic number $\chi(M^*)=26$, the volume of $M^*$ can be obtained immediately by the Gauss-Bonnet formula as well.
Moreover, Ratcliffe and Tschantz calculated the homology of the Davis manifold and it is the first closed hyperbolic $4$-manifold example for which the signature and intersection form have been determined \cite{rt:2001}.

For real toric manifolds, Choi and Park have studied the  multiplication structure of the cohomology ring in \cite{choi_park:2017}. Suppose $M(P,\lambda)$ is a real toric manifold and $\Lambda$ is the characteristic matrix of $\lambda$. We have 

\begin{theorem}  \label{theorem: ChoiP2} \emph{(Choi-Park \cite{choi_park:2017})} Suppose $\mathcal{R}$ is a commutative ring with unity in which $2$ is a unit.There is a ring isomorphism
$$H^*(M(P,\lambda);\mathcal{R})\cong \underset{\omega\in row \Lambda}\oplus\widetilde{H}^{*-1}(K_\omega;\mathcal{R}),$$
where the product structure on $\underset{\omega\in row \Lambda}\oplus\widetilde{H}^{*-1}(K_\omega;\mathcal{R})$ is given by the
canonical maps
$$\widetilde{H}^{k-1}(K_{\omega_1};\mathcal{R})\otimes \widetilde{H}^{l-1}(K_{\omega_2};\mathcal{R})\rightarrow \widetilde{H}^{k+l-1}(K_{\omega_1+\omega_2};\mathcal{R}),$$
which are induced by simplicial maps $K_{\omega_1+\omega_2} \rightarrow K_
{\omega_1} \star K_
{\omega_2},$ where $\star$ denotes
the simplicial join and $K_{\omega_1+\omega_2}$ is the full sub-complex of $K=(\partial P)^*$ by restricting to $\varphi(\varphi^{-1}(\omega_1)+\varphi^{-1}(\omega_2)) \subseteq \mathcal{V}$.
\end{theorem}

By Theorem \ref{theorem: ChoiP2}, we have the following theorem.
\begin{theorem} (The intersection part of Theorem \ref{theorem:manifold})
	Let $M$ be a closed, oriented, hyperbolic, real toric 4-manifold with second Betti number $\beta_2$. Then the intersection form of $M$ is equivalent over 
	$\mathbb{Z}$ to $\oplus_{\beta_2/2}\begin{pmatrix}
	0 & 1\\
	1&0\\
	\end{pmatrix}$. 
\end{theorem}

\begin{proof}
	By the Hirzebruch signature formula, every closed oriented hyperbolic 4-manifold has zero signature (see \cite{longR:2000} for instance). Thus $M$ has zero signature. In addition, one of the referees also point out one more straightforward approach specific to the real toric manifolds as follow: an orientable 4-dimensional real toric manifold, not necessarily hyperbolic, always admits an orientation-reversing involution, obtained by reflecting $M$ along a hypersurface corresponding to a facet of the orbit polytope. Thus the signature of $M$ vanishes.
	 
	Suppose a basis of $\underset{\omega_i\in row~ \Lambda}\oplus\widetilde{H}^{1}(K_{\omega_i};\mathbb{Q})$ is $e_1,e_2,\cdots, e_l$. From the product structure on  $ \underset{\omega_i,\omega_j\in row \Lambda}\oplus\widetilde {H}^{3}(K_{\omega_i+\omega_j};\mathbb{Q})$, it is obvious that $e_i\otimes e_j$ and $e_j\otimes e_i$ map to the same element in $\underset{\omega_i,\omega_j\in row \Lambda}\oplus\widetilde {H}^{3}(K_{\omega_i+\omega_j};\mathbb{Q})$. In addition, $\widetilde{H}^1(K_{\omega_i};\mathbb{Q})\otimes \widetilde{H}^1(K_{\omega_i};\mathbb{Q})\rightarrow \widetilde{H}^{3}(K_{\omega_i+\omega_i};\mathbb{Q})$ is always a zero mapping for $K_{\omega_i+\omega_i}=\emptyset$. Then, it can be inferred that the diagonals of the matrix of the intersection form of $M$ are all zero. We then complete the proof by \cite[Exercise 1.2.7]{kirby_calculus}.
\end{proof}

\begin{appendices}
\vspace*{-0.5cm}	
{\renewcommand\baselinestretch{1.28}\selectfont

 \begin{table}[H]
 	\begin{floatrow}
 		\capbtabbox{
 			\footnotesize
 			\begin{tabular}{|c|c|c|}
 				\hline
 				layer &  $F_i$&  $q_i$   \\
 				\hline
 				1 & $F_1$ &(1, 0, 0, 0)   \\
 				\hline
 				&$F_2$&$\frac{1}{2}(\phi,~~0,~~\phi^{-1},~~1)$\\
 				\cline{2-3}
 				& $F_3$ &$\frac{1}{2}(\phi,~~0,~~-\phi^{-1},~~1)$\\
 				\cline{2-3}
 				& $F_4$ &$\frac{1}{2}(\phi,~~1,~~0,~~\phi^{-1})$\\
 				\cline{2-3}
 				& $F_5$ &$\frac{1}{2}(\phi,~~0,~~\phi^{-1},~~-1)$\\
 				\cline{2-3}
 				& $F_6 $ &$\frac{1}{2}(\phi,~~0,~~-\phi^{-1},~~-1)$\\
 				\cline{2-3}
 				2 & $F_7$ &$\frac{1}{2}(\phi,~~1,~~0,~~-\phi^{-1})$\\
 				\cline{2-3}
 				& $F_8$ &$\frac{1}{2}(\phi,~~-1,~~0,~~\phi^{-1})$\\
 				\cline{2-3}
 				& $F_9$ &$\frac{1}{2}(\phi,~~-1,~~0,~~-\phi^{-1})$\\
 				\cline{2-3}
 				& $F_{10}$ &$\frac{1}{2}(\phi,~~\phi^{-1},~~1,~~0)$\\
 				\cline{2-3}
 				& $F_{11}$ &$\frac{1}{2}(\phi,~~\phi^{-1},~~-1,~~0)$\\
 				\cline{2-3}
 				& $F_{12}$ &$\frac{1}{2}(\phi,~~-\phi^{-1},~~1,~~0)$\\
 				\cline{2-3}
 				& $F_{13}$ &$\frac{1}{2}(\phi,~~-\phi^{-1},~~-1,~~0)$\\
 				\hline
 				 & $F_{14}$ &$\frac{1}{2}(1,~~1,~~1,~~1)$\\
 				\cline{2-3}
 				& $F_{15}$ &$\frac{1}{2}(1,~~1,~~1,~~-1)$\\
 				\cline{2-3}
 				& $F_{16}$ &$\frac{1}{2}(1,~~1,~~-1,~~1)$\\
 				\cline{2-3}
 				& $F_{17}$ & $\frac{1}{2}(1,~~1,~~-1,~~-1)$\\
 				\cline{2-3}
 				& $F_{18}$ &$\frac{1}{2}(1,~~-1,~~1,~~1)$\\
 				\cline{2-3}
 				& $F_{19}$ &$\frac{1}{2}(1,~~-1,~~1,~~-1)$\\
 				\cline{2-3}
 				& $F_{20}$ &$\frac{1}{2}(1,~~-1,~~-1,~~1)$\\
 				\cline{2-3}
 				& $F_{21}$ &$\frac{1}{2}(1,~~-1,~~-1,~~-1)$\\	
 				\cline{2-3}
 				& $F_{22}$ &$\frac{1}{2}(1,~~\phi,~~~~\phi^{-1},~~0)$\\
 				\cline{2-3}
 				& $F_{23}$ &$\frac{1}{2}(1,~~\phi,~~~~-\phi^{-1},~~0)$\\
 				\cline{2-3}
 				3& $F_{24}$ &$\frac{1}{2}(1,~~-\phi,~~~~\phi^{-1},~~0)$\\
 				\cline{2-3}
 				& $F_{25}$ &$\frac{1}{2}(1,~~-\phi,~~~~-\phi^{-1},~~0)$\\
 				\cline{2-3}
 				
 				& $F_{26}$ &$\frac{1}{2}(1,~~\phi^{-1},~~0,~~\phi)$\\
 				\cline{2-3}
 				& $F_{27}$ &$\frac{1}{2}(1,~~\phi^{-1},~~0,~~-\phi)$\\
 				\cline{2-3}
 				& $F_{28}$ &$\frac{1}{2}(1,~~-\phi^{-1},~~0,~~\phi)$\\
 				\cline{2-3}
 				& $F_{29}$ &$\frac{1}{2}(1,~~-\phi^{-1},~~0,~~-\phi)$\\
 				\cline{2-3}
 				& $F_{30}$ &$\frac{1}{2}(1,~~0,~~\phi,~~\phi^{-1})$\\
 				\cline{2-3}
 				& $F_{31}$ &$\frac{1}{2}(1,~~0,~~\phi,~~-\phi^{-1})$\\
 				\cline{2-3}
 				& $F_{32}$ &$\frac{1}{2}(1,~~0,~~-\phi,~~\phi^{-1})$\\
 				\cline{2-3}
 				& $F_{33}$ &$\frac{1}{2}(1,~~0,~~-\phi,~~-\phi^{-1})$\\

 				\hline
 					& $F_{34}$ &$\frac{1}{2}(\phi^{-1},~~0,~~1,~~\phi)$\\
 					\cline{2-3}
 					& $F_{35}$ &$\frac{1}{2}(\phi^{-1},~~0,~~1,~~-\phi)$\\	
 					\cline{2-3}
 					 & $F_{36}$ &$\frac{1}{2}(\phi^{-1},~~0,~~-1,~~\phi)$\\	
 					\cline{2-3}
 					& $F_{37}$ &$\frac{1}{2}(\phi^{-1},~~0,~~-1,~~-\phi)$\\	
 					\cline{2-3}
 					& $F_{38}$ &$\frac{1}{2}(\phi^{-1},~~\phi,~~0,~~1)$\\
 						\cline{2-3}	
 						& $F_{39}$ &$\frac{1}{2}(\phi^{-1},~~\phi,~~0,~~-1)$\\	
 					
 \cline{2-3}	
 4& $F_{40}$ &$\frac{1}{2}(\phi^{-1},~~-\phi,~~0,~~1)$\\	
 \cline{2-3}
 & $F_{41}$ &$\frac{1}{2}(\phi^{-1},~~-\phi,~~0,~~-1)$\\	
 \cline{2-3}
 & $F_{42}$ &$\frac{1}{2}(\phi^{-1},~~1,~~\phi,~~0)$\\
 \cline{2-3}	
 & $F_{43}$ &$\frac{1}{2}(\phi^{-1},~~1,~~-\phi,~~0)$\\	
 \cline{2-3}
 & $F_{44}$ &$\frac{1}{2}(\phi^{-1},~~-1,~~\phi,~~0)$\\	
 \cline{2-3}
 & $F_{45}$ &$\frac{1}{2}(\phi^{-1},~~-1,~~-\phi,~~0)$\\
 \hline
 			\end{tabular}
 		}{
 		\caption{}
 		\label{table:coor1}
 	}

 	\capbtabbox{
 		\footnotesize
 		\begin{tabular}{|c|c|c|}
 			\hline
 			layer &  $F_i$&  $q_i$   \\	
 			\hline
 			& $F_{46}$ &$(0,~~1,~~0,~~0)$\\
 			\cline{2-3}
 			& $F_{47}$ &$(0,~~0,~~1,~~0)$ \\	
 			\cline{2-3}
 			& $F_{48}$ &$(0,~~0,~~0,~~1)$\\	
 			\cline{2-3}
 			& $F_{49}$ &$\frac{1}{2}(0,~~1,~~\phi^{-1},~~\phi)$\\	
 			\cline{2-3}
 			& $F_{50}$ &$\frac{1}{2}(0,~~1,~~\phi^{-1},~~-\phi)$\\	
 			\cline{2-3}
 			& $F_{51}$ &$\frac{1}{2}(0,~~1,~~-\phi^{-1},~~\phi)$\\	
 			\cline{2-3}
 			& $F_{52}$ &$\frac{1}{2}(0,~~1,~~-\phi^{-1},~~-\phi)$\\	
 			\cline{2-3}
 			& $F_{53}$ &$\frac{1}{2}(0,~~\phi,~~1,~~\phi^{-1})$ \\	
 			\cline{2-3}
 			& $F_{54}$ &$\frac{1}{2}(0,~~\phi,~~1,~~-\phi^{-1})$ \\	
 			\cline{2-3}
 			& $F_{55}$ &$\frac{1}{2}(0,~~\phi,~~-1,~~\phi^{-1})$ \\	
 			\cline{2-3}
 			& $F_{56}$ &$\frac{1}{2}(0,~~\phi,~~-1,~~-\phi^{-1})$\\	
 			\cline{2-3}
 			& $F_{57}$ &$\frac{1}{2}(0,~~\phi^{-1},~~\phi,~~1)$\\	
 			\cline{2-3}
 			& $F_{58}$ &$\frac{1}{2}(0,~~\phi^{-1},~~\phi,~~-1)$\\	
 			\cline{2-3}
 			& $F_{59}$ &$\frac{1}{2}(0,~~\phi^{-1},~~-\phi,~~1)$ \\	
 			\cline{2-3}
 			& $F_{60}$ &$\frac{1}{2}(0,~~\phi^{-1},~~-\phi,~~-1)$\\			
 			\cline{2-3}
 			5& $F_{61}$ &$(0,~~-1,~~0,~~0)$ \\	
 			\cline{2-3}
 			& $F_{62}$ & $(0,~~0,~~-1,~~0)$ \\	
 			\cline{2-3}
 			& $F_{63}$ & $(0,~~0,~~0,~~-1)$ \\	
 			\cline{2-3}
 			& $F_{64}$ & $\frac{1}{2}(0,~~-1,~~\phi^{-1},~~\phi)$\\	
 			\cline{2-3}
 			& $F_{65}$ & $\frac{1}{2}(0,~~-1,~~\phi^{-1},~~-\phi)$ \\	
 			\cline{2-3}
 			& $F_{66}$ & $\frac{1}{2}(0,~~-1,~~-\phi^{-1},~~\phi)$ \\	
 			\cline{2-3}
 			& $F_{67}$ & $\frac{1}{2}(0,~~-1,~~-\phi^{-1},~~-\phi)$\\	
 			\cline{2-3}
 			& $F_{68}$ & $\frac{1}{2}(0,~~-\phi,~~1,~~\phi^{-1})$\\	
 			\cline{2-3}
 			& $F_{69}$ &$\frac{1}{2}(0,~~-\phi,~~1,~~-\phi^{-1})$\\	
 			\cline{2-3}
 			& $F_{70}$ &$\frac{1}{2}(0,~~-\phi,~~-1,~~\phi^{-1})$\\	
 			\cline{2-3}
 			& $F_{71}$ &$\frac{1}{2}(0,~~-\phi,~~-1,~~-\phi^{-1})$\\	
 			\cline{2-3}
 			& $F_{72}$ &$\frac{1}{2}(0,~~-\phi^{-1},~~\phi,~~1)$\\	
 			\cline{2-3}
 			& $F_{73}$ & $\frac{1}{2}(0,~~-\phi^{-1},~~\phi,~~-1)$ \\	
 			\cline{2-3}
 			& $F_{74}$ & $\frac{1}{2}(0,~~-\phi^{-1},~~-\phi,~~1)$\\	
 			\cline{2-3}
 			& $F_{75}$ &$\frac{1}{2}(0,~~-\phi^{-1},~~-\phi,~~-1)$ \\	
 			\hline	
 & $F_{76}$ &$\frac{1}{2}(-\phi^{-1},~~0,~~1,~~\phi)$\\
 \cline{2-3}
 & $F_{77}$ &$\frac{1}{2}(-\phi^{-1},~~0,~~1,~~-\phi)$\\	
 \cline{2-3}
 & $F_{78}$ &$\frac{1}{2}(-\phi^{-1},~~0,~~-1,~~\phi)$\\	
 \cline{2-3}
 & $F_{79}$ &$\frac{1}{2}(-\phi^{-1},~~0,~~-1,~~-\phi)$\\	
 \cline{2-3}
 & $F_{80}$ &$\frac{1}{2}(-\phi^{-1},~~\phi,~~0,~~1)$\\
 \cline{2-3}	
 & $F_{81}$ &$\frac{1}{2}(-\phi^{-1},~~\phi,~~0,~~-1)$\\	
 \cline{2-3}
 6& $F_{82}$ &$\frac{1}{2}(-\phi^{-1},~~-\phi,~~0,~~1)$\\	
 \cline{2-3}
 & $F_{83}$ &$\frac{1}{2}(-\phi^{-1},~~-\phi,~~0,~~-1)$\\	
 \cline{2-3}
 & $F_{84}$ &$\frac{1}{2}(-\phi^{-1},~~1,~~\phi,~~0)$\\
 \cline{2-3}	
 & $F_{85}$ &$\frac{1}{2}(-\phi^{-1},~~1,~~-\phi,~~0)$\\	
 \cline{2-3}
 & $F_{86}$ &$\frac{1}{2}(-\phi^{-1},~~-1,~~\phi,~~0)$\\	
 \cline{2-3}
 & $F_{87}$ &$\frac{1}{2}(-\phi^{-1},~~-1,~~-\phi,~~0)$\\
 \hline	
 		\end{tabular}
 	}{\caption{}
 	\label{table:coor2}
 }

 \capbtabbox{
 	\footnotesize
 	\begin{tabular}{|c|c|c|}
 		\hline
 		layer &  $F_i$&  $q_i$   \\
 		\hline	  		
 		& $F_{88}$ &$\frac{1}{2}(-1,~~1,~~1,~~1)$\\
 		\cline{2-3}
 		& $F_{89}$ &$\frac{1}{2}(-1,~~1,~~1,~~-1)$\\
 		\cline{2-3}
 		& $F_{90}$ &$\frac{1}{2}(-1,~~1,~~-1,~~1)$\\
 		\cline{2-3}
 		& $F_{91}$ &$\frac{1}{2}(-1,~~1,~~-1,~~-1)$\\
 		\cline{2-3}
 		& $F_{92}$ &$\frac{1}{2}(-1,~~-1,~~1,~~1)$\\
 		\cline{2-3}
 		& $F_{93}$ &$\frac{1}{2}(-1,~~-1,~~1,~~-1)$\\
 		\cline{2-3}
 		& $F_{94}$ &$\frac{1}{2}(-1,~~-1,~~-1,~~1)$\\
 		\cline{2-3}
 		& $F_{95}$ &$\frac{1}{2}(-1,~~-1,~~-1,~~-1)$\\	
 		\cline{2-3}
 		& $F_{96}$ &$\frac{1}{2}(-1,~~\phi,~~~~\phi^{-1},~~0)$\\
 		\cline{2-3}
 		& $F_{97}$ &$\frac{1}{2}(-1,~~\phi,~~~~-\phi^{-1},~~0)$\\
 		\cline{2-3}
 		7& $F_{98}$ &$\frac{1}{2}(-1,~~-\phi,~~~~\phi^{-1},~~0)$\\
 		\cline{2-3}
 		& $F_{99}$ &$\frac{1}{2}(-1,~~-\phi,~~~~-\phi^{-1},~~0)$\\
 		\cline{2-3}
 		& $F_{100}$ &$\frac{1}{2}(-1,~~\phi^{-1},~~0,~~\phi)$\\
 		\cline{2-3}
 		& $F_{101}$ &$\frac{1}{2}(-1,~~\phi^{-1},~~0,~~-\phi)$\\
 		\cline{2-3}
 		& $F_{102}$ &$\frac{1}{2}(-1,~~-\phi^{-1},~~0,~~\phi)$\\
 		\cline{2-3}
 		& $F_{103}$ &$\frac{1}{2}(-1,~~-\phi^{-1},~~0,~~-\phi)$\\
 		\cline{2-3}
 		& $F_{104}$ &$\frac{1}{2}(-1,~~0,~~\phi,~~\phi^{-1})$\\
 		\cline{2-3}
 		& $F_{105}$ &$\frac{1}{2}(-1,~~0,~~\phi,~~-\phi^{-1})$\\
 		\cline{2-3}
 		& $F_{106}$ &$\frac{1}{2}(-1,~~0,~~-\phi,~~\phi^{-1})$\\
 		\cline{2-3}
 		& $F_{107}$ &$\frac{1}{2}(-1,~~0,~~-\phi,~~-\phi^{-1})$\\
 		\hline
 		& $F_{108}$ &$\frac{1}{2}(-\phi,~~0,~~\phi^{-1},~~1)$\\
 		\cline{2-3}
 		& $F_{109}$ &$\frac{1}{2}(-\phi,~~0,~~\phi^{-1},~~-1)$\\
 		\cline{2-3}
 		& $F_{110}$ &$\frac{1}{2}(-\phi,~~0,~~-\phi^{-1},~~1)$\\
 		\cline{2-3}
 		& $F_{111}$ &$\frac{1}{2}(-\phi,~~0,~~-\phi^{-1},~~-1)$\\
 		\cline{2-3}
 		& $F_{112}$ &$\frac{1}{2}(-\phi,~~1,~~0,~~\phi^{-1})$\\
 		\cline{2-3}
 		8& $F_{113}$ &$\frac{1}{2}(-\phi,~~1,~~0,~~-\phi^{-1})$\\
 		\cline{2-3}
 		& $F_{114}$ &$\frac{1}{2}(-\phi,~~-1,~~0,~~\phi^{-1})$\\
 		\cline{2-3}
 		& $F_{115}$ &$\frac{1}{2}(-\phi,~~-1,~~0,~~-\phi^{-1})$\\
 		\cline{2-3}
 		& $F_{116}$ &$\frac{1}{2}(-\phi,~~\phi^{-1},~~1,~~0)$\\
 		\cline{2-3}
 		& $F_{117}$ &$\frac{1}{2}(-\phi,~~\phi^{-1},~~-1,~~0)$\\
 		\cline{2-3}
 		& $F_{118}$ &$\frac{1}{2}(-\phi,~~-\phi^{-1},~~1,~~0)$\\
 		\cline{2-3}
 		& $F_{119}$ &$\frac{1}{2}(-\phi,~~-\phi^{-1},~~-1,~~0)$\\
 		\hline
 		
 		9 & $F_{120}$ &$(-1,~~0,~~0,~~0)$  \\
 		\hline		
 	\end{tabular}
 }{\caption{}
 \label{table:coor3}
}
\end{floatrow}
\end{table}
}

\begin{figure}[H]
	\scalebox{0.33}[0.33]{\includegraphics {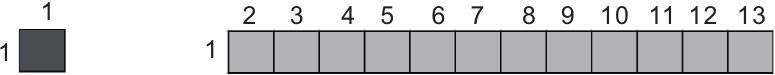}}
	\caption{ $X_{11} (=X_{99})$ and $M_{12}$}\label{figure:aj1}
\end{figure}

\begin{figure}[H]
	\scalebox{0.33}[0.33]{\includegraphics {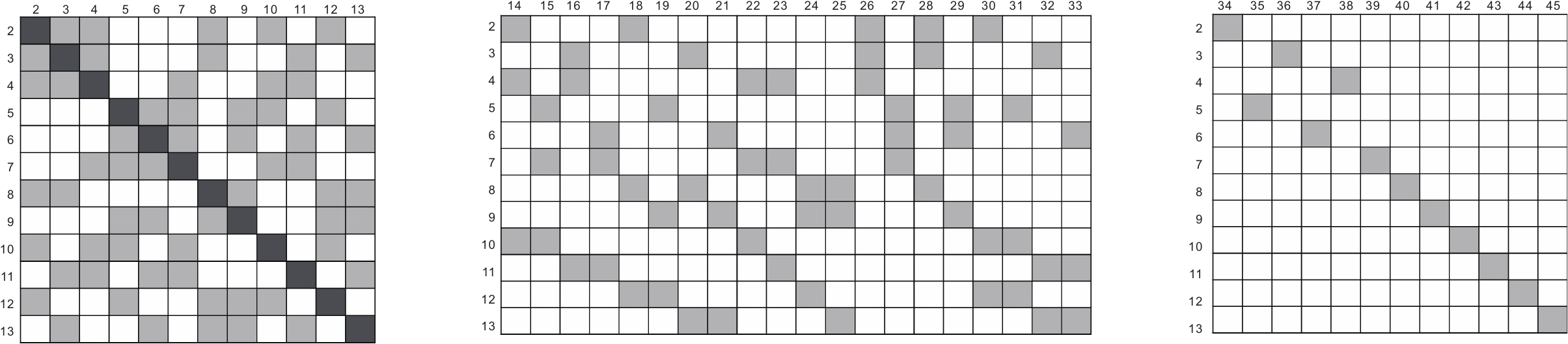}}
	\caption{$X_{22}$, $X_{23}$ and $X_{24}$}\label{figure:aj2}
\end{figure}

\begin{figure}[H]
	\scalebox{0.33}[0.33]{\includegraphics {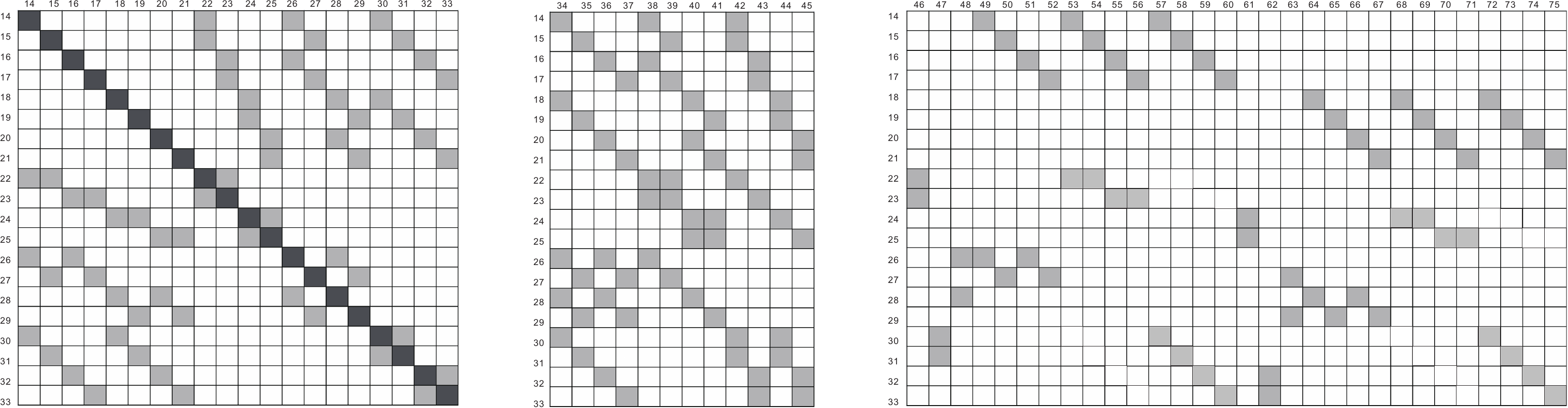}}
	\caption{$X_{33}$, $X_{34}$ and $X_{35}$}\label{figure:aj3}
\end{figure}

\begin{figure}[H]
	\scalebox{0.33}[0.33]{\includegraphics {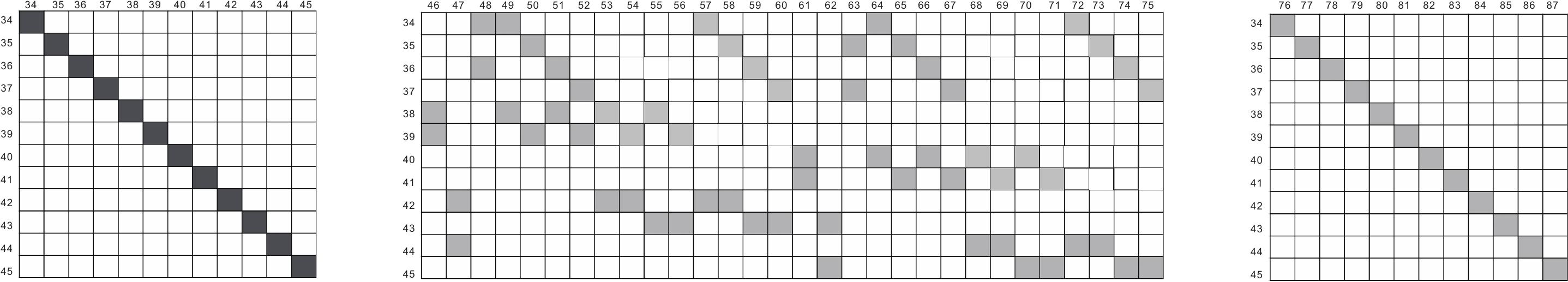}}
	\caption{$X_{44}$, $X_{45}$ and $X_{46}$}\label{figure:aj4}
\end{figure}

	\begin{figure}[H]
		\scalebox{0.33}[0.33]{\includegraphics {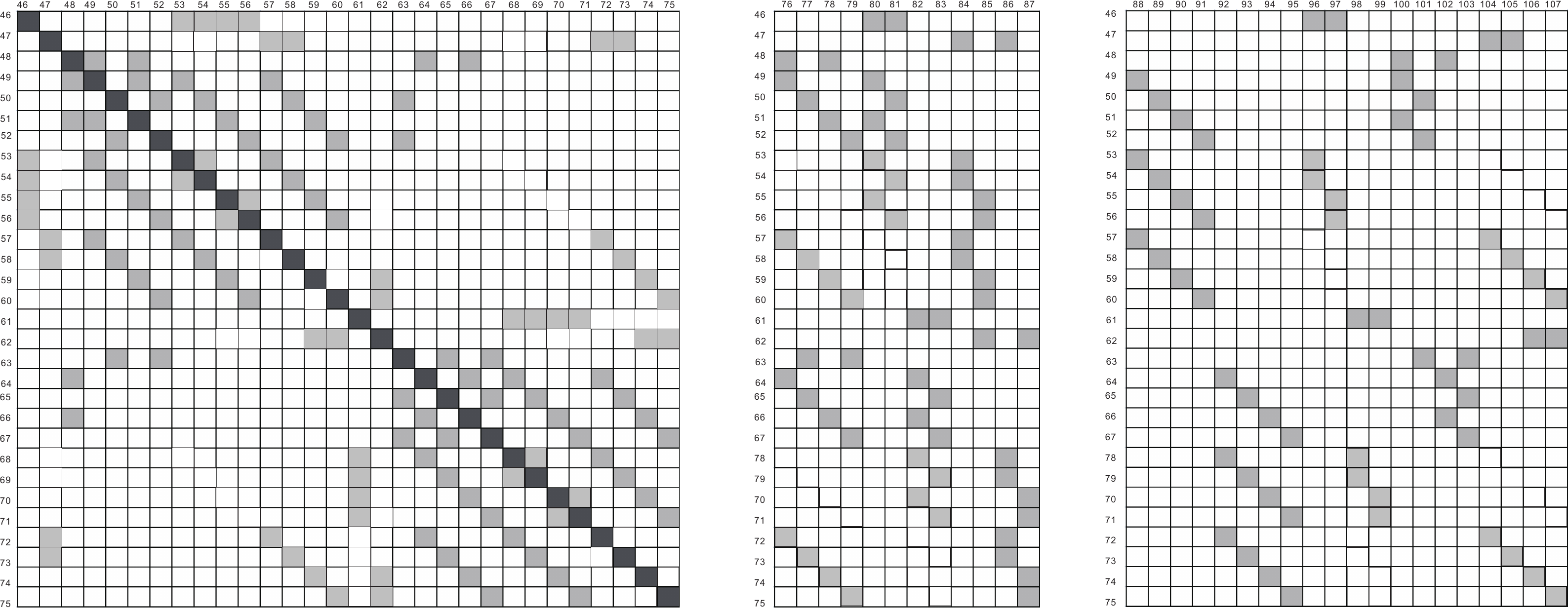}}
		\caption{$X_{55}$, $X_{56}$ and $X_{57}$}\label{figure:aj5}
	\end{figure}
	
	\begin{figure}[H]
		\scalebox{0.33}[0.33]{\includegraphics {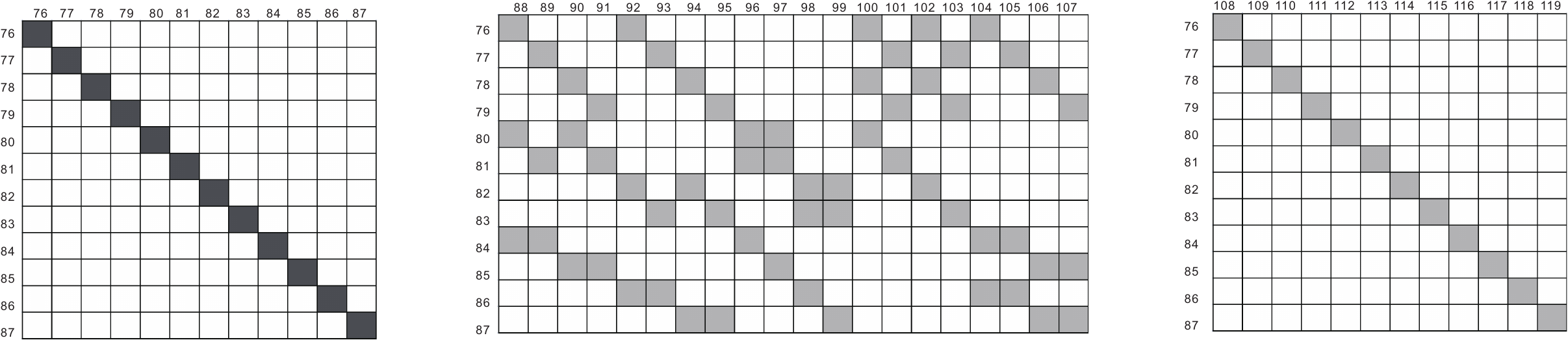}}
		\caption{$X_{66}$, $X_{67}$ and $X_{68}$}\label{figure:aj6}
	\end{figure}
	
	\begin{figure}[H]
		\scalebox{0.33}[0.33]{\includegraphics {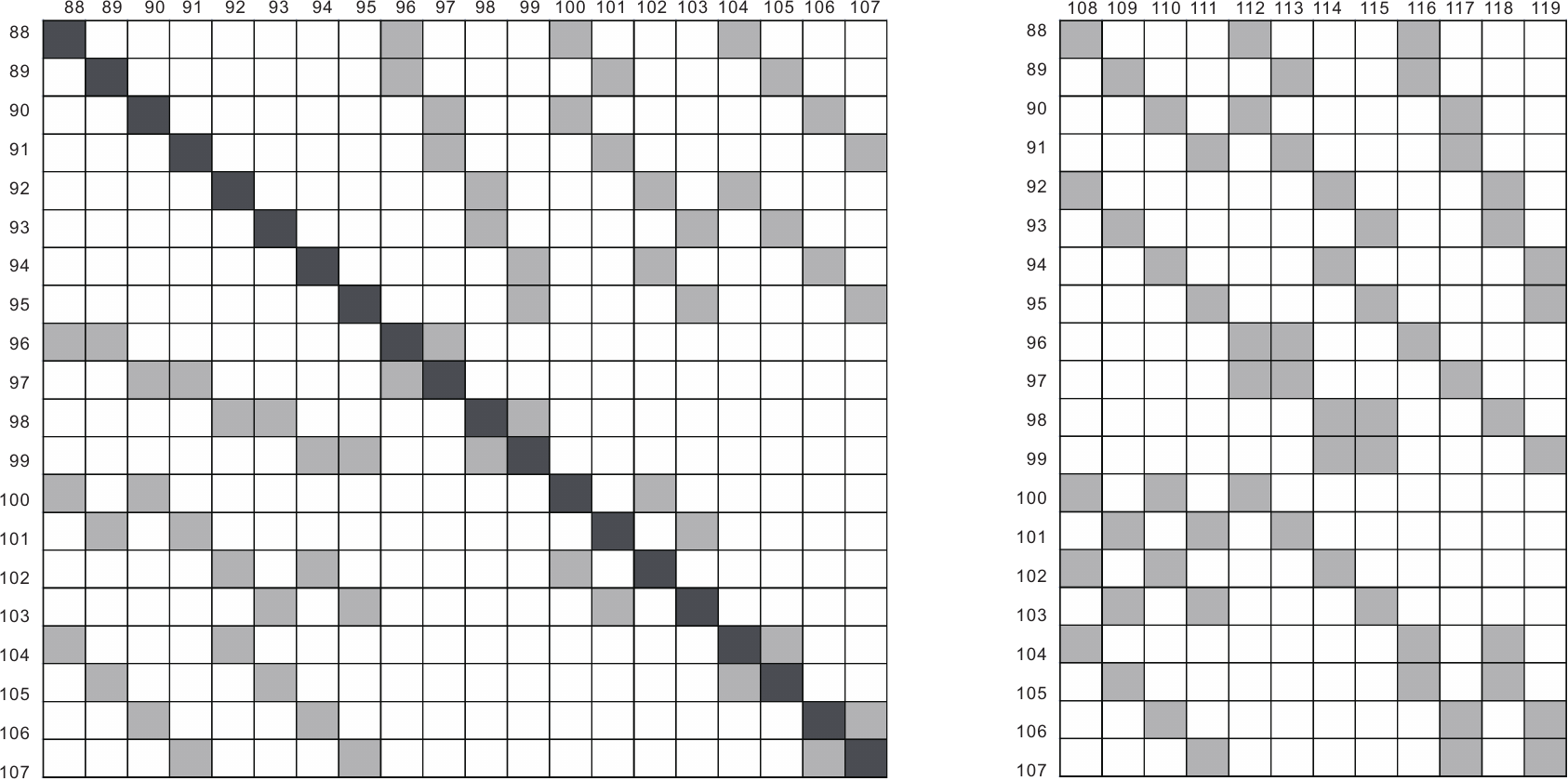}}
		\caption{$X_{77}$ and $X_{78}$}\label{figure:aj7}
	\end{figure}
	
	\begin{figure}[H]
		\scalebox{0.33}[0.33]{\includegraphics {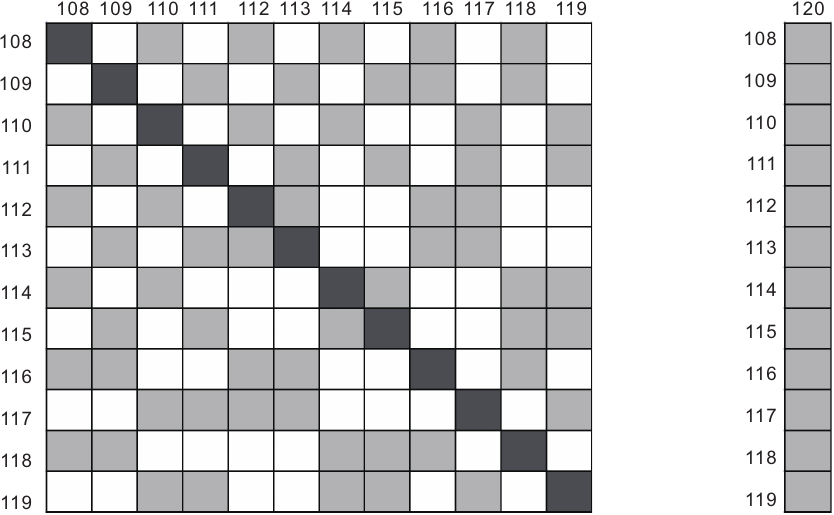}}
		\caption{$X_{88}$, $X_{89}$}\label{figure:aj8}
	\end{figure}

{\footnotesize	
	\renewcommand\baselinestretch{1.5}\selectfont
	\begin{longtable}{|c|c|p{11.3cm}|}
		
		\captionsetup{font=normalsize}
		\caption{Classification of k-coloring and representatives}\label{table:rep}\\
		\hline 
	$k$ & $\#$ representative &  representatives of added color\\
	\hline	
	5 & 3 & (3), (7), (15)\\
	\hline	
	6 & 7 & (3,5), (3,7), (3,12), (3,13), (3,15), (7,11), (7,15)\\
	\hline
	&  &(3,5,6), (3,5,7), (3,5,9), (3,5,10), (3,5,11), (3,5,14), \\
	7& 16 &(3,5,15), (3,7,11), (3,7,12), (3,7,13), (3,7,15)\\
	& & (3,12,15), (3,13,14), (3,13,15), (7,11,13), (7,11,15)\\
	\hline
	& & (3,5,6,7), (3,5,6,9), (3,5,6,11), (3,5,6,15), (3,5,7,9),\\
	& & (3,5,7,10), (3,5,7,11), (3,5,7,14), (3,5,7,15), (3,5,9,14), \\
	&  & (3,5,9,15), (3,5,10,12), (3,5,10,13), (3,5,10,15), (3,5,11,13), \\
	8&28 & (3,5,11,14), (3,5,11,15), (3,5,14,15), (3,7,11,12), (3,7,11,13),\\
	& & (3,7,11,15), (3,7,12,13), (3,7,12,15), (3,7,13,14), (3,7,13,15), \\
	& & (3,13,14,15), (7,11,13,14), (7,11,13,15)\\
	\hline
	& & (3,5,6,7,9), (3,5,6,7,11), (3,5,6,7,15), (3,5,6,9,10), (3,5,6,9,11),\\
	& & (3,5,6,9,14), (3,5,6,9,15), (3,5,6,11,13), (3,5,6,11,15), \\
	& & (3,5,7,9,11), (3,5,7,9,14), (3,5,7,9,15), (3,5,7,10,11), (3,5,7,10,12), \\
	&  &  (3,5,7,10,13), (3,5,7,10,14), (3,5,7,10,15), (3,5,7,11,13),  \\
	9& 35 &(3,5,7,11,14), (3,5,7,11,15), (3,5,7,14,15), (3,5,9,14,15),   \\
	& &(3,5,10,12,15), (3,5,10,13,14), (3,5,10,13,15), (3,5,11,13,14), \\
	& &(3,5,11,13,15), (3,5,11,14,15), (3,7,11,12,13), (3,7,11,12,15),\\
	&&(3,7,11,13,14), (3,7,11,13,15), (3,7,12,13,15), (3,7,13,14,15),\\
	&& (7,11,13,14,15)\\
	\hline
	& & (3,5,6,7,9,10), (3,5,6,7,9,11), (3,5,6,7,9,14), (3,5,6,7,9,15), \\
	& & (3,5,6,7,11,13), (3,5,6,7,11,15), (3,5,6,9,10,12), (3,5,6,9,10,13),\\
	& &(3,5,6,9,10,15), (3,5,6,9,11,13), (3,5,6,9,11,14), (3,5,6,9,11,15),  \\
	&  & 	(3,5,6,9,14,15), (3,5,6,11,13,14), (3,5,6,11,13,15), (3,5,7,9,11,13),  \\
	10& 35& 	(3,5,7,9,11,14), (3,5,7,9,11,15), (3,5,7,9,14,15), (3,5,7,10,11,12), \\
	& &(3,5,7,10,11,13), (3,5,7,10,11,15), (3,5,7,10,12,14),  \\
	& & (3,5,7,10,12,15), (3,5,7,10,13,14), (3,5,7,10,13,15),  \\
	& &(3,5,7,10,14,15), (3,5,7,11,13,14), (3,5,7,11,13,15),\\
	
	& & (3,5,7,11,14,15), (3,5,10,13,14,15), (3,5,11,13,14,15), \\
	& &(3,7,11,12,13,14), (3,7,11,12,13,15), (3,7,11,13,14,15)\\
	\hline
	&&(3,5,6,7,9,10,11), (3,5,6,7,9,10,12), (3,5,6,7,9,10,13), \\
	&&(3,5,6,7,9,10,15),  (3,5,6,7,9,11,13), (3,5,6,7,9,11,14),  \\
	&&(3,5,6,7,9,11,15), (3,5,6,7,9,14,15), (3,5,6,7,11,13,14),\\
	
	&& (3,5,6,7,11,13,15), (3,5,6,9,10,12,15), (3,5,6,9,10,13,14),\\
	11&28&(3,5,6,9,10,13,15), (3,5,6,9,11,13,14), (3,5,6,9,11,13,15),  \\
	&&(3,5,6,9,11,14,15), (3,5,6,11,13,14,15), (3,5,7,9,11,13,14), \\
	&&(3,5,7,9,11,13,15), (3,5,7,9,11,14,15), (3,5,7,10,11,12,13),\\
	
	&& (3,5,7,10,11,12,15), (3,5,7,10,11,13,14), (3,5,7,10,11,13,15),\\
	&&(3,5,7,10,12,14,15), (3,5,7,10,13,14,15), (3,5,7,11,13,14,15), \\
	\hline
	&&(3,7,11,12,13,14,15), (3,5,6,7,9,10,11,12),  \\
	&&(3,5,6,7,9,10,11,13), (3,5,6,7,9,10,12,15), (3,5,6,7,9,10,13,14),\\
	&& (3,5,6,7,9,10,13,15), (3,5,6,7,9,11,13,14), (3,5,6,7,9,11,13,15),  \\
	12&16&(3,5,6,7,9,11,14,15), (3,5,6,7,11,13,14,15), (3,5,6,9,10,13,14,15), \\
	&&(3,5,6,9,11,13,14,15), (3,5,7,9,11,13,14,15), (3,5,7,10,11,12,13,14), \\
	&&(3,5,7,10,11,12,13,15), (3,5,7,10,11,13,14,15)\\
	\hline
	
	&&(3,5,6,7,9,10,11,12,13), (3,5,6,7,9,10,11,12,15),\\ &&(3,5,6,7,9,10,11,13,14), (3,5,6,7,9,10,11,13,15),\\ 13&7&(3,5,6,7,9,10,13,14,15), (3,5,6,7,9,11,13,14,15), \\
	&&(3,5,7,10,11,12,13,14,15)\\
	\hline
	&&(3,5,6,7,9,10,11,12,13,14), (3,5,6,7,9,10,11,12,13,15),\\
	14&3&(3,5,6,7,9,10,11,13,14,15)\\
	\hline
	15&1&(3,5,6,7,9,10,11,12,13,14,15)\\
	\hline
\end{longtable}}

\end{appendices}


\end{document}